\documentclass[%
reprint,
amsmath,amssymb,
aps,
pre,
physrev,
]{revtex4-2}

\usepackage{graphicx}% Include figure files
\usepackage[graphicx]{realboxes}
\usepackage{dcolumn}% Align table columns on decimal point
\usepackage{bm}
\usepackage{lipsum}
\usepackage[justification=justified, singlelinecheck=off]{caption}
\usepackage{subcaption}
\usepackage{soul}
\usepackage{color}
\usepackage{rotating}

\DeclareMathOperator*{\argmax}{argmax}

\newtheorem{Conjecture}{Conjecture}

\usepackage{hyperref}% add hypertext capabilities
\hypersetup{
	colorlinks=true,
	linkcolor=black,
	filecolor=black,
	citecolor=black,      
	urlcolor=black,
}
\usepackage[mathlines]{lineno}% Enable numbering of text and display math
\begin{document}
	
	\preprint{APS/123-QED}
	
	\title{Densest plane group packings of regular polygons}
	
	\author{Miloslav Torda}
	\email{miloslav.torda@liverpool.ac.uk}
	\affiliation{
		Leverhulme Research Centre for Functional Materials Design, University of Liverpool, Liverpool, UK
	}
	\affiliation{
		Department of Computer Science, University of Liverpool, Liverpool, UK
	}
	\author{John Y. Goulermas}
	\email{j.y.goulermas@liverpool.ac.uk}
	\affiliation{
		Department of Computer Science, University of Liverpool, Liverpool, UK
	}
	\author{Vitaliy Kurlin}
	\email{vitaliy.kurlin@liverpool.ac.uk}
	\affiliation{
		Department of Computer Science, University of Liverpool, Liverpool, UK
	}
	\author{Graeme M. Day}
	\email{G.M.Day@soton.ac.uk}
	\affiliation{
		School of Chemistry, University of Southampton, Southampton, UK
	}
	
	%\date{\today}
	
	\begin{abstract}
		Packings of regular convex polygons ($n$-gons) that are sufficiently dense have been studied extensively in the context of modeling physical and biological systems as well as discrete and computational geometry. Former results were mainly regarding densest lattice or double-lattice configurations. Here we consider all two-dimensional crystallographic symmetry groups (plane groups) by restricting the configuration space of the general packing problem of congruent copies of a compact subset of the two-dimensional Euclidean space to particular isomorphism classes of the discrete group of isometries. We formulate the plane group packing problem as a nonlinear constrained optimization problem. By means of the Entropic Trust Region Packing Algorithm that approximately solves this problem, we examine some known and unknown densest packings of various 
		$n$-gons in all $17$ plane groups and state conjectures about common symmetries of the densest plane group packings for every $n$-gon.
	\end{abstract}
	
	\maketitle

	\section{Introduction}
	
	Understanding the packing properties of crystalline solids has important implications for solid state physics modeling \cite{torquato2018}, materials science \cite{boles2016,sanchez2022}, and biophysics \cite{liang2001}. In two-dimensional Euclidean space, crystal structures based on the densely packed representations of a molecule by a regular convex polygon ($n$-gon) were found to be adequate models for virus structures \cite{tarnai1995} or self-assembly of organic molecules on metal surfaces \cite{zoppi2012}. Compared to the densest packings, lower density but higher symmetry crystal structures of complex noncovalent molecular systems on surface substrates \cite{barth2007}, monolayer covalent organic frameworks \cite{cui2018}, or two-dimensional crystallization of proteins on lipid monolayers \cite{wang1999} can also be regarded as densest packings, although among a particular isomorphism class of periodic structures.
	
	Moreover, the crystallization problem and the disc packing problem are identical in the Euclidean space of dimensions two for some energy potentials \cite{theil2006,weinan2009}. Thus, fast ways of identifying dense packings could accelerate predictions of molecular crystal structures, where the usual approach is to search for the lowest energy configurations \cite{woodley2020}.
	
	The packing problem is well-studied in discrete and computational geometry. Substantial consideration has been given to densest lattice \cite{courant1965least,toth2013} and double-lattice packings \cite{kuperbergSup,hales} as special cases of the general densest packing configurations. On the other hand, little is known about densest packings when configurations are restricted to the remaining isomorphism classes of discrete groups of isometries of the two-dimensional Euclidean space. 
	
	\begin{figure*}[t]
		\centering
		\includegraphics[trim={0 0 0 0},width=\textwidth]{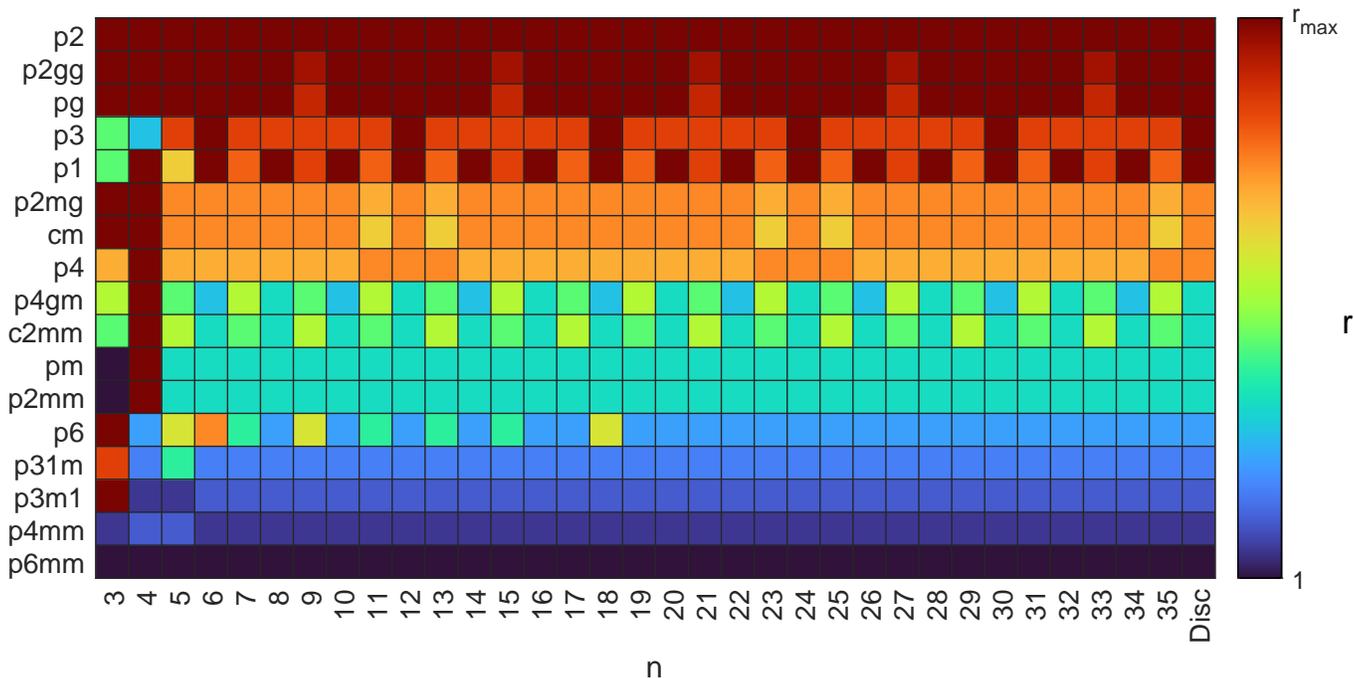}
		\caption{\label{fig:tabelFig6} The colored rank table of plane groups in relation to the number of vertices $n$ of an $n$-gon. For every $n=3,\ldots, 35$ plane groups are ordered according to densities in Table~\ref{tab:table1}, and a color is assigned based on rank $r$ ranging from one to $r_{max}$. The value of $r_{max}$ depends on a specific $n$-gon.}
	\end{figure*}
	
\begin{figure*}
		\centering
		
		\begin{minipage}{3.4cm}
			%\vspace{-6cm}
			\centering
			\includegraphics[trim={0.1 0.1 0.1 0.1},clip,height=2.3cm]{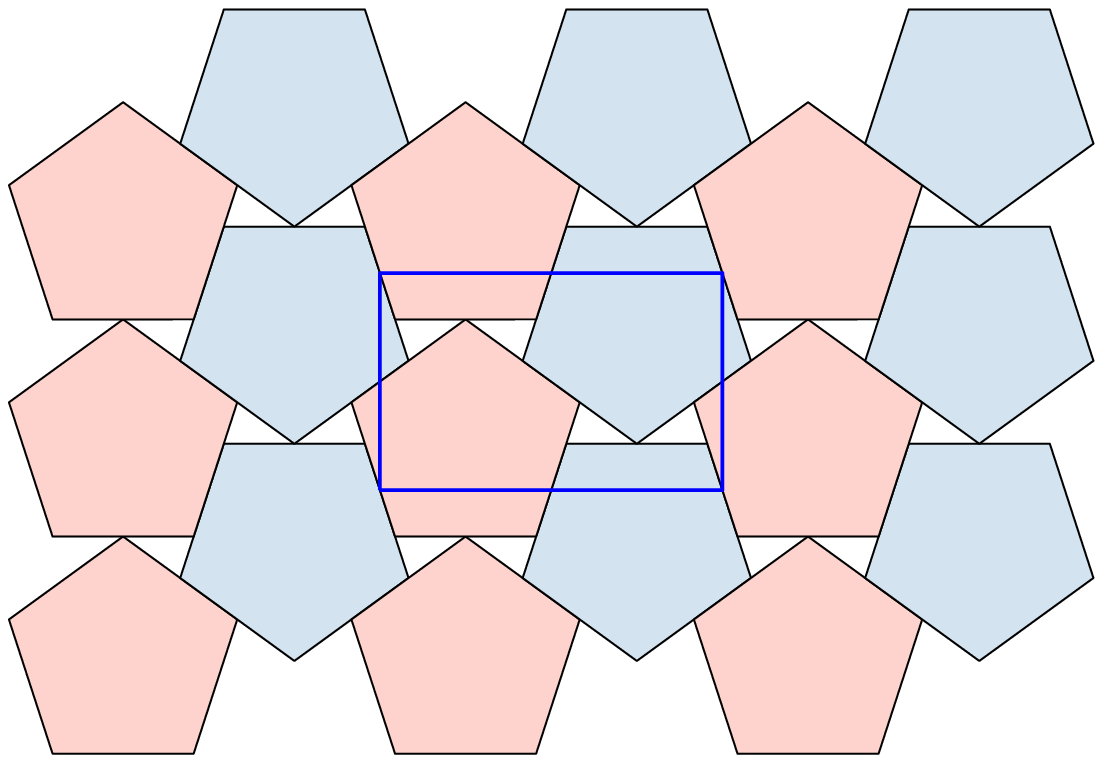}
			
			\includegraphics[trim={0.1 0.1 0.1 0.1},clip,height=2.3cm]{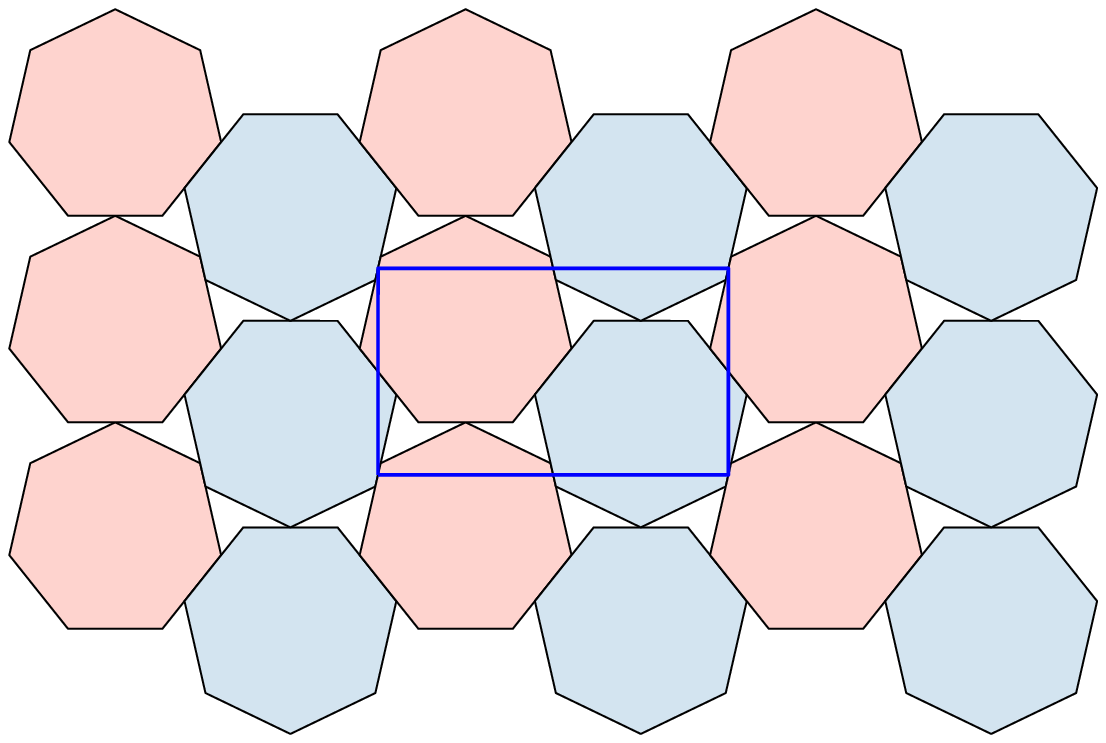}
			
			\includegraphics[trim={0.1 0.1 0.1 0.1},clip,height=2.3cm]{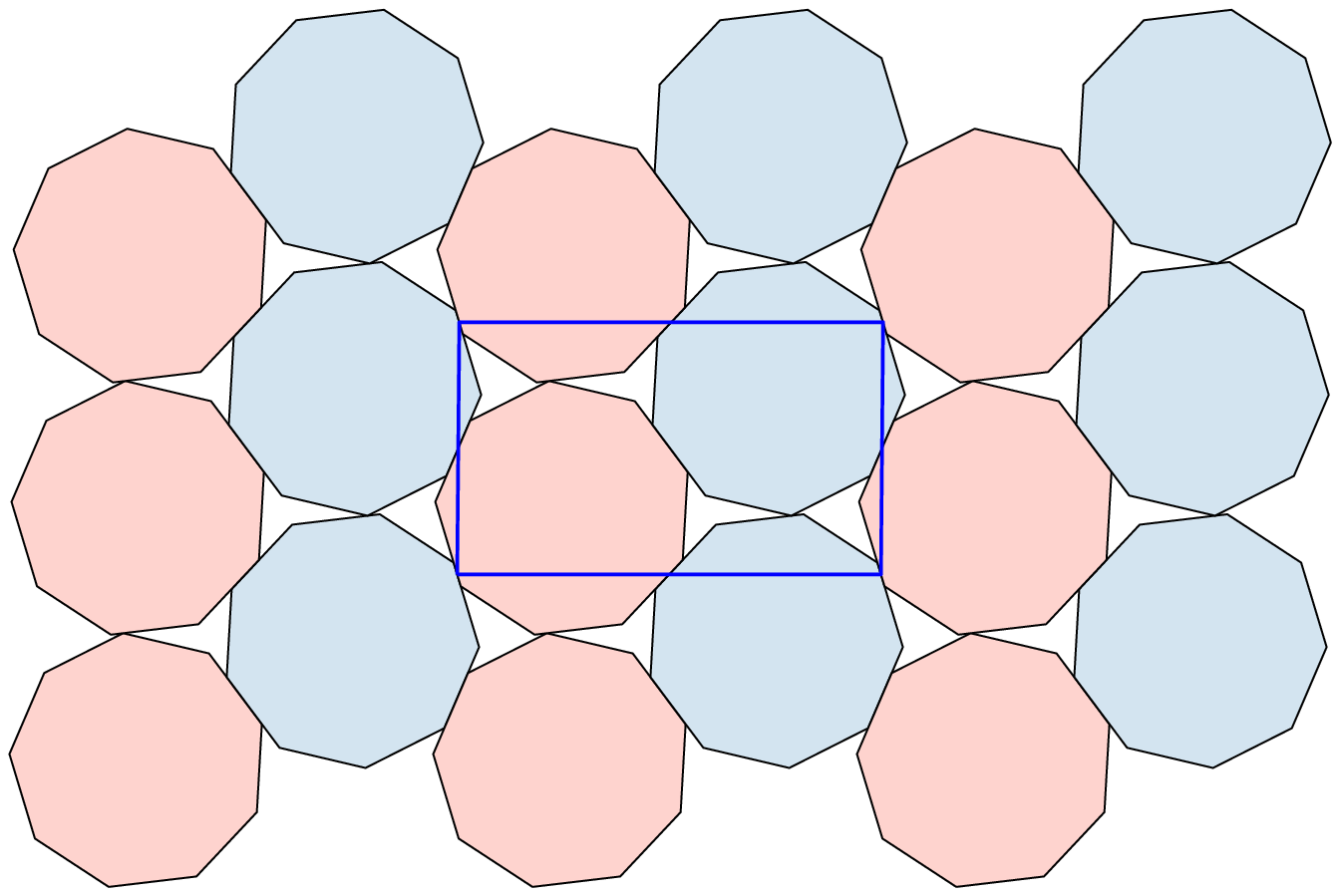}
			%\vspace{-0.5cm}
			
			\includegraphics[trim={0.1 0.1 0.1 0.1},clip,height=2.3cm]{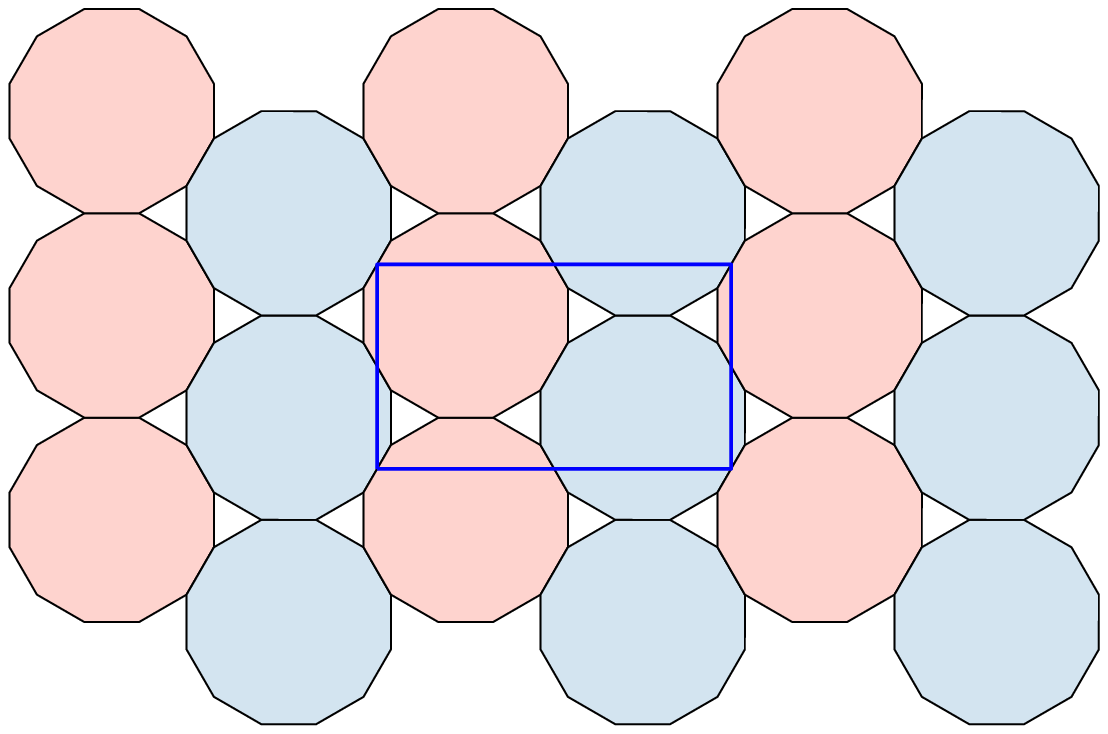}
			\subcaption{p2}
			%\vspace{-1.4cm} 
		\end{minipage}	
		\begin{minipage}{3.4cm}
			%\vspace{-6cm}
			\centering
			\includegraphics[trim={0.1 0.1 0.1 0.1},clip,height=2.3cm]{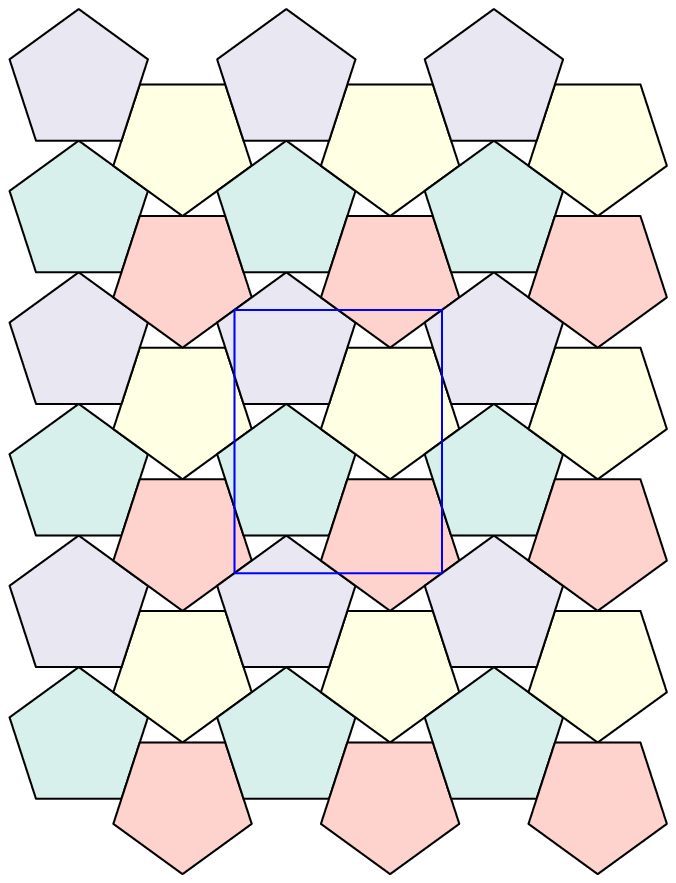}
			
			\includegraphics[trim={0.1 0.1 0.1 0.1},clip,height=2.3cm]{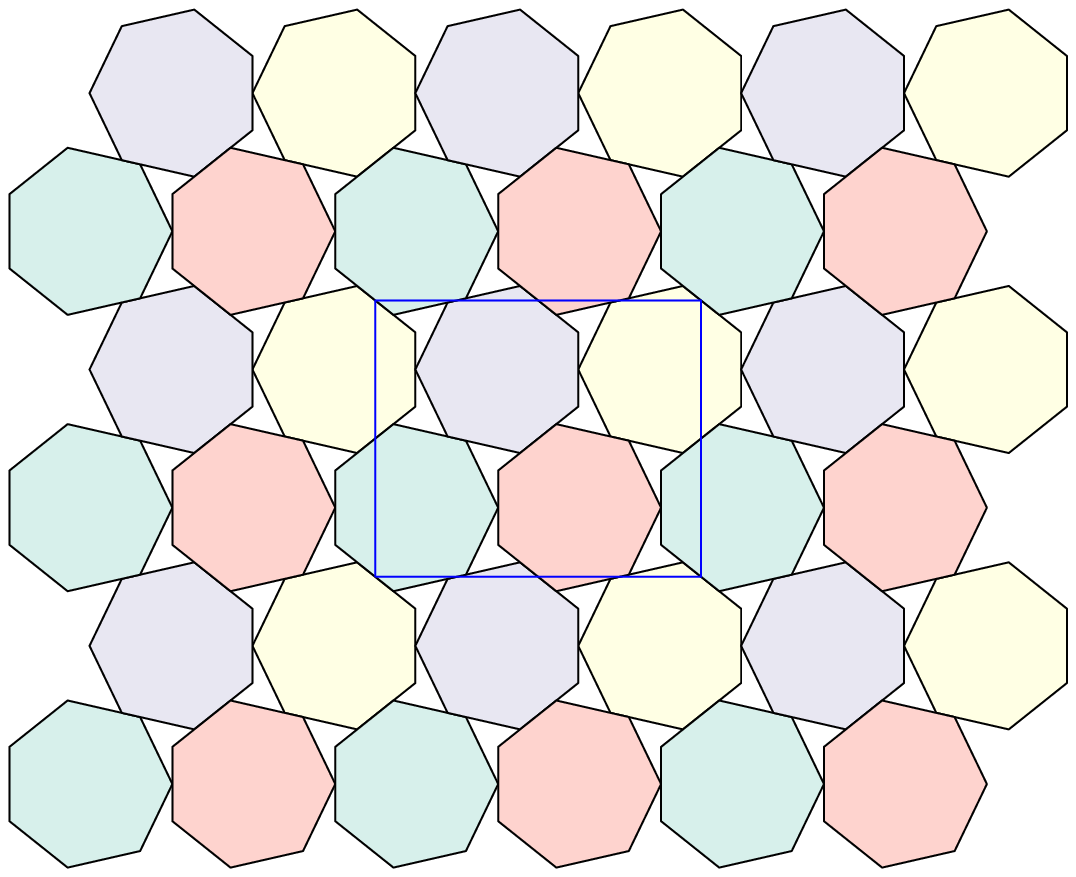}
			
			\includegraphics[trim={0.1 0.1 0.1 0.1},clip,height=2.3cm]{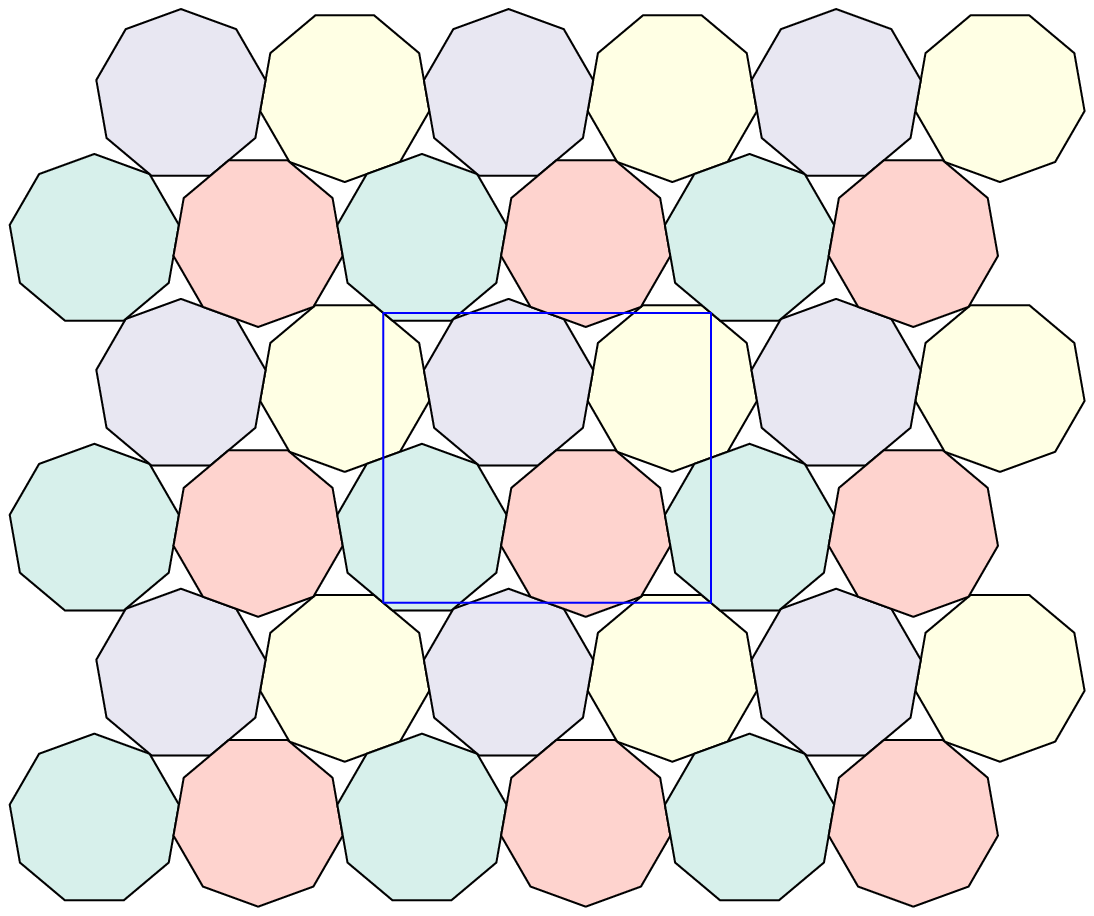}
			
			\includegraphics[trim={0.1 0.1 0.1 0.1},clip,height=2.3cm]{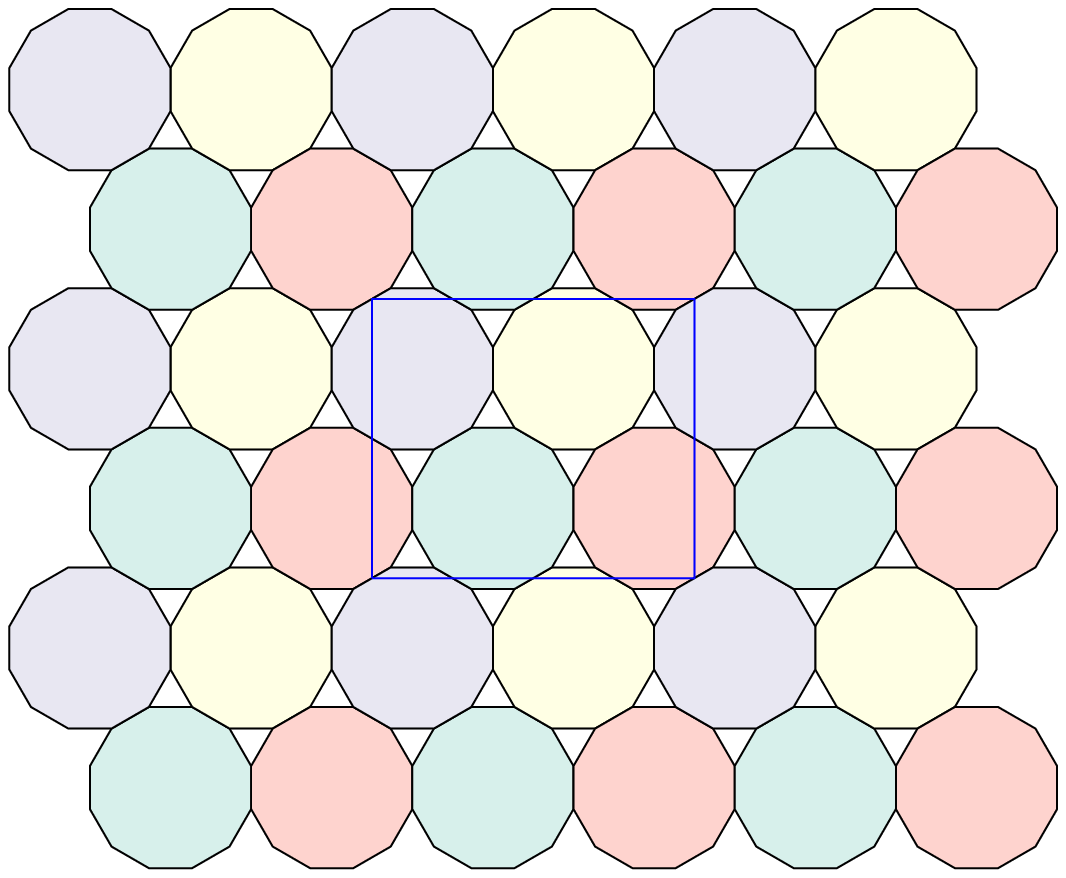}
			
			\subcaption{p2gg}
			
		\end{minipage}
		\begin{minipage}{3.4cm}
			%\vspace{-6cm}
			\centering
			\includegraphics[trim={0.1 0.1 0.1 0.1},clip,height=2.3cm]{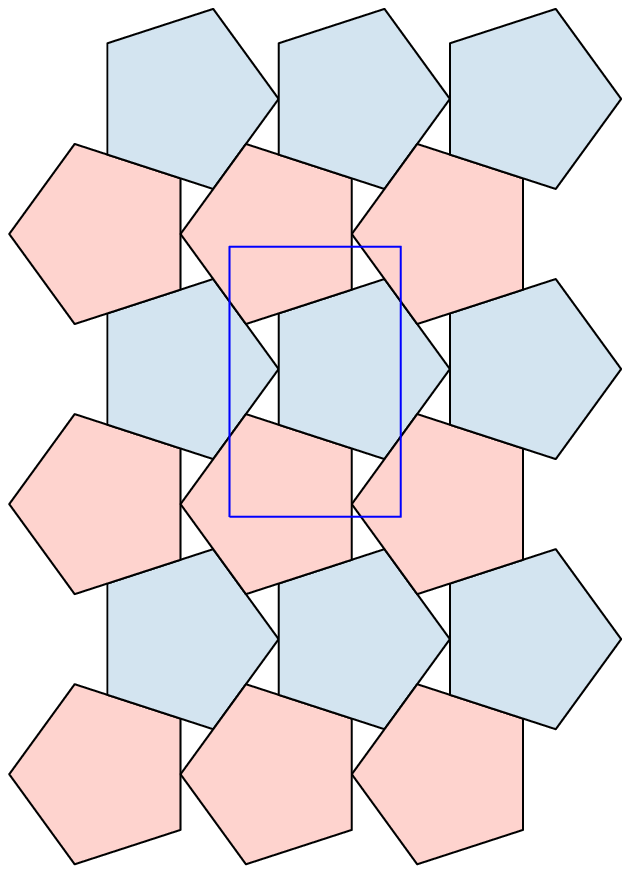}
			
			\includegraphics[trim={0.1 0.1 0.1 0.1},clip,height=2.3cm]{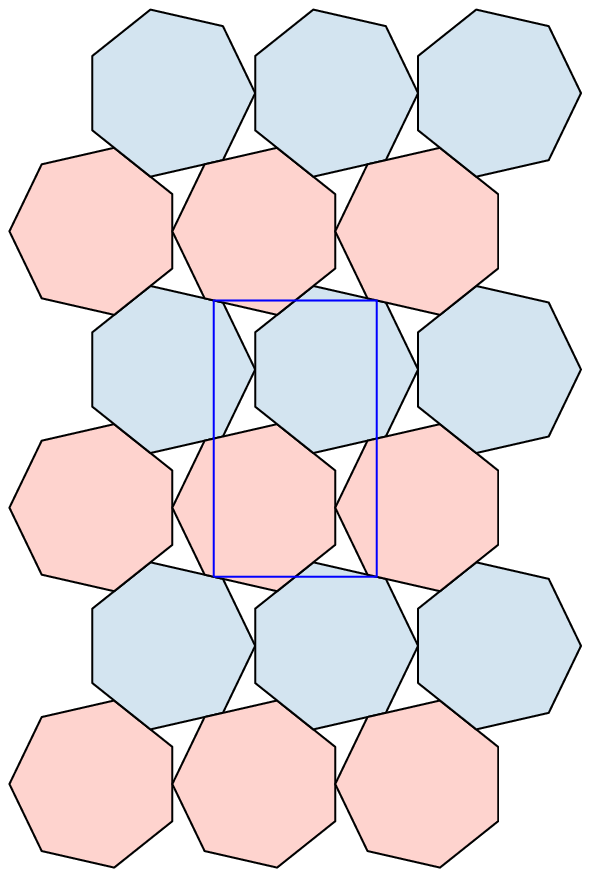}
			
			\includegraphics[trim={0.1 0.1 0.1 0.1},clip,height=2.3cm]{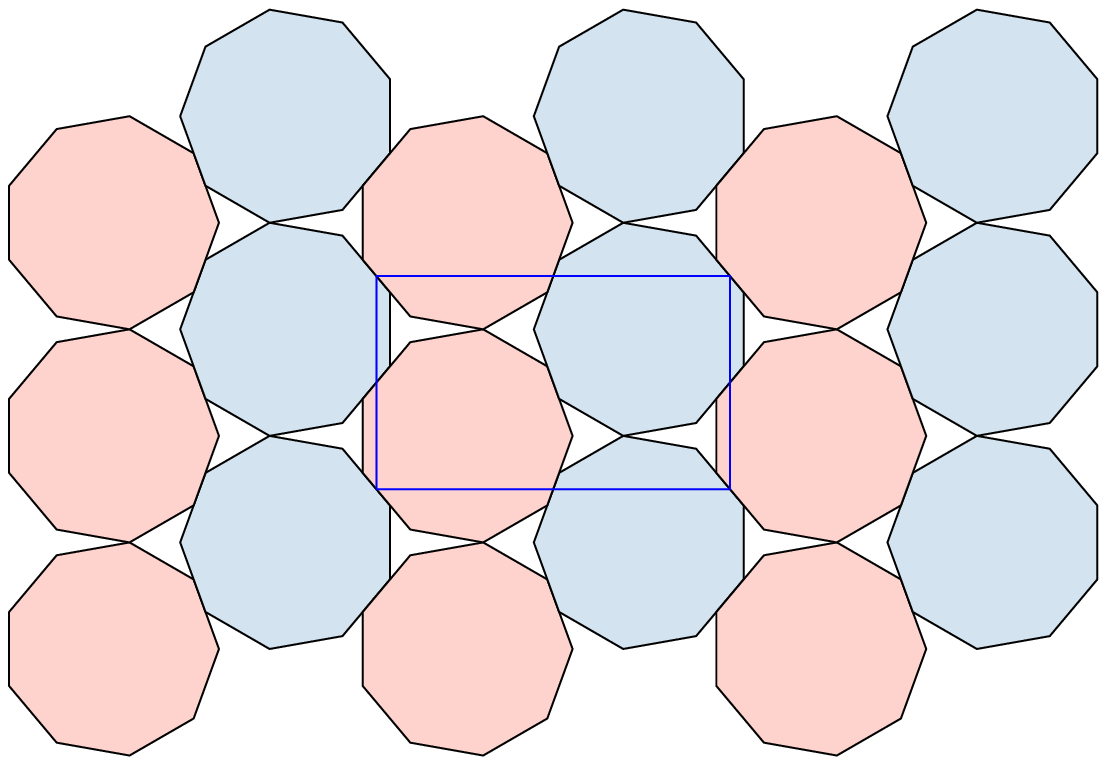}	
			
			\includegraphics[trim={0.1 0.1 0.1 0.1},clip,height=2.3cm]{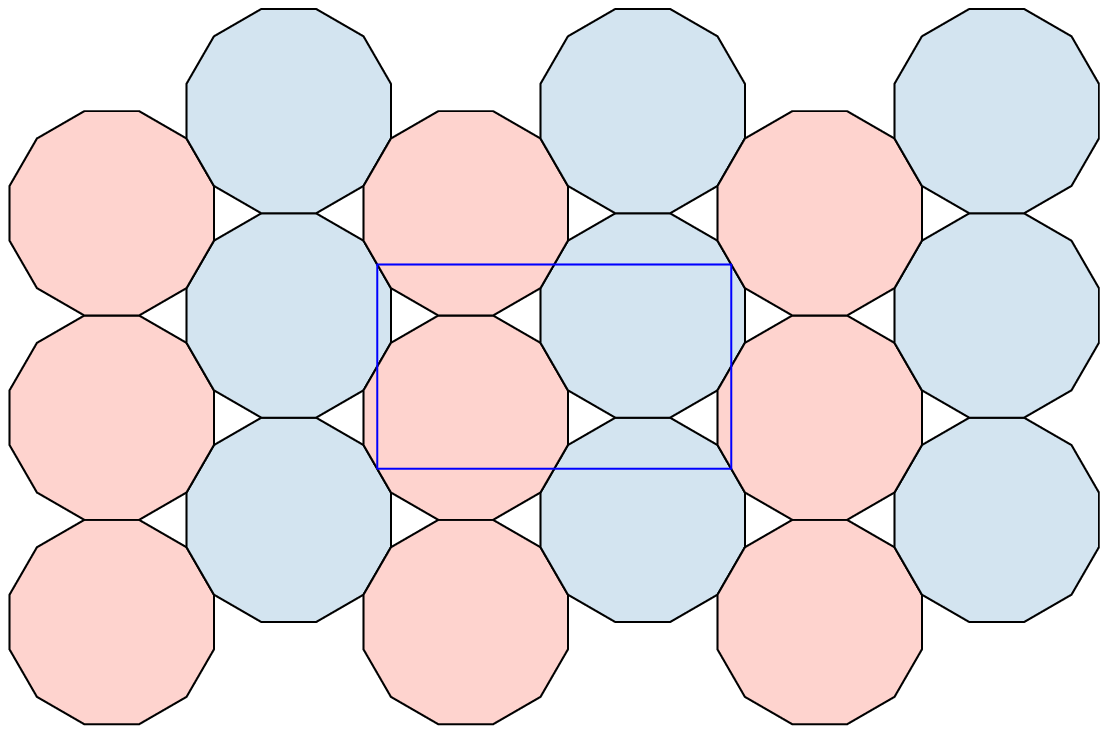}
			\subcaption{pg} 
		\end{minipage}
		\begin{minipage}{3.4cm}
			%\vspace{-6cm}
			\centering
			\includegraphics[trim={0.1 0.1 0.1 0.1},clip,height=2.3cm]{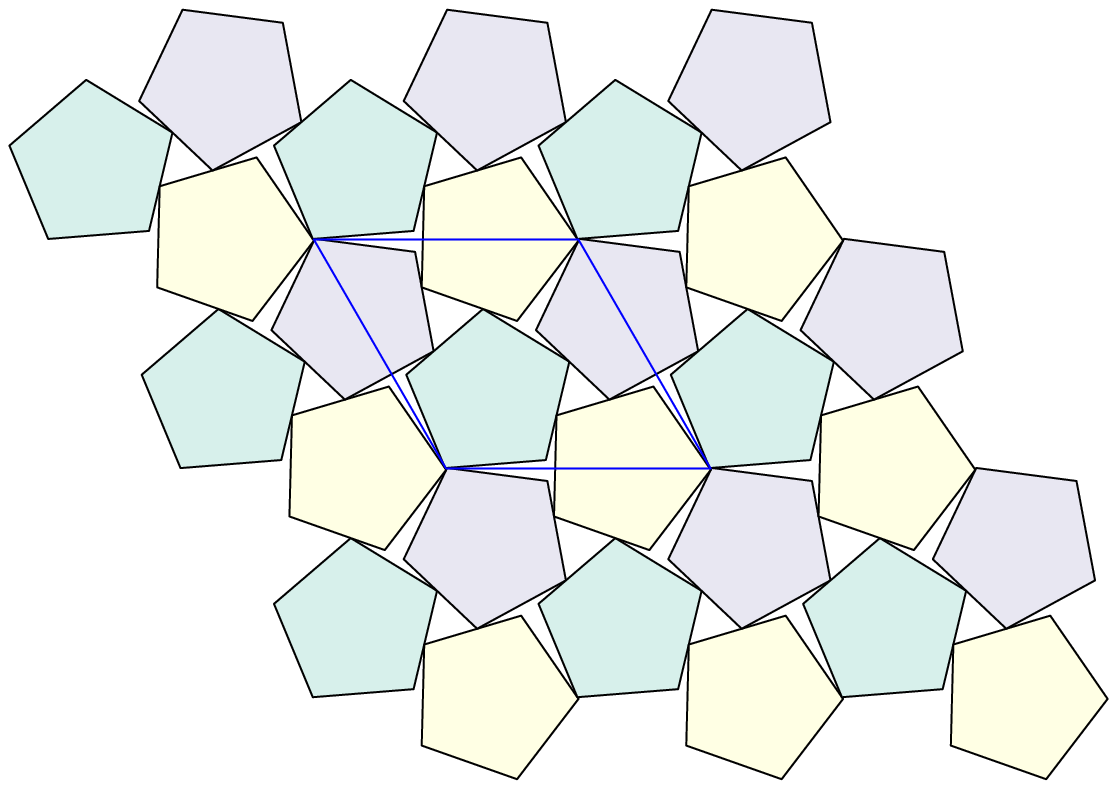}
			
			\includegraphics[trim={0.1 0.1 0.1 0.1},clip,height=2.3cm]{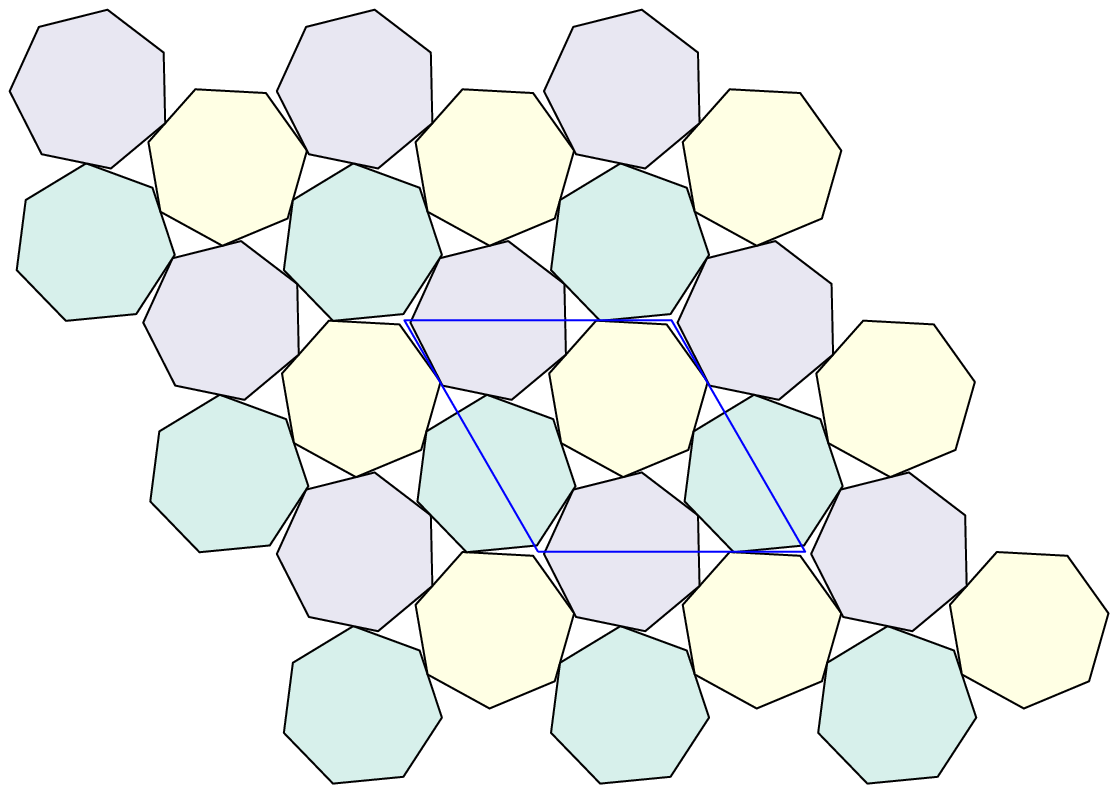}
			
			\includegraphics[trim={0.1 0.1 0.1 0.1},clip,height=2.3cm]{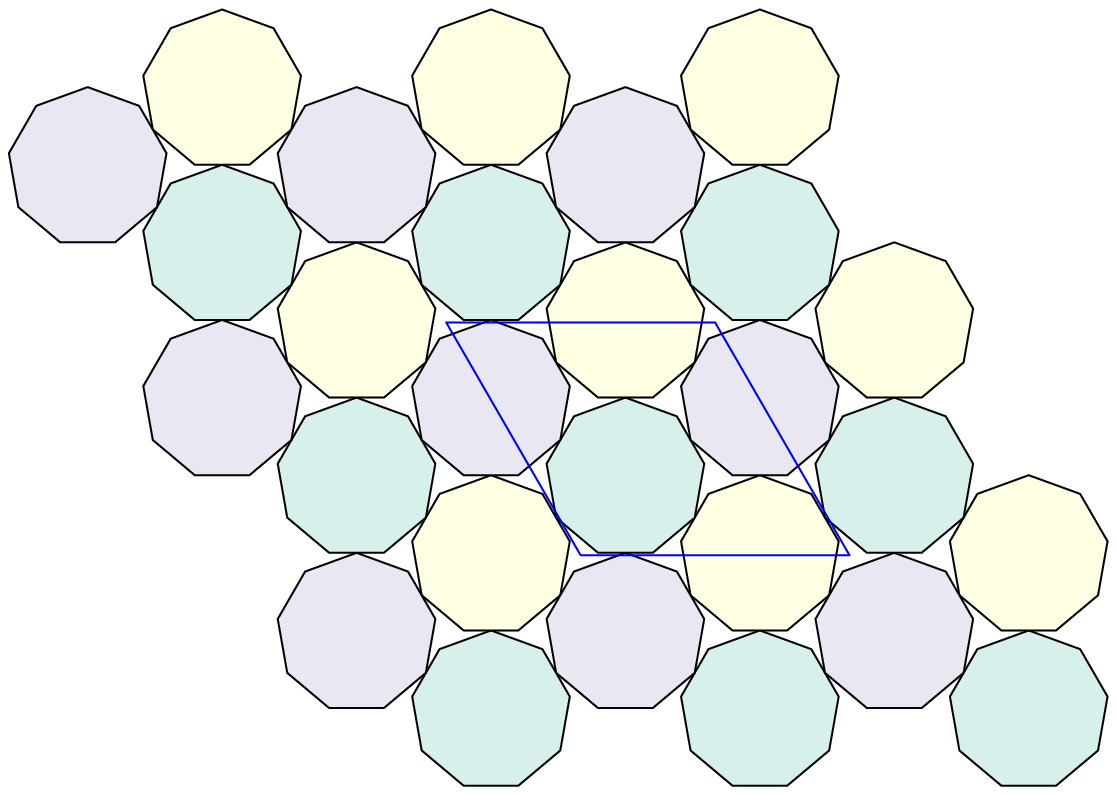}
			
			\includegraphics[trim={0.1 0.1 0.1 0.1},clip,height=2.3cm]{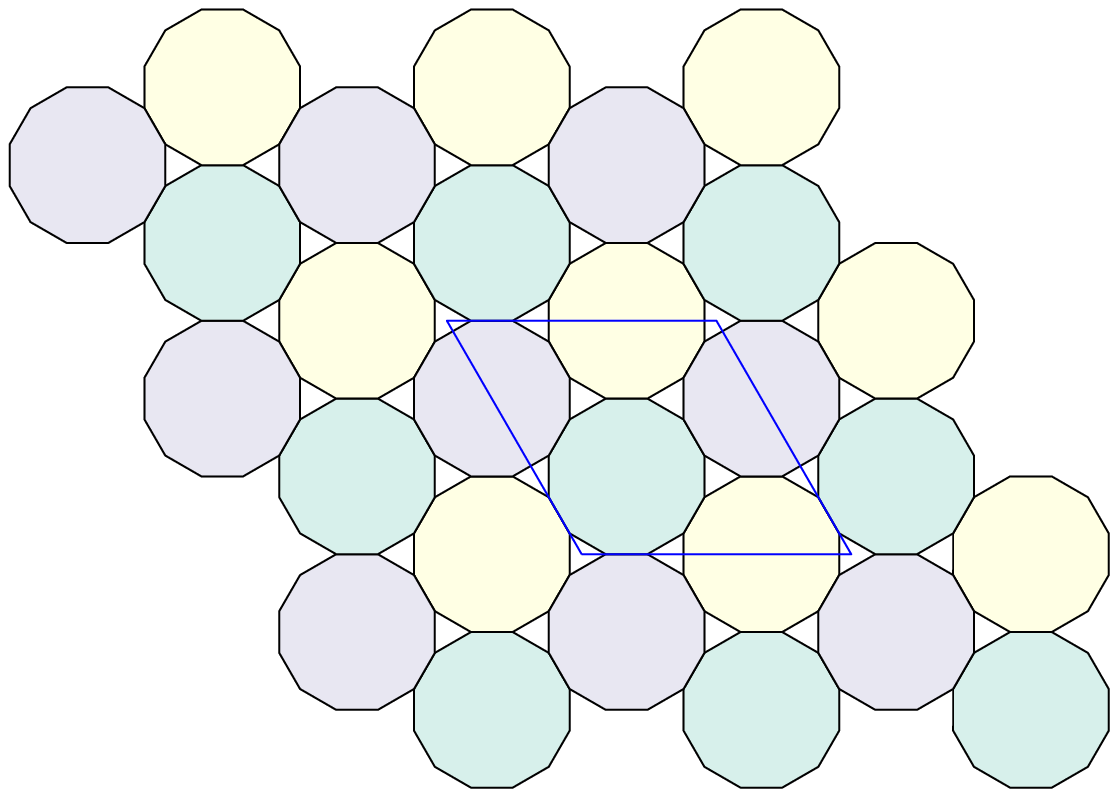}
			\subcaption{p3} 
		\end{minipage}
		\begin{minipage}{3.4cm}
			%\vspace{-6cm}
			\centering
			\includegraphics[trim={0.1 0.1 0.1 0.1},clip,height=2.3cm]{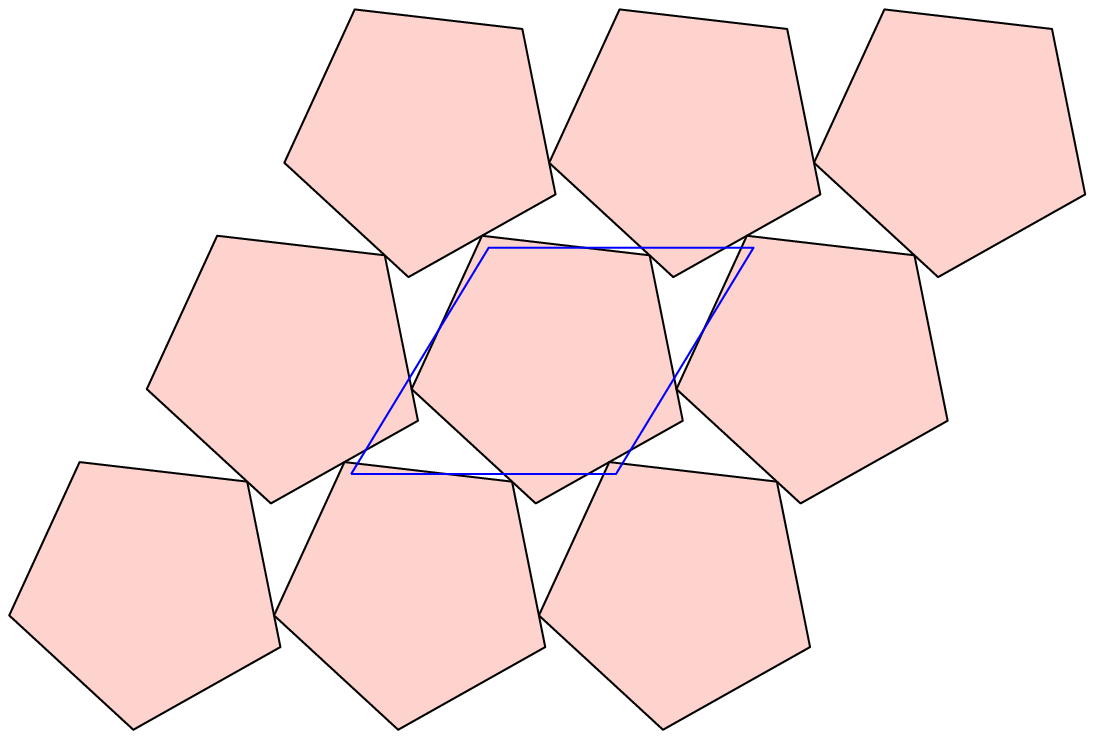}
			
			\includegraphics[trim={0.1 0.1 0.1 0.1},clip,height=2.3cm]{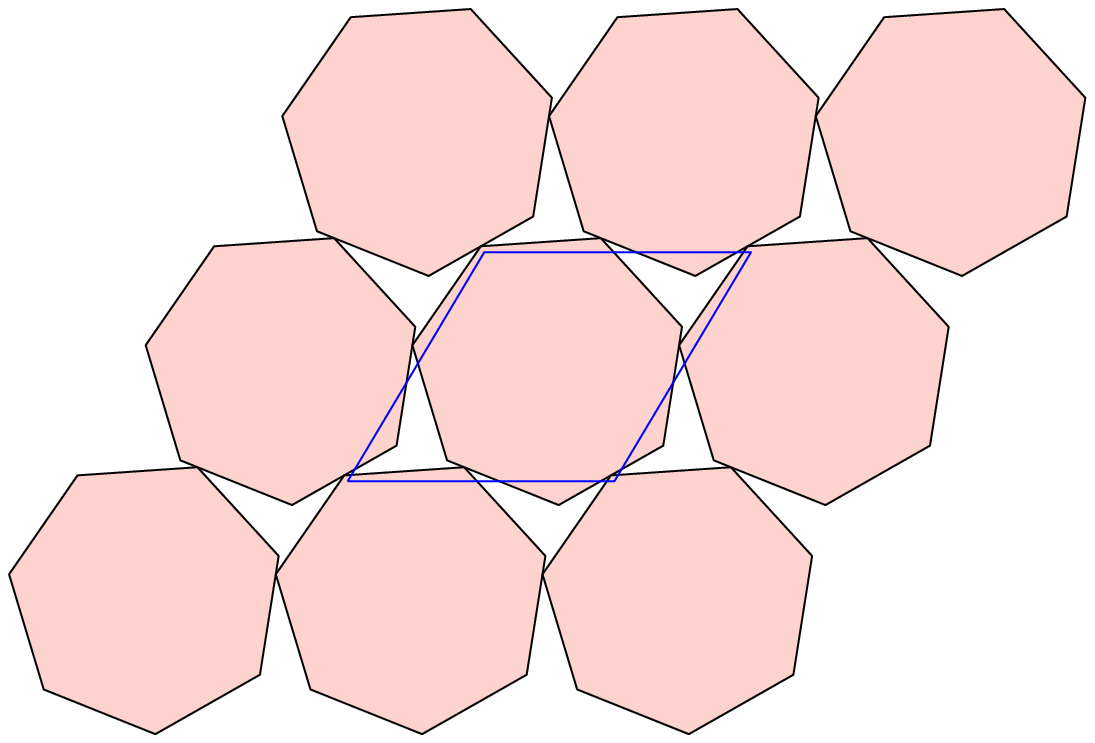}
			
			\includegraphics[trim={0.1 0.1 0.1 0.1},clip,height=2.3cm]{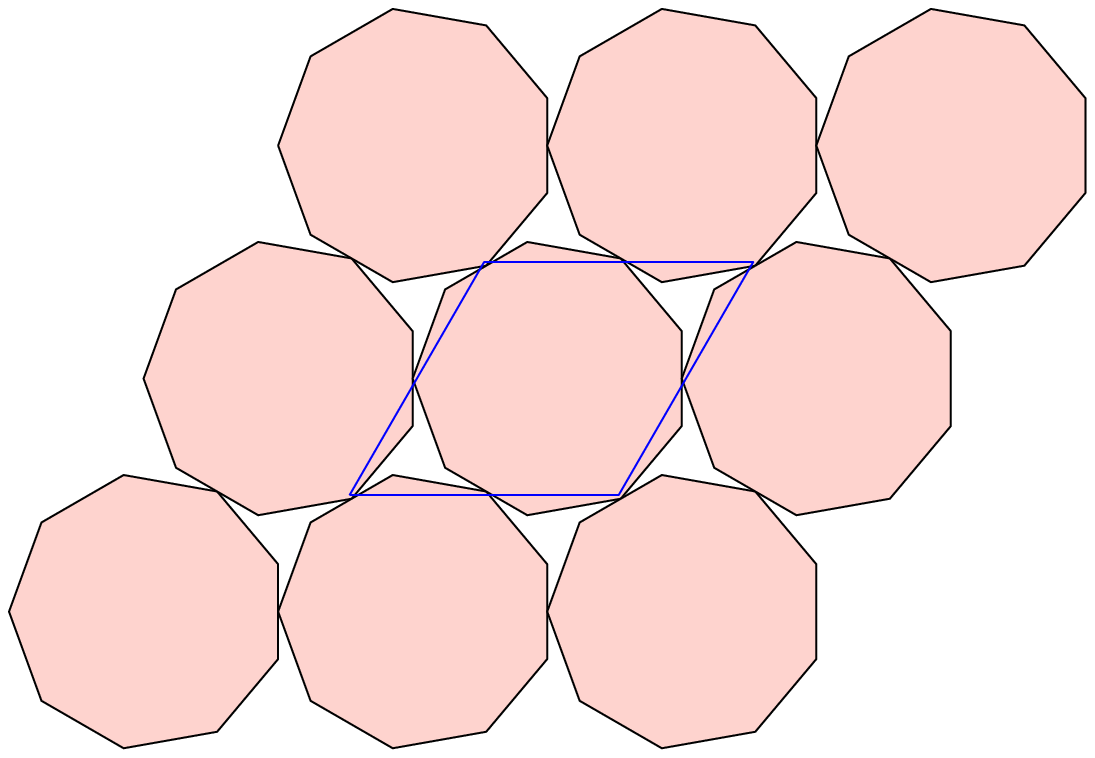}
			
			\includegraphics[trim={0.1 0.1 0.1 0.1},clip,height=2.3cm]{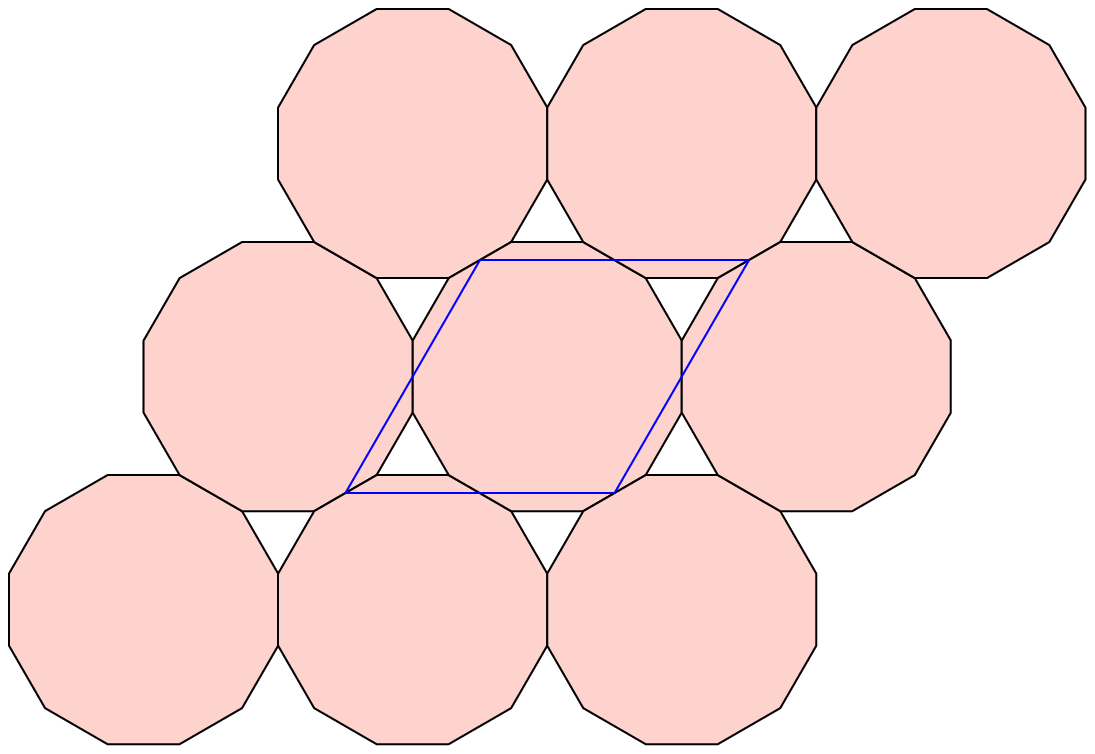}
			\subcaption{p1} 
		\end{minipage}
		
		\caption{\label{fig:heptagonFig1} 
			Densest configurations of (from top to bottom) pentagon, heptagon, enneagon, and dodecagon in plane groups $p2$, $p2gg$, $pg$, $p3$, and $p1$ with the following densities: pentagon in $p2 / p2gg / pg \approxeq 0.92131$, $p3 \approxeq 0.87048$ and $p1 \approxeq 0.81725$; heptagon in $p2 / p2gg / pg \approxeq 0.89269$, $p3 \approxeq 0.88085$ and $p1 \approxeq 0.86019$; enneagon in $p2 \approxeq 0.90103$, $p2gg \approxeq 0.89989$, $pg \approxeq 0.89860$ and $p3 / p1 \approxeq 0.88773$; dodecagon in $p2 / p2gg /pg / p3 /p1 \approxeq 0.92820$. The blue parallelogram denotes the primitive cell of the respective configuration. Colors represent symmetry operations modulo lattice translations.}
		%\vspace{-11pt}
	\end{figure*}

	We examine the densest packings of various $n$-gons where the packing configurations are restricted to one of the $17$ isomorphism classes of the discrete group of isometries of the two-dimensional Euclidean space containing a lattice subgroup, in literature also referred to as plane groups or wallpaper groups. FIG.~\ref{fig:heptagonFig1} illustrates plane group packings on densest $p2$, $p2gg$, $pg$, $p3$ and $p1$ configurations of a pentagon, heptagon, enneagon, and dodecagon. We consider finding the densest plane group packing as a nonlinear, constrained, and bounded optimization problem. Using the Entropic Trust Region Packing Algorithm \cite{torda2022entropic}, developed specifically to search for densest crystallographic symmetry group packings of arbitrary dimensions, we successfully recover approximations of known densest lattice and double-lattice packing configurations including a disc, regarded as a limiting $n$-gon when the number of vertices approaches infinity. Additionally, we obtained the previously unknown highest density packings of the $n$-gons for all $17$ plane groups and $n$
	equal to $3, 4, \ldots, 25, 30, 35, 36, 37, 39, 42, 55, 89$.
	
	Our experiments suggest the following relationships between symmetries of $n$-gons and shared symmetries of their respective densest plane group packing configurations divided into three classes. In the $p2/p2gg/pg/p3/p1$ plane group class,
	\begin{enumerate}
		\item except for centrally nonsymmetric $n$-gons containing a three-fold rotational symmetry with number of vertices higher or equal to nine, densities of the densest $p2$, $p2gg$ and $pg$ configurations are equal,
		\item for centrally symmetric $n$-gons, densities of the densest $p2$, $p2gg$, $pg$ and $p1$ configurations are equal,
		\item for centrally nonsymmetric $n$-gons containing a three-fold rotational symmetry, densities of the densest $p3$ and $p1$ configurations are equal,
		\item for centrally symmetric $n$-gons containing a three-fold rotational symmetry, densities of the densest $p2$, $p2gg$, $pg$, $p3$ and $p1$ configurations are equal,
	\end{enumerate}
	in the $p2mg/cm/p4$ plane group class,
	\begin{enumerate}
		\item except for $n$-gons with vertices equal to $12k-1$ and $12k+1$, where $k$ is an integer, densities of densest$p2mg$ and $cm$configurations are equal,
		\item for all $n$-gons with $12$-fold rotational symmetry, densities of the densest $p2mg$, $cm$ and $p4$ configurations are equal,
	\end{enumerate}
	and in the $p4gm/c2mm/pm/p2mm$ plane group class,
	\begin{enumerate}
		\item densities of the densest $pm$ and $p2mm$ configurations are equal for all $n$-gons,
		\item for centrally symmetric $n$-gons, densities of the densest $pm$, $p2mm$, and $c2mm$ configurations are equal,
		\item for $n$-gons containing a four-fold rotational symmetry, densities of the densest $pm$, $p2mm$, $c2mm$ and $p4gm$ configurations are equal,
	\end{enumerate}	
	Consequently, the densest known packings of a pentagon and a heptagon have higher symmetries than that of a double-lattice configuration.

	FIG.~\ref{fig:tabelFig6} visually summarizes our results. For each $n$-gon, all groups are ranked according to the density of the respective densest plane group packing, either obtained experimentally or extrapolated from the rankings of $n$-gons with similar symmetries. For instance, the ranking of densities of a $33$-gon is based on rankings of $9$-gon, $15$-gon, $21$-gon and $39$-gon.
	
	The manuscript is organized in the following way. Section~\ref{subsec:spacegrouppacking} introduces the plane group packing and the underlying densest plane group packing problem. In Section~\ref{sec:Reuslts}, we present the densest plane group packings of regular $n$-gons. We examine the symmetries of packing configurations in distinct plane group classes based on symmetries of densest plane group packings of a disc. Section~\ref{sec:conclusions} summarizes our experimental results in the form of multiple conjectures about common symmetries of the densest plane group packings of $n$-gons for arbitrary $n$.
	
	\section{Plane group packing}
	\label{subsec:spacegrouppacking}
	
	We consider the two-dimensional Crystallographic Symmetry Group (CSG) $G$, which is a discrete subgroup of the group of isometries of the two-dimensional Euclidean space containing a lattice subgroup. The parallelogram spanned by the generators of the lattice $L$
	associated with CSG $G$ is called the primitive cell and is denoted by $U_L$. An asymmetric unit is a subset of the primitive cell such that the whole two-dimensional Euclidean space is filled when the CSG symmetry operations are applied.
	
	It has to be noted that the term CSG has two distinct meanings in literature; one referring to an actual group and the other to a group isomorphism class. We refer to isomorphism classes of two-dimensional CSGs as plane groups. All two-dimensional CSGs are classified into $17$ plane groups. Each class is assigned to one of the four maximal crystallographic point group conjugacy classes, referred to as the crystal system. In all the following, we use the International Union of Crystallography plane group notation \cite{brock2016}.
	
	Given a two-dimensional CSG $G$, an element of a plane group $G$, and a polygon $K$ whose centroid lies in the asymmetric unit of $G$, by a CSG packing $\mathcal{K}_{G}$, we mean a collection of nonoverlapping orbits of $K$ with respect to the action of $G$ on the Euclidean plane.
	
	Since a CSG packing $\mathcal{K}_{G}$ is a periodic system of rotated and translated copies of polygon $K$, following the formula for the density of packing of a periodic system \cite{rogers}, the density of a CSG packing has a simple closed form expression
	\begin{equation*}\label{eq:density}
		\rho \left(\mathcal{K}_G \right) = \frac{N \text{area}  (K)}{\text{area}(U_{L})},
	\end{equation*}
	where $N$ is the number of symmetry operations in CSG $G$
	modulo translations by the lattice associated with $G$, $U_{L}$ is the primitive cell, and $\text{area}(\cdot)$ denotes area.
	
	Given a plane group $G$ and a polygon $K$, the densest plane group packing of $K$ is a CSG packing $\mathcal{K}_{G_{\max}}$ that maximizes packing density over the whole isomorphism class $\mathcal{G}$. Formally expressed,
	\begin{equation*}\label{problemPack}
		\mathcal{K}_{G_{\max}}=\argmax_{\mathcal{K}_{G \in \mathcal{G}}}\rho \left(\mathcal{K}_G \right).
	\end{equation*}
	Here we search not only over the whole plane group $G$
	but also over all rotations and translations of $K$, whose centroid lies in the asymmetric unit of $G$ such that the resulting configuration is a CSG packing.
	
	FIG.~\ref{fig:heptagonFig1} presents examples of the densest plane group packings of a pentagon, heptagon, enneagon, and dodecagon in plane groups $p2$, $p2gg$, $pg$, $p3$, and $p1$.
	
	We consider the densest plane group packing as a nonlinear, constrained, and bounded optimization problem. We transform this problem via stochastic relaxation \cite{geman1984} to the problem of estimation of a probability distribution with the probability mass concentrated on the maxima of the optimization landscape. The two-dimensional lattice group $L$
	associated with a CSG induces a quotient space $\mathbb{R}^2 / L$, which is homeomorphic to a two-dimensional torus. Therefore, we define a parametric family of probability distributions based on the multivariate von Mises distribution \cite{mardia2008}, the Extended Multivariate von Mises distribution (EMvM), and perform a stochastic trust region search on the functional space induced by this family of toroidal probability distributions where the Kullback-Leibler divergence \cite{kullback1951} defines the trust region radius. The resulting Entropic Trust Region Packing Algorithm (ETRPA) is a variant of the natural gradient method \cite{amari}.
	
	At the first iteration, $N$ samples are drawn from a uniform distribution on an $n$-dimensional torus. The number of samples used and dimensionality of torus depend on the crystal system of the plane group where the search is performed. For instance, in the case of the $p2$ group there are six configuration parameters, two for fractional coordinates of the polygon's centroid in the asymmetric unit, one for the angle of rotation of the polygon, and three parameters defining the shape of the primitive cell, that is lengths of lattice generators and an angle between them. For a six-dimensional torus, the dimensionality of the EMvM parametric space is $p=72$, and number of samples drawn at each iteration is set such that $\frac{p}{N}<0.07$, meaning that for the search in $p2$ group $N=1040$. An overlap constraint violation is evaluated for each sampled configuration, and for configurations with no intersections between polygons, packing density is computed. Afterwards, a new batch of samples is generated from the EMvM based on distribution parameters estimates of the largest increase in density and lowest constraint violation. This process is repeated until the algorithm converges to a point distribution or after $8000$ iterations.
	
	There are two main benefits of ETRPA. First, the search does not depend on initial configurations, a caveat of many stochastic search methods. The initial sampling from the uniform distribution on an $n$-dimensional torus provides a satisfactory overview of the optimization landscape induced by the configuration space. Second, by performing the search on an $n$-dimensional torus, we cover problematic instances when the optimal solution lies on the boundary.
	
	To further improve the accuracy of the approximate solutions, after the initial run of the ETRPA search, we perform a refining process by creating progressively smaller $\epsilon$-neighborhoods around the best solution found and employ ETRPA with new boundaries defined by the $\epsilon$-neighborhoods. To prevent the escape of the optimal solution from the $\epsilon$-neighborhood, if a solution with higher density than in previous runs is found, the $\epsilon$-neighborhood is not decreased, only recentered on this solution. Otherwise, the $\epsilon$-neighborhood is decreased. The lowering of $\epsilon$ is repeated $30$ times.
	
	It has to be noted that the entropic trust region is not a physical simulation since particles in the system are allowed to overlap during the search. A detailed treatment of ETRPA is presented for the interested reader in \cite{torda2022entropic}.
	
	\begin{figure}[t]
		\centering
		\includegraphics[trim={20 0 20 10},clip,width=\linewidth]{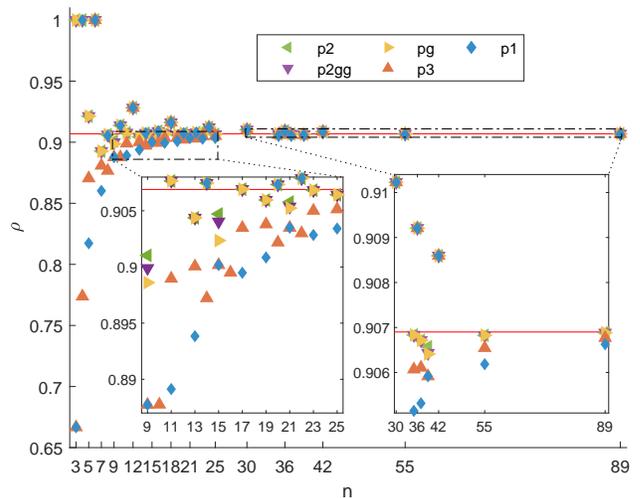}
		\caption{Densities of the densest packings of examined $n$-gons in the$p2 / p2gg / pg / p3 / p1$ plane groups class. The red line denotes the density of the densest disc packing in this plane group class.}
		\label{fig:p2Evo}
	\end{figure}
	
	\section{Results}
	\label{sec:Reuslts}
	
	\begin{turnpage}
		
		\begin{table*}
			\begin{ruledtabular}
				%\tiny
				\setlength{\tabcolsep}{0.5pt}
				\renewcommand{\arraystretch}{1.2}
				\begin{tabular}{ccccccccccccccccccccccccc}
					& 		3   &  4  &  5  &    6  &    7  &    8  &    9  &  10 &  11 &    12 &     13 & 14  & 15  &   16 & 	17 &  18 & 19 & 20 & 21 & 22  & 23 &    24 &    25 &  Disc\footnote{We consider the disc as a limiting n-gon where the number of vertices n approaches infinity.}   \\ \hline
					p2 &    .99999 & .99999 & 	.92131 &	.99999 &	.89269 &	.90616 &	.90103 & .91371 &	.90766 &	.92820 &   	.90437  & .90746 & 	.90471  &	.90901 & 	.90692 & .91622 & .90595 & .90733 &.90577 &	.90789	& .90683 &    .91211 &	.90642 & .90689	\\
					p2gg &  .99999 & .99999 &	.92131 &	.99999 &	.89269 &	.90616 &	.89989 & .91371 &	.90766 &	.92820 & 	.90437  & .90746 & 	.90403  &	.90901 & 	.90692 & .91622 & .90595 & .90733 &.90536 &	.90789	& .90683 &	.91211 &	.90642 & .90689	\\ 
					pg   &  .99999 & .99999 &	.92131 &	.99999 &	.89269 &	.90616 &	.89860 & .91371 &	.90766 &	.92820 & 	.90437  & .90746 & 	.90369  &	.90901 &	.90692 & .91622 & .90595 & .90733 &.90523 &	.90789	& .90683 &	.91211 &	.90642 & .90689	\\
					p3  &   .66666 & .77376 & 	.87048 &	.99999 &	.88085 &	.87665 &	.88773 & .88775 &	.89896 &	.92820 &  	.90005  & .89722 & 	.90017	&	.89950 &	.90348 & .91622 & .90379 & .90218 &.90349 &	.90300	& .90500 &	.91211 &	.90513 & .90689	\\
					p1 &    .66666 & .99999 &	.81725 &	.99999 &	.86019 &	.90616 &	.88773 & .91371 &	.88912 &	.92820 &    .89383  & .90746 & 	.90017  &	.90901 &	.89947 & .91622 & .90084 & .90733 &.90349 &	.90789	& .90284 &	.91211 &	.90341 & .90689	\\ \hline 
					p2mg & 	.99999 & .99999 &	.85410 &	.85714 &	.84226 &	.86555 &	.83419 & .83722 &	.83116 &	.86156 & 	.83306  & .83823 & 	.83993  &	.84856 &	.84215 & .84346 & .84189 & .84571 &.84052 &	.84068	& .83929 &	.84662 &	.83952 & .84178 \\
					cm   & 	.99999 & .99999 &	.85410 &	.85714 &	.84226 &	.86555 &	.83419 & .83722 &	.82795 &	.86156 & 	.83212  & .83823 & 	.83993  &	.84856 &	.84215 & .84346 & .84189 & .84571 &.84052 &	.84068	& .83865 &	.84662 &	.83916 & .84178 \\ 
					p4   &  .71281 & .99999 &	.84211 &	.81776 &	.84219 &	.85031 &	.83030 & .83004 &	.83780 &	.86156 & 	.83691  & .83527 & 	.83765  &	.84613 &	.84190 & .83945 & .84177 & .84367 &.83971 &	.83933	& .84072 &	.84662 &	.84061 & .84178 \\ \hline  
					p4gm & 	.69615 & .99999 &	.71119 &	.74613 &	.76477 &	.82842 &	.76593 & .77205 &	.77628 &	.80384 & 	.77662	& .77869 & 	.78028  &	.79564 &	.78043 & .78137 & .78213 & .79192 &.78220 &	.78271	& .78313 &	.78991 &    .78317 & .78539 \\
					c2mm & 	.66666 & .99999 &	.71714 &	.74999 &	.76253 &	.82842 &	.76697 & .77254 &	.77570 &	.80384 & 	.77697	& .77882 & 	.78006  &	.79564 &	.78058 & .78141 & .78202 & .79192 &.78229 &	.78273	& .78307 &	.78991 & 	.78322 & .78539 \\  
					pm   & 	.49999 & .99999 &	.69098 &	.74999 &	.73825 &	.82842 &	.75712 & .77254 &	.76655 &	.80384 & 	.77193	& .77882 & 	.77530  &	.79564 &	.77754 & .78141 & .77911 & .79192 &.78025 &	.78273	& .78111 &	.78991 &    .78177 & .78539 \\ 
					p2mm &  .49999 & .99999 &	.69098 &	.74999 &	.73825 &	.82842 &	.75712 & .77254 &	.76655 &	.80384 & 	.77193  & .77882 & 	.77530  &	.79564 &	.77754 & .78141 & .77911 & .79192 &.78025 &	.78273	& .78111 &	.78991 &	.78177 & .78539 \\\hline
					p6 &    .99999 & .72193 &	.75933 &	.85714 &	.75740 &	.76438 &	.78535 & .76932 &	.77293 &	.79560 &    .77254  & .77326 & 	.77997	&	.77425 &	.77536 & .78533 & .77523 & .77536 &.77863 &	.77571	& .77622 &	.78181 &	.77616 & .77734	\\
					p31m & 	.74999 & .66323 &	.70166 &	.71999 &	.72084 &	.71565 &	.73410 & .71653 &	.72182 &	.74613 & 	.72205  & .72380 & 	.73087  &	.72583 &	.72735 & .72825 & .72787 & .72698 &.72996 &	.72662	& .72722 &	.73320 &	.72725 & .72900 \\ 
					p3m1 & 	.99999 & .49742 &	.53854 &	.66666 &	.57196 &	.57980 &	.63041 & .58887 &	.59164 &	.61880 &	.59535	& .59664 & 	.61359  &	.59851 &	.59921 & .61081 & .60029 & .60071 &.60915 &	.60139	& .60166 &	.60807 &	.60211 & .60459	\\ 
					p4mm & 	.57735 & .49999 &	.55537 &	.52148 &	.51446 &	.56854 &	.52443 & .53314 &	.54240 &	.53589 &    .54144	& .53607 & 	.53385  &	.54604 &	.53500 & .53725 & .54014 & .53792 &.53994 & .53783	& .53683 &   .54211 &    .53716 & .53901 \\ 
					p6mm & 	.49999 & .46410 &	.49372 &	.47999 &	.49025 &	.48235 &	.48940 & .48305 &	.47750 &	.49742 &    .47995	& .48454 & 	.48724  &	.48518 &	.48675 & .48550 & .48660 & .48548 &.48664 &	.48542	& .48409 &	.48880 &	.48418 & .48600 \\ 
				\end{tabular}
			\end{ruledtabular}
			\caption{\label{tab:table1} 
				Densities of the densest plane group packing configurations of $n$-gons obtained by ETRPA for $n=3,4,\ldots,25$. The densities are truncated at the fifth decimal place.}
		\end{table*}
	\end{turnpage}
	
	Using the ETRPA, we recovered known, as well as obtained previously unknown densest packings of $n$-gons for $n = 3, 4, \ldots, 25, 30, 35, 36, 37, 39, 42, 55, 89$ in all $17$ plane groups, including a disc regarded as a limiting $n$-gon when the number of vertices $n$ approaches infinity. The densities of the densest configurations obtained in our experiments for $n=3,4,\ldots,25$ are shown in Table~\ref{tab:table1} \footnote{ See Supplemental Material at  \href{https://milotorda.net/wp-content/uploads/supplemental_material.pdf}{milotorda.net} for a complete table of densest plane group packings of $n$-gons, including configuration parameters and visualizations of respective structures and for exact densities of a discs densest plane group packings and additional information on relationships between densest plane group packings of a disc and $n$-gons with $24$-fold rotational symmetry}. All the values are truncated at the fifth decimal place.	
	
	\begin{figure*}[t]
		\centering
		
		\begin{minipage}[b]{7.4cm}
			\centering
			\includegraphics[trim={0.1 0.1 0.1 0.1},clip,height=2.5cm]{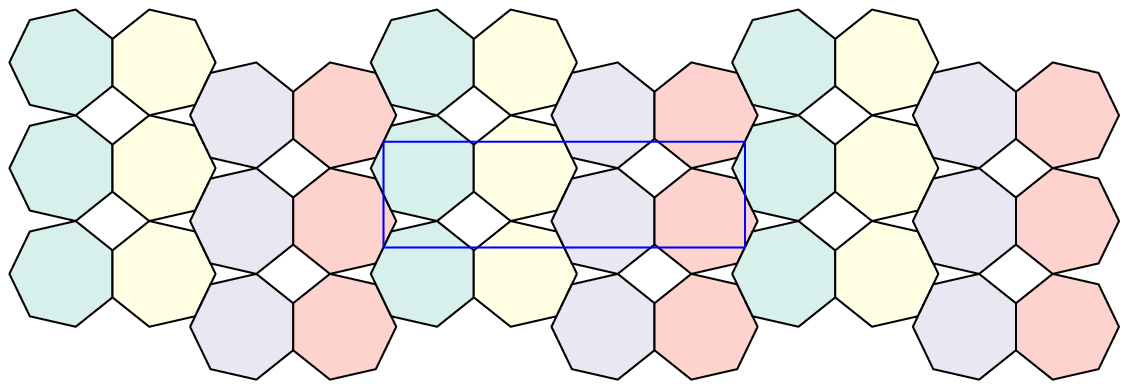}
			
			\includegraphics[trim={0.1 0.1 0.1 0.1},clip,height=2.3cm]{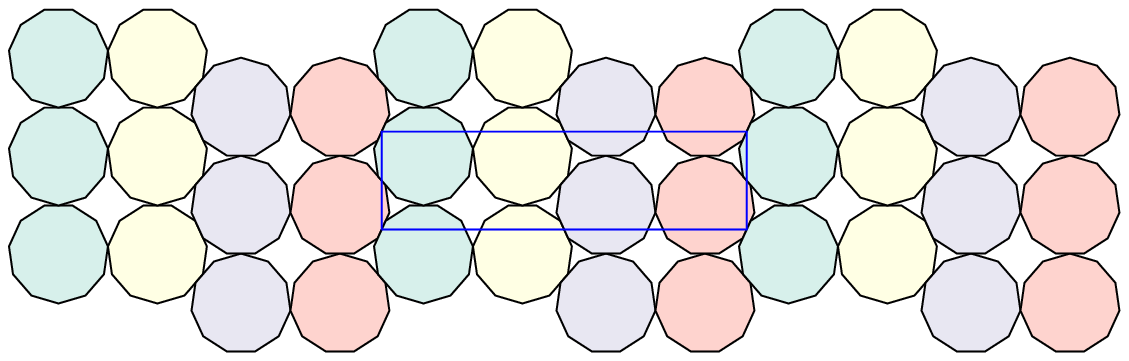}
			
			\includegraphics[trim={0.1 0.1 0.1 0.1},clip,height=2.3cm]{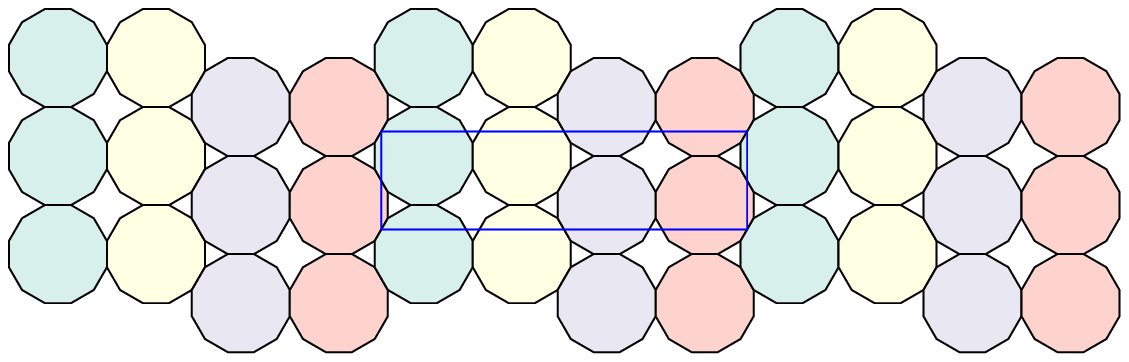}
			\subcaption{p2mg}
		\end{minipage}	
		\begin{minipage}[b]{7.1cm}
			%\vspace{-6cm}
			\centering
			\includegraphics[trim={0.1 0.1 0.1 0.1},clip,height=2.5cm]{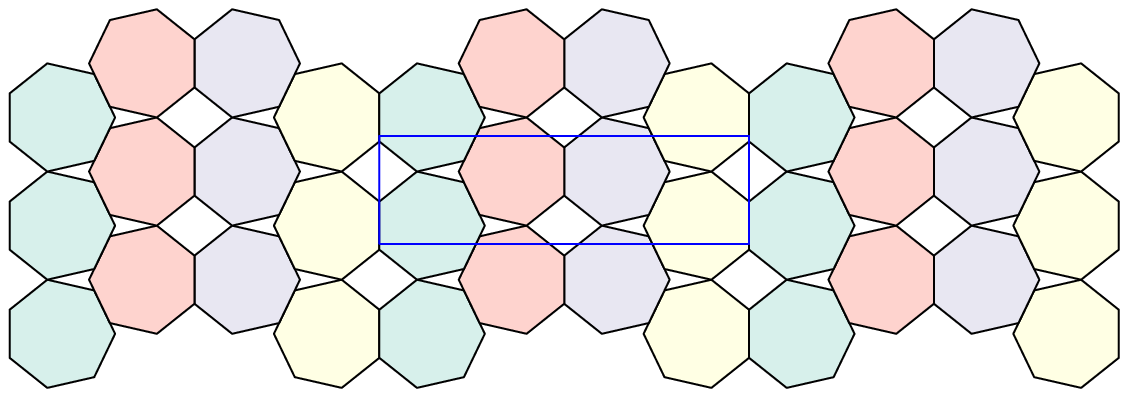}
			
			\includegraphics[trim={0.1 0.1 0.1 0.1},clip,height=2.4cm]{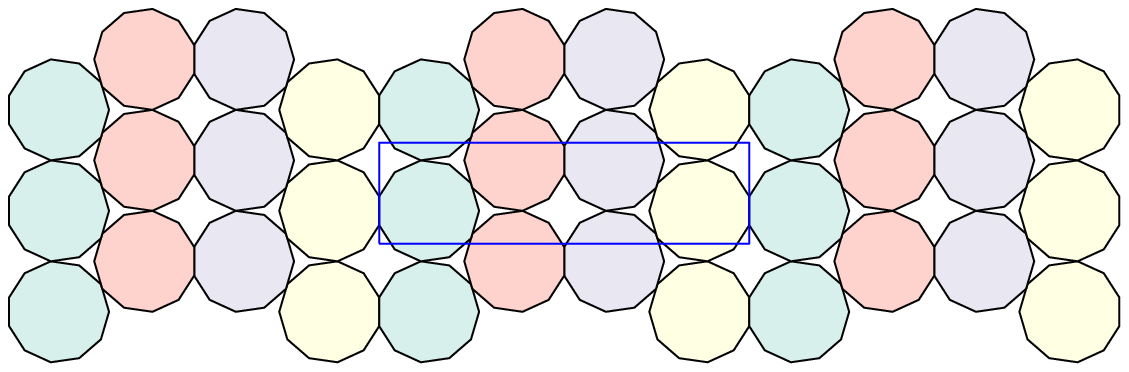}
			
			\includegraphics[trim={0.1 0.1 0.1 0.1},clip,height=2.3cm]{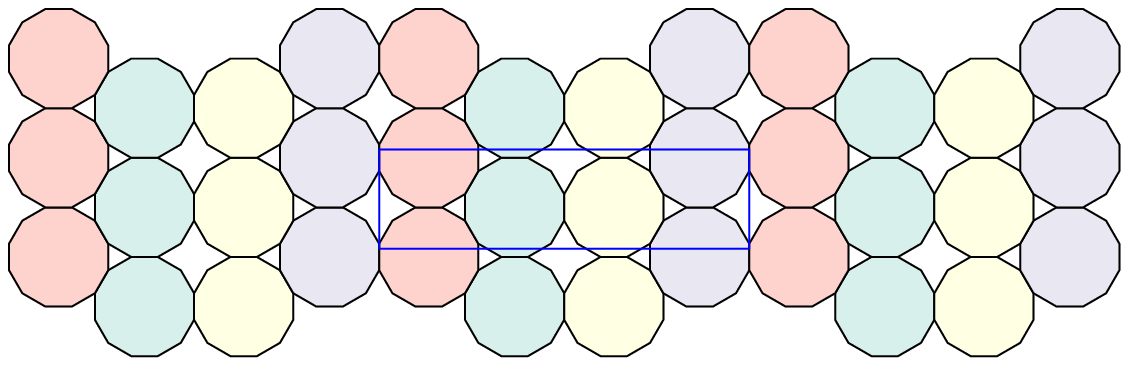}
			\subcaption{cm} 
		\end{minipage}		
		\begin{minipage}[b]{3cm}
			%\vspace{-6cm}
			\centering
			\includegraphics[trim={0.1 0.1 0.1 0.1},clip,height=2.4cm]{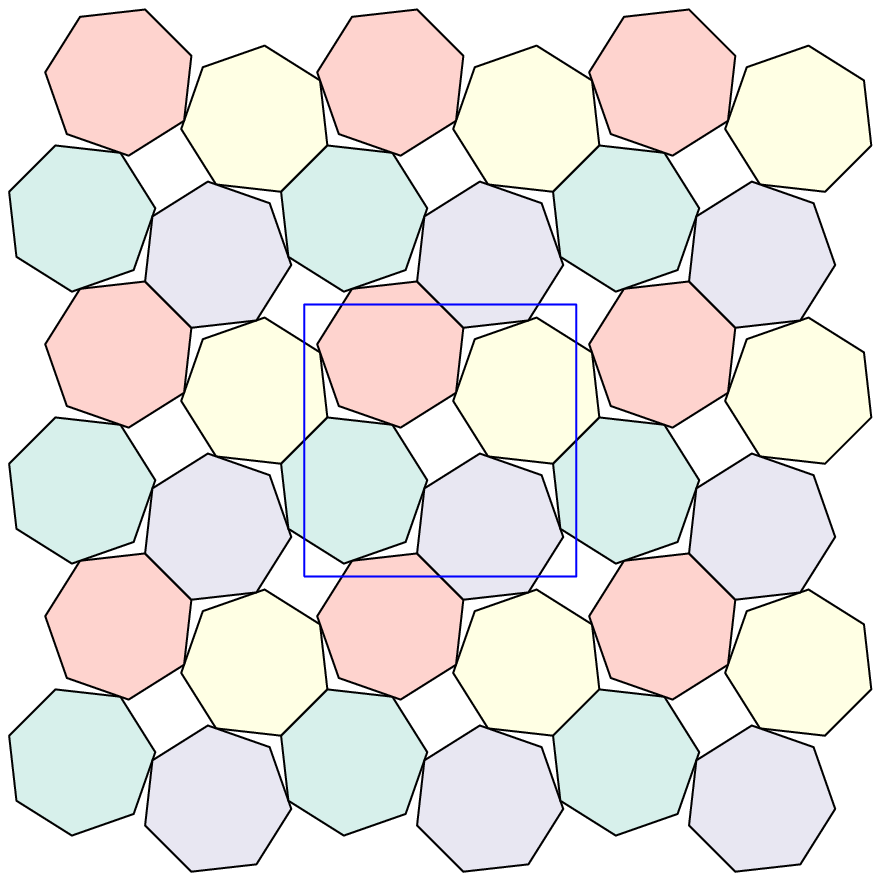}	
			
			\includegraphics[trim={0.1 0.1 0.1 0.1},clip,height=2.4cm]{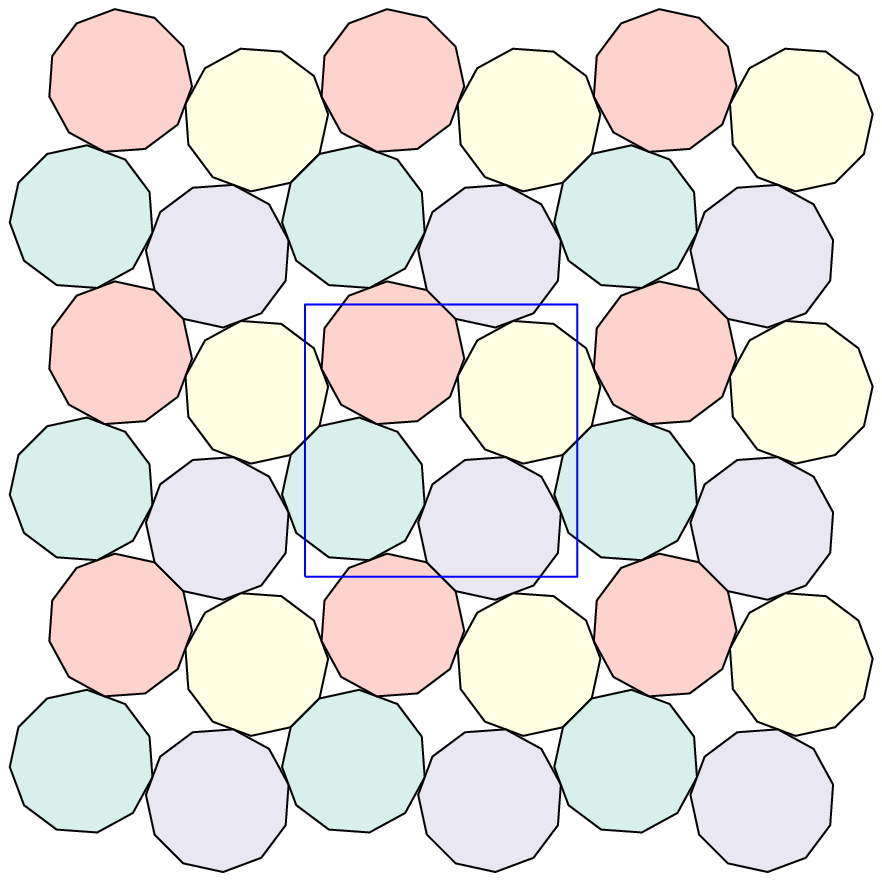}
			
			\includegraphics[trim={0.1 0.1 0.1 0.1},clip,height=2.4cm]{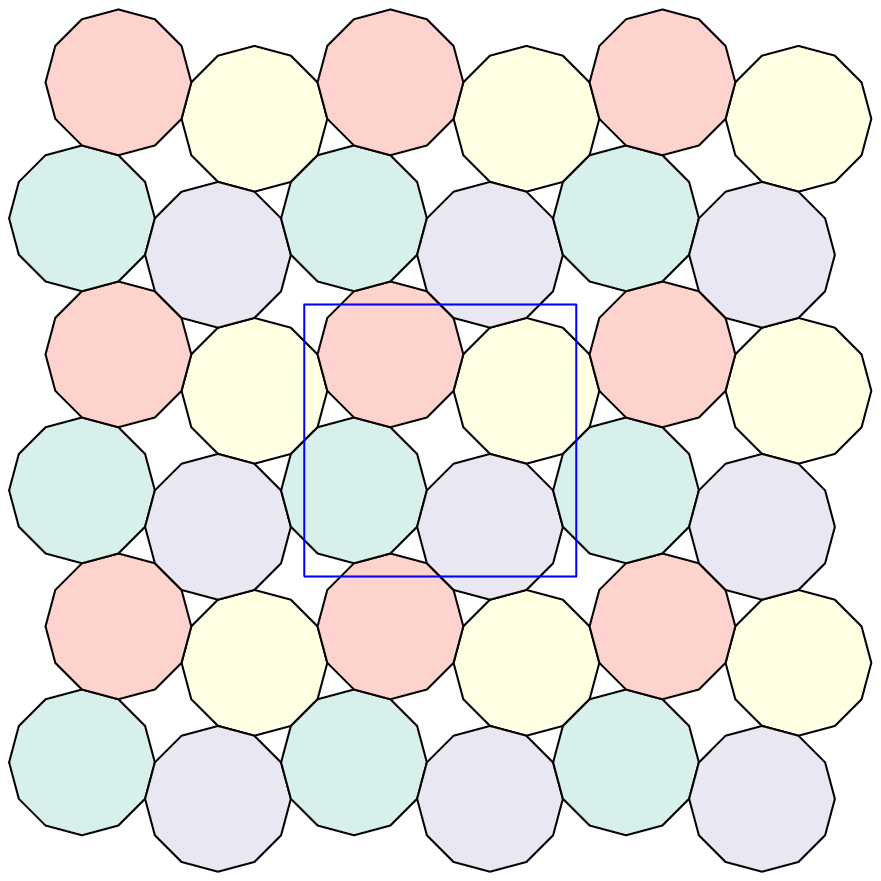}
			\subcaption{p4} 
		\end{minipage}
		
		\caption{\label{fig:heptagonFig2} Densest configurations of (top) heptagon, (middle) endecagon, and (bottom) dodecagon in plane groups $p2mg$, $cm$, and $p4$ with the following densities: heptagon in $p2mg / cm \approxeq 0.84226$ and $p4 \approxeq 0.84219$; endecagon in $p2mg \approxeq 0.83116$, $cm \approxeq 0.82795$ and $p4 \approxeq 0.83780$; dodecagon in $p2mg / cm / p4 \approxeq 0.86156$. The blue parallelogram denotes the primitive cell of the respective configuration. Colors represent symmetry operations modulo lattice translations.}
	\end{figure*}

	In the following paragraphs, we examine these configurations classified according to the disc's densest plane group packings. Since the symmetry group of a circle contains symmetries of all $n$-gons, our results indicate that the $24$-fold rotational symmetry of an $n$-gon is sufficient to constitute the optimal plane group configurations of a disc and further suggests a relationship with the plane group symmetries. Notably, for the packing densities to be equal in the plane group class $p2/p2gg/pg/p3/p1$, two-fold and three-fold rotational symmetries are necessary, examined in Section~\ref{sec: p2/p2gg/pg/p3/p1}. In the class $p2mg/cm/p4$, a four-fold rotational and a local three-fold symmetry is necessary, examined in Section~\ref{sec:p4gm/c2mm/pm/p2mm}, and in the class $p4gm/c2mm/pm/p2mm$, a four-fold rotational symmetry is necessary, examined in Section~\ref{sec:p4gm/c2mm/pm/p2mm}. In the context of the crystallographic restriction theorem \cite{elliott1998}, which states that periodic crystals can only have two-fold, three-fold, four-fold, and six-fold rotational symmetries, minimal rotational symmetry containing all preceding is $12$-fold. However, a $12$-fold rotational symmetry of an $n$-gon does not cover a local eight-fold rotational symmetry that is present in a disc's optimal $p4mm$ packing, as demonstrated in Section~\ref{sec:p6/p31m/p3m1/p4mm/p6mm}. Therefore, minimal symmetry containing the symmetries mentioned earlier is a $24$-fold rotational symmetry.
	
	\subsection{Densest $p2$, $pg$, $p2gg$, $p3$, and $p1$ packings}
	\label{sec: p2/p2gg/pg/p3/p1}
	
	It is known that the packing density of the densest packing of a disc is$\frac{\pi}{\sqrt{12}} \approx  0.9068996 \dots$ \cite{toth2013}. This density was attained as the densest plane group packing of a disc in groups$p2$, $p2gg$, $pg$, $p3$, and $p1$.
	
	The highest packing densities among all plane groups were observed in the plane group $p2$ or all examined $n$-gons. The $p2$ group is sometimes referred to as a double lattice since it can be viewed as a collection of two lattices related by a two-fold rotational symmetry.
	
	Our results are consistent with known densest packings of polygons, that is, the uniform triangular, square, and hexagonal tilings, densest known packing of a pentagon \cite{hales}, heptagon \cite{kuperbergSup}, octagon \cite{atkinson}, enneagon \cite{de2011}, and disc \cite{fejes1942dichteste}. Moreover, given that the double lattice packing is at least locally optimal for convex polygons in the space of all packings \cite{kallus2016local} and that the plane group packings are inherently periodic, our results support the optimality of $p2$ packing among all plane group packings. 
	
	Additionally, for the space of plane group packing configurations, we experimentally verified the conjecture that the densest $p2$ packing of the heptagon is less than any other shape \cite{kallus2015pessimal}. The extremality of the heptagon can be generalized to all $n$-gons since the $p2$ packing density converges to the optimal packing density of a disc when the number of vertices is increased, as is shown in FIG.~\ref{fig:p2Evo}. Moreover, our results suggest that for every $n$-gon such that $n > 6k+1=m$ and $k \in \mathbb{N}$ the densest $p2$ $n$-gon packing is strictly higher than the densest $p2$ $m$-gon packing.   

Further, our results show that for some $n$-gons, the densest plane group packings in groups$p2$, $p2gg$, and $pg$ are equal. Thus the densest known configurations of pentagon and heptagon have higher symmetry than a double lattice, and these configurations can be realized using a glide reflection instead of a rotation by $\pi$ around the center of symmetry of the $p2$ group. This observation holds for every $n$-gon we examined, where the number of vertices $n$ of a given $n$-gon is not equal to $3k$ where $k=2,3,\ldots$.

It is known that the group $p2gg$ has three maximal nonisomorphic subgroups, one with $p2$ symmetry and two with $pg$ symmetry \cite{brock2016}. Due to the mirror symmetry of $n$-gons and the equality of densest $p2$ and $pg$ configurations, in the densest $p2gg$ configuration, the $p2$ subgroup induces an additional $pg$ symmetry, and the glide reflection plane of one of the $pg$ subgroups induces an additional $p2$ symmetry. Thus, the densest $p2gg$ configuration coincides with the densest $p2$ and $pg$ configurations in these cases.

The only instances where the equality between maximal densities in groups $p2$, $p2gg$, and $pg$ does not hold is for $n$-gons with three-fold but no central symmetry where $n\geq 9$. The lowest densities of these three packing configurations were attained in the $pg$ group. Furthermore, any densest $pg$ configuration of an $n$-gon can be easily converted to a $p2gg$ and $p2$ packing with the same density, which means that the densest $pg$ packings for a fixed $n$
can serve as a lower bound for $p2$ and $p2gg$ densest packings. In fact, this lower bound was attained as densest $p2$ and $p2gg$ packings for all but centrally non-symmetric $n$-gons containing a three-fold rotational symmetry and $n\geq 9$. The additional symmetry operations in $p2gg$ group and the additional degree of freedom of the $p2$ crystal system (oblique) allow for higher density packing than that of $pg$ in these cases.

The densest packing of a convex compact subset of the two-dimensional Euclidean space with central symmetry is that of a lattice packing \cite{rogers1951}. In the crystallographic setting, it is the plane group $p1$, which is a group containing only lattice translations. Moreover, any lattice packing of a centrally symmetric convex polygon can be easily converted to a $p2$ packing with the same density \cite{mount1991}, meaning that for centrally symmetric $n$-gons, densities of the densest $p1$ and $p2$ packings are equal. Indeed, in our experiments, the densest $p2$ packings of centrally symmetric $n$-gons attained the same approximate highest density as in $p1$. Moreover, we obtained the same densities as in the densest $p2gg$ and $pg$ packings. These observations suggest that optimal packings of centrally symmetric polygons can also be realized as either two lattices related to each other by a glide reflection or as four lattices related to each other by two glide reflections and two-fold rotational symmetry. Consequently, optimal packings of centrally symmetric $n$-gons have higher symmetry than that of a lattice or double lattice packing.

\begin{figure}[t]
	\centering
	\includegraphics[trim={0 0 0 0},clip,height=0.13\textwidth]{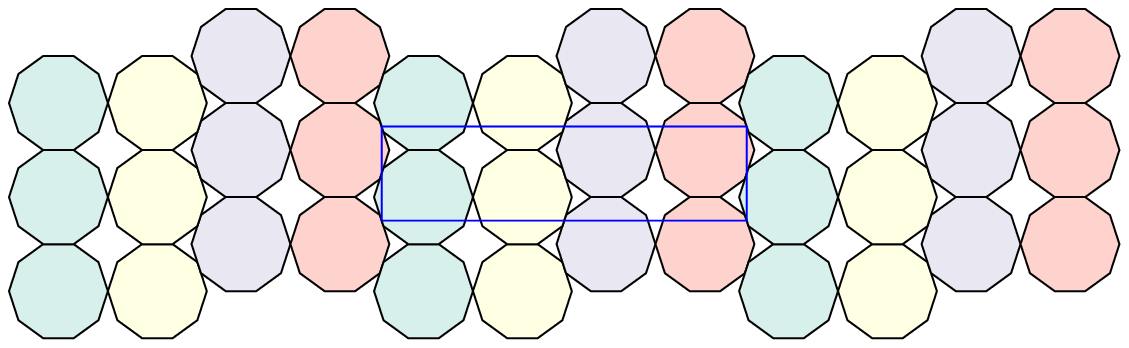}
	
	\includegraphics[trim={54 100 41 100},clip,height=0.156\textwidth]{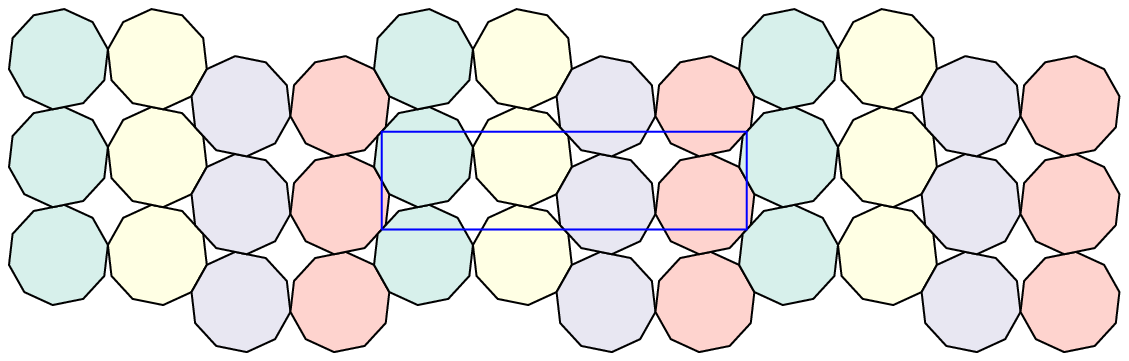}
	\caption{\label{fig:dekagonP2MGalt} Two nonisomorphic densest $p2mg$ packings of a decagon obtained by ETRPA with packing density of approximately $0.83722$.}
\end{figure}

Concerning the densest $p3$ and $p1$
packings of $n$-gons with three-fold rotational symmetry, our results show that the densities in these instances are equal. Moreover, combining two-fold and three-fold rotational symmetries for $n$-gons, where the number of vertices $n$
is divisible by six, the densities of the densest packings in groups $p2$, $p2gg$, $pg$, $p3$, and $p1$ are equal. Consequently, in our experiments, the symmetries of optimal packing configurations of $n$-gons with a six-fold rotational symmetry coincided with the symmetries of optimal packing configurations of a disc.

Additionally, it is known that the densest $p1$
packing of a regular triangle has the lowest density among all densest $p1$
configurations of two-dimensional convex shapes \cite{courant1965least}. Combined with the observation that for $n$-gons without central symmetry, densities of the densest $p3$
packings are greater or equal to their respective densest $p1$
packing densities and supported by the convergence of densest $p3$
packing densities to the optimal packing density of a disc, shown in Figure \ref{fig:p2Evo}, suggests that the regular triangle also minimizes the maximum density in the group $p3$.

\begin{figure}[!t]
	\centering
	\includegraphics[trim={20 0 20 10},clip,width=\linewidth]{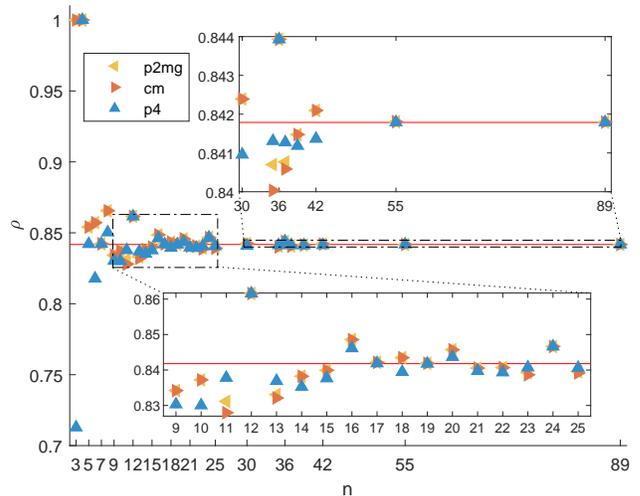}
	\caption{
	Densities of the densest packings of examined $n$-gons in the $p2mg / cm / p4$ plane groups class. The red line denotes the density of the in this plane group class.}
	\label{fig:p2mgEvo}
\end{figure}

\subsection{Densest $p2mg$, $cm$, and $p4$ packings}
\label{sec:p2mg/cm/p4}
	
	\begin{figure*}[!t]
		\centering	
		\begin{minipage}[b]{2.8cm}
			%\hspace{-2cm}
			\centering
			\includegraphics[trim={0.1 0.1 0.1 0.1},clip,height=2.8cm]{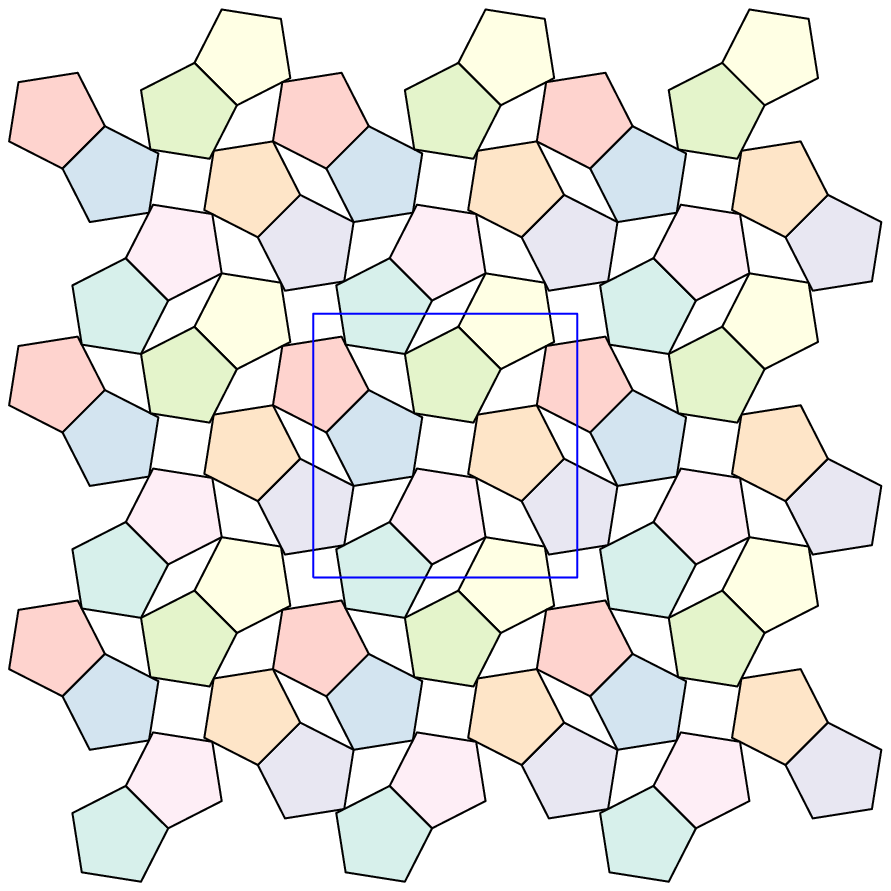}
			
			\includegraphics[trim={0.1 0.1 0.1 0.1},clip,height=2.8cm]{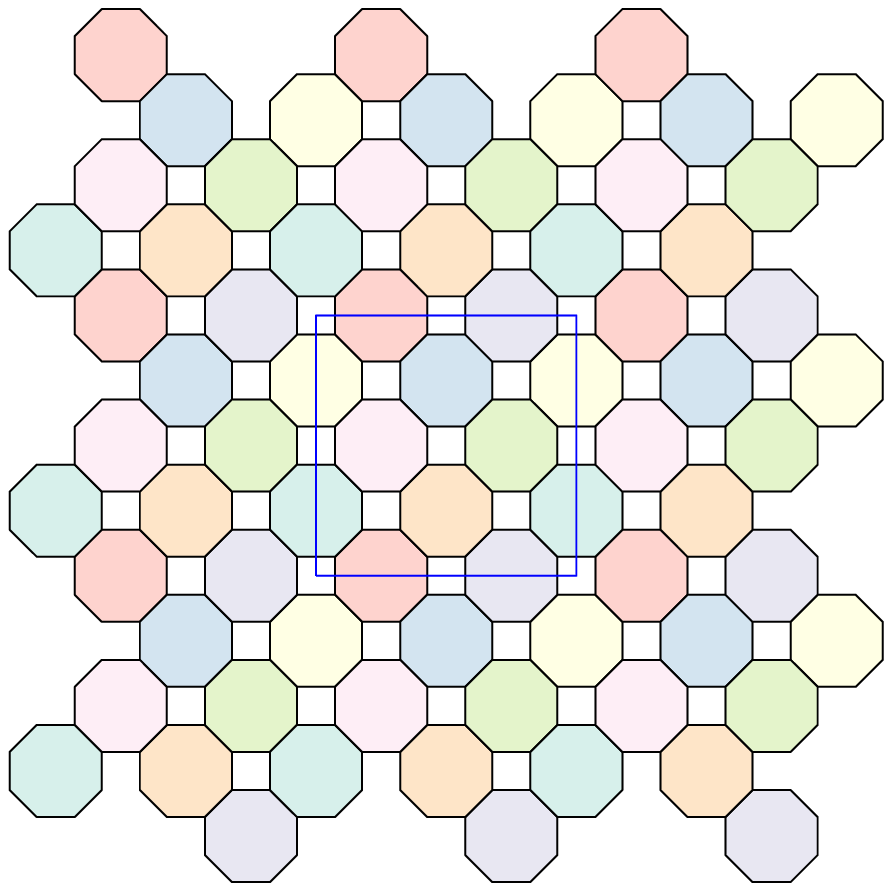}
			
			\includegraphics[trim={0.1 0.1 0.1 0.1},clip,height=2.8cm]{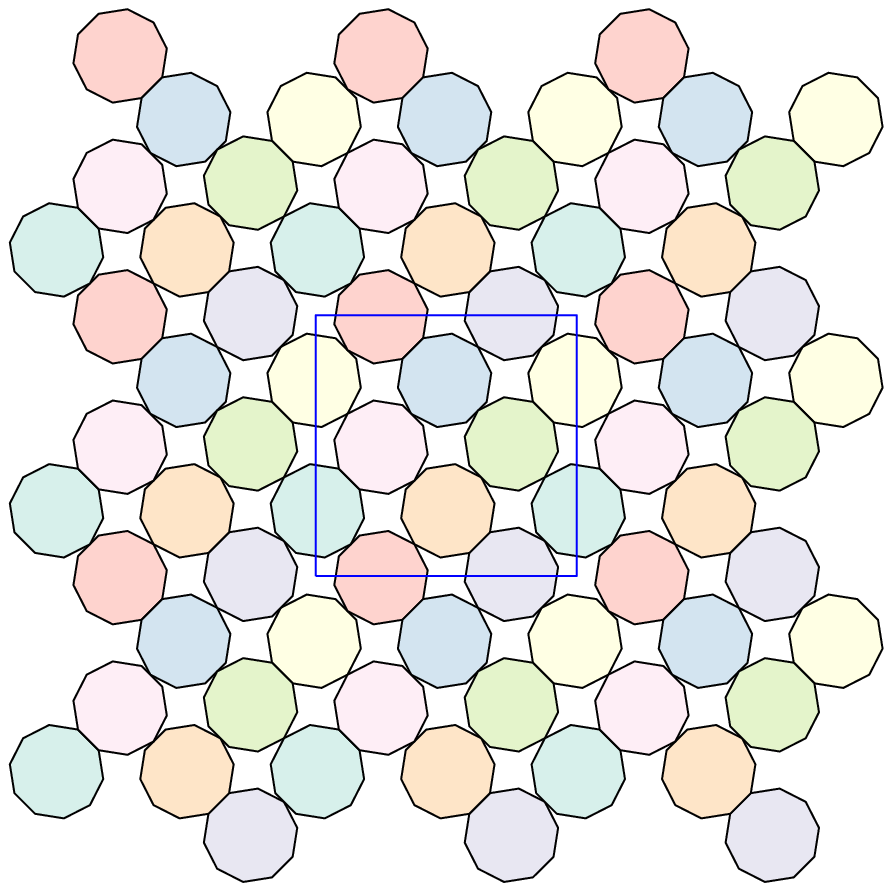}
			\subcaption{p4gm}
		\end{minipage}
		\begin{minipage}[b]{5.9cm}
			%\vspace{-6cm}
			\centering
			\includegraphics[trim={0.1 0.1 0.1 0.1},clip,height=2.8cm]{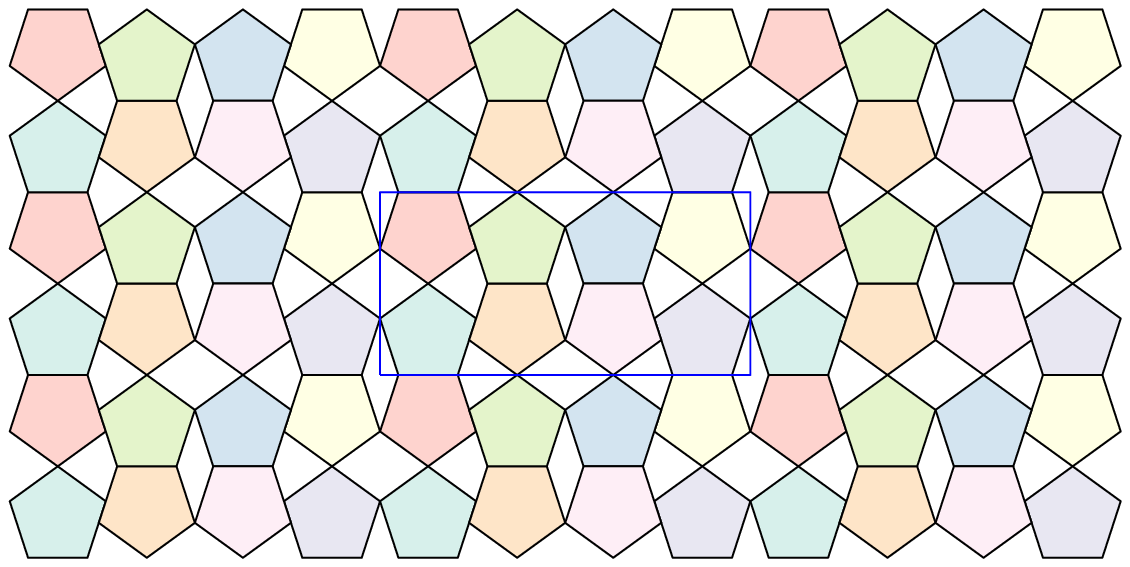}
			
			\includegraphics[trim={0.1 0.1 0.1 0.1},clip,height=2.8cm]{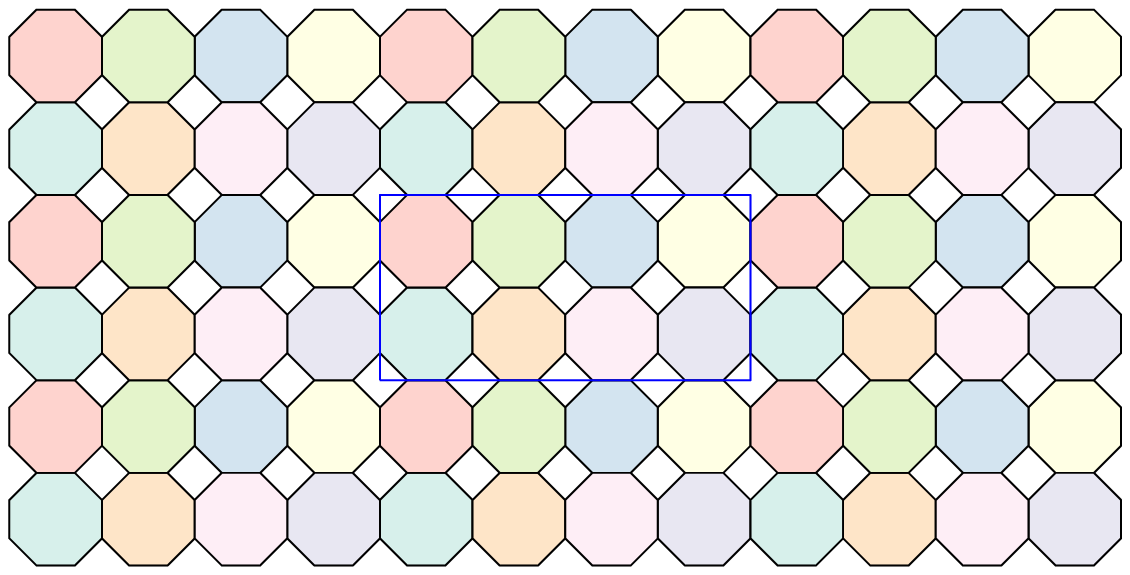}
			
			\includegraphics[trim={0.1 0.1 0.1 0.1},clip,height=2.8cm]{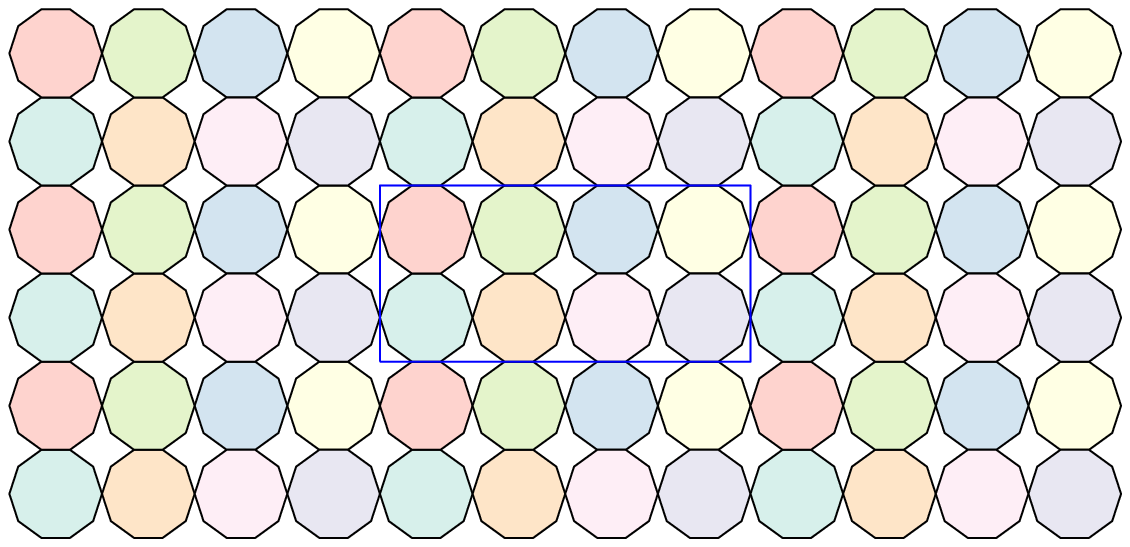}
			\subcaption{c2mm} 
		\end{minipage}		
		\begin{minipage}[b]{3cm}
			%\hspace{-2cm}
			\centering
			\includegraphics[trim={0.1 0.1 0.1 0.1},clip,height=2.8cm]{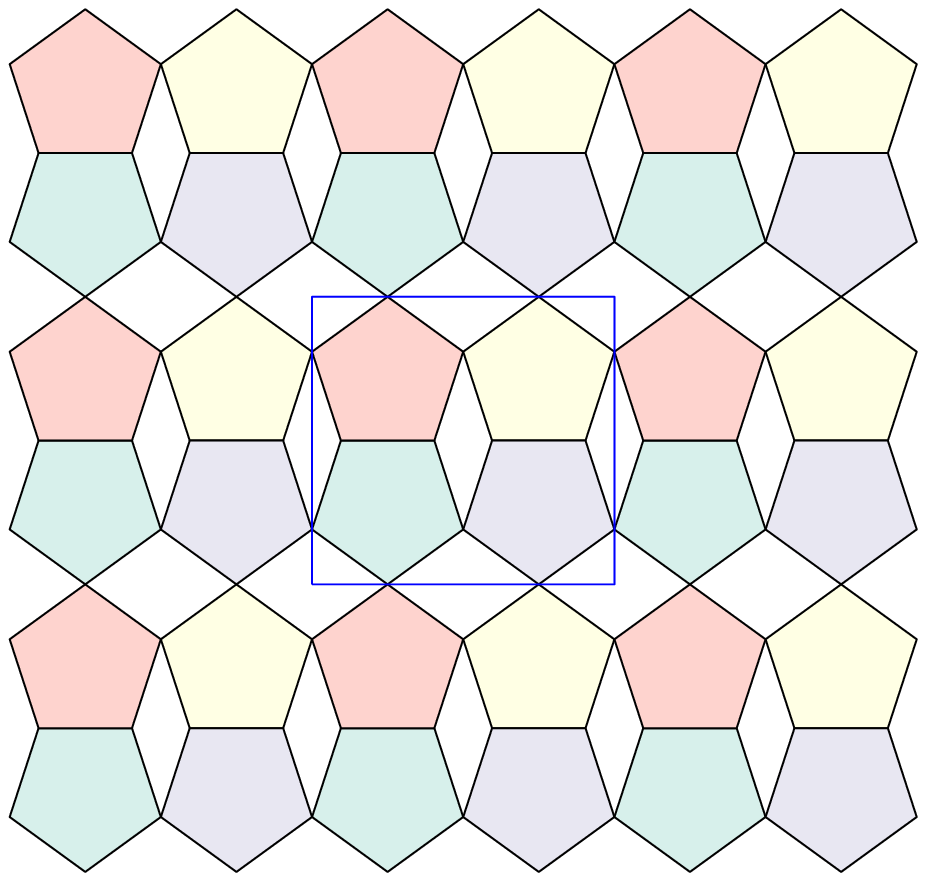}
			
			\includegraphics[trim={0.1 0.1 0.1 0.1},clip,height=2.8cm]{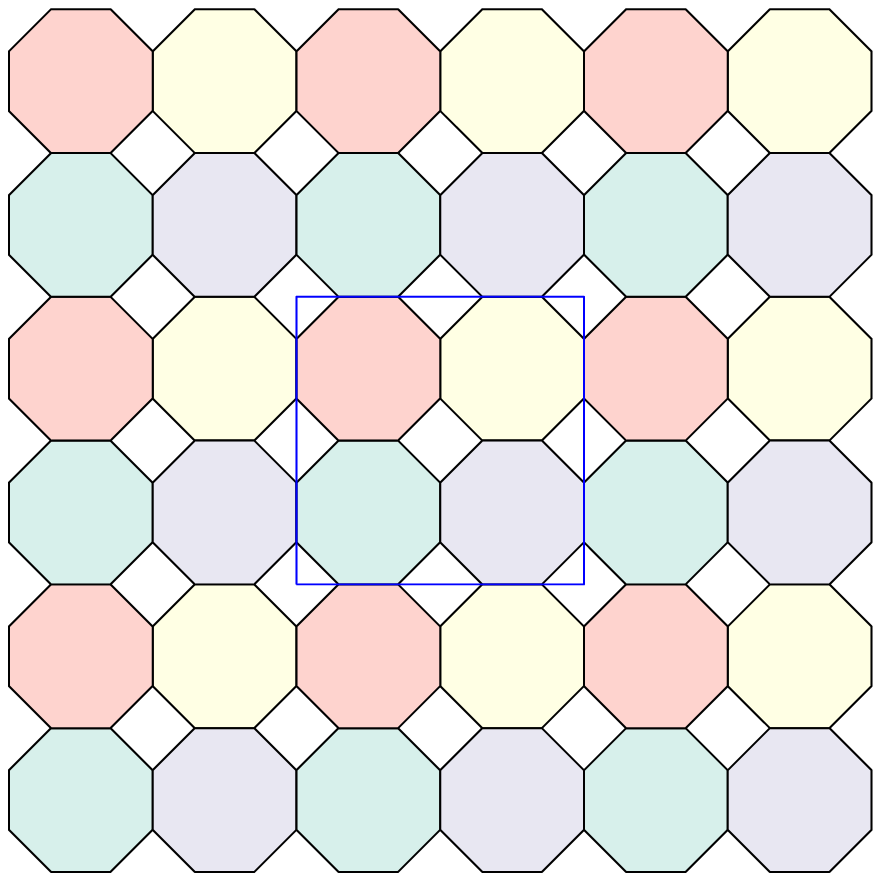}
			
			\includegraphics[trim={0.1 0.1 0.1 0.1},clip,height=2.8cm]{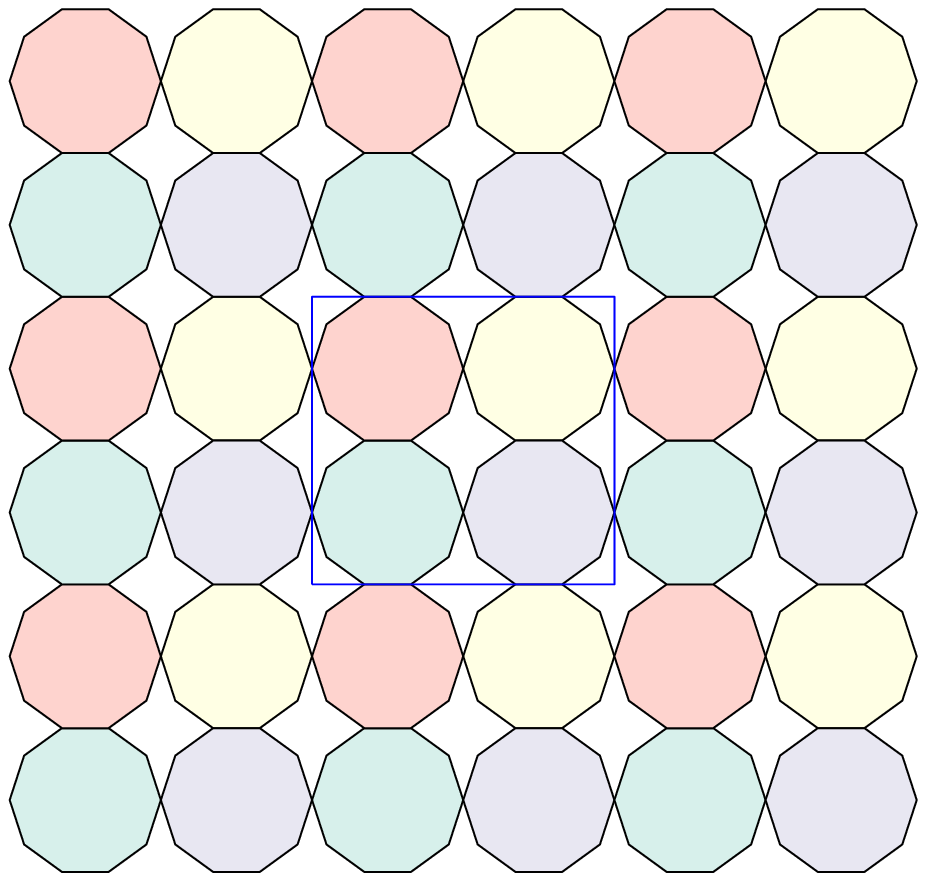}
			\subcaption{p2mm} 
			\label{subfig:p2mmFig}
		\end{minipage}
		\begin{minipage}[b]{5cm}
			%\vspace{-6cm}
			\centering
			\includegraphics[trim={0.1 0.1 0.1 0.1},clip,height=2.8cm]{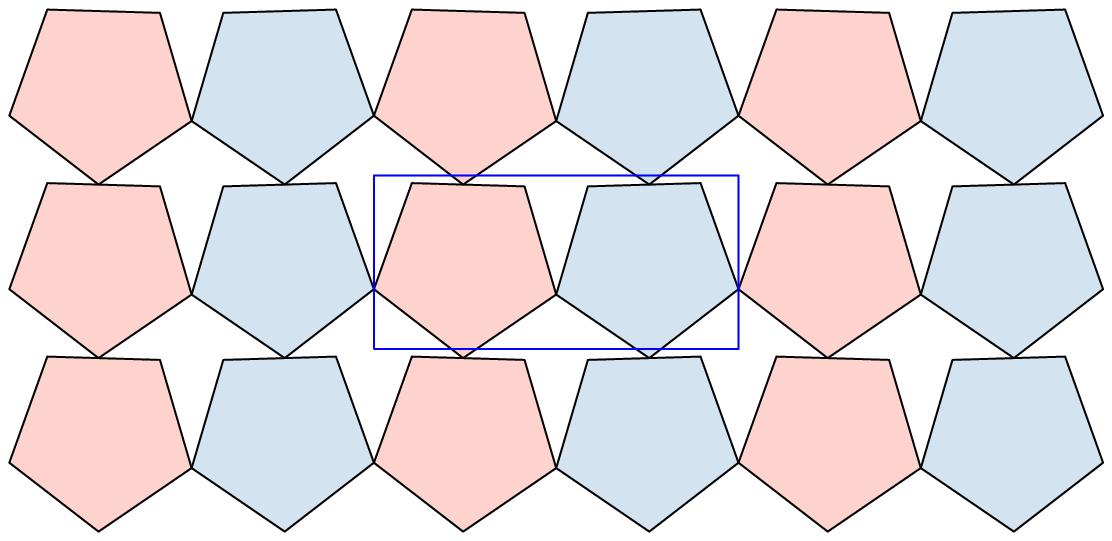}
			
			\includegraphics[trim={0.1 0.1 0.1 0.1},clip,height=2.8cm]{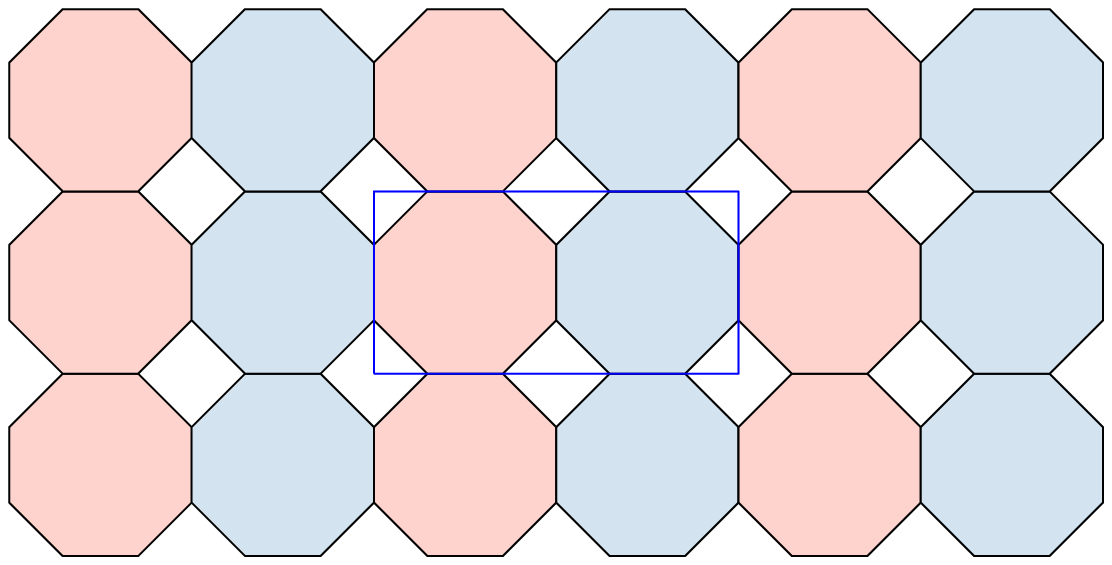}
			
			\includegraphics[trim={0.1 0.1 0.1 0.1},clip,height=2.8cm]{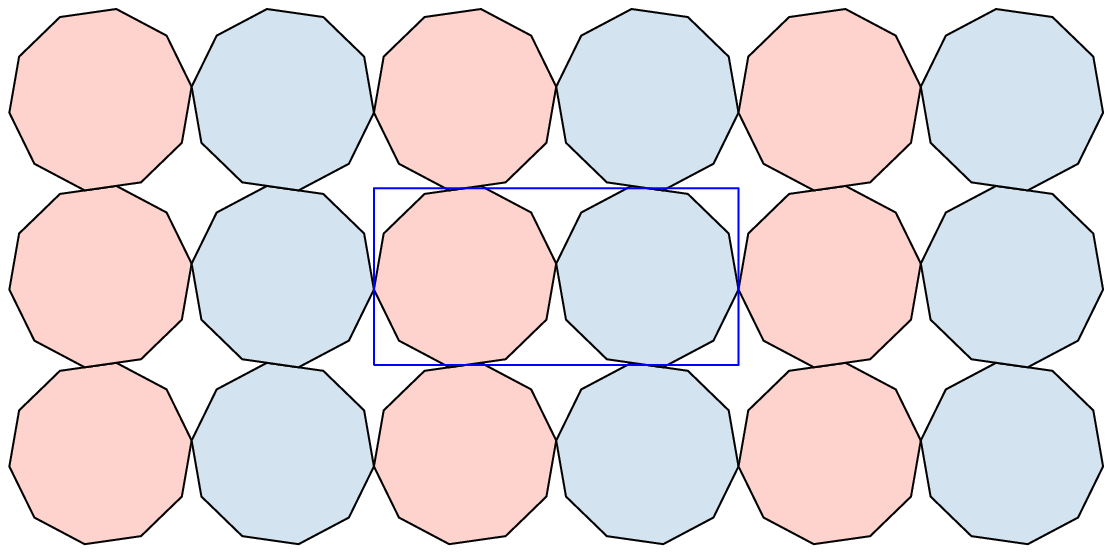}
			\subcaption{pm}
			\label{subfig:heptagonFig3}
		\end{minipage}
		
		\caption{\label{fig:heptagonFig3} 
		Densest configurations of (top) pentagon, (middle) octagon, and (bottom) decagon in plane groups $p4gm$, $c2mm$, $p2mm$ and $pm$ with the following densities: pentagon in $p4gm \approxeq 0.71119$, $c2mm \approxeq 0.71714$ and $p2mm / pm \approxeq 0.69098$; octagon in $p4gm / c2mm / p2mm / pm \approxeq 0.82842 $; decagon in $p4gm \approxeq 0.77205$ and $c2mm / p2mm / pm \approxeq 0.77254 $. The blue parallelogram denotes the primitive cell of the respective configuration. Colors represent symmetry operations modulo lattice translations.}
	\end{figure*}

	The packing densities of a disc in plane groups $p2mg$, $cm$, and $p4$
	are equal, with a density of approximately $0.8938363$. Our results suggest this is also true for $n$-gons with $12$-fold rotational symmetry. The difference between these $p2mg$/$cm$, and $p4$ packing configurations is that the $p4$ and $p2mg/cm$ packing configurations are not isometric, as can be observed by visual comparison in the example of a dodecagon in FIG.~\ref{fig:heptagonFig2}. However, there is a local three-fold rotational symmetry of dodecagonal trimers and four-fold rotational symmetry of dodecagonal tetramers present in $p2mg$, $cm$ and $p4$ configurations. The densities of densest $p2mg$ and $cm$ packings are equal for all 
	$n$-gons we examined except for those where the number of vertices is close to a polygon with $12$-fold rotational symmetry. Precisely, for $12k-1$ and $12k+1$, where $k \in \mathbb{N}$. In fact, they are isometric. On the other hand, for centrally symmetric $n$-gons there are at least two nonisometric densest $p2mg$ configurations shown on the example of decagon in FIG.~\ref{fig:dekagonP2MGalt}. Consequently, for centrally symmetric $n$-gons without four-fold rotational symmetry there exist densest nonisometric $p2mg$ and $cm$ configurations with equal density.
	
	Visually comparing the heptagon $p2mg$ and $cm$ packings in FIG.~\ref{fig:heptagonFig2}, the structure of both packing configurations is similar in that the mirror symmetry planes of polygons are orthogonal to the mirror symmetry planes of the $cm$ group. This is not the case for the $p2mg$ packing configuration of the endecagon where the planes of mirror symmetries of polygons are parallel to the $p2mg$ mirror symmetry planes of reflection. In fact, it is not hard to convert the densest $cm$ packings in Table~\ref{tab:table1} to $p2mg$ packings of equal density, indicating that the densest $cm$ packing configurations of $n$-gons are lower or equal to the densest $p2mg$ packing configurations. The $12k-1$ or $12k+1$ symmetry of the polygon where $k \in \mathbb{N}$ allows for a slightly higher density in $p2mg$ group than in the $cm$ group.
	
	Furthermore, the densities of the densest $p2mg/cm$ packings are greater or equal to $p4$ packings with equality for $n$-gons with $12$-fold rotational symmetry, except in cases of $n$-gons with $12k-1$ and $12k+1$ symmetries, where densest $p4$ packing densities are above corresponding $p2mg$/$cm$ packing densities. These intricacies of $p2mg/cm/p4$ packing configurations are visually presented as a colored rank table in FIG.~\ref{fig:tabelFig6}.
	
	The highest packing density in groups $p2mg$ and $cm$ was attained by the triangle and square where both polygons tile the two-dimensional Euclidean plane, and the highest density $p4$ packing configuration was attained by one of the uniform tilings by a square. The lowest packing density in plane groups $p2mg$ and $cm$ was observed in the case of the endecagon. The lowest density $p4$ packing was attained by the triangle, although higher than the densities of the densest packings of the triangle in the $p1$ and $p3$ plane groups. From the convergence of $n$-gon packings to the densest $p2mg/cm/p4$ configurations of a disc, presented in FIG.~\ref{fig:p2mgEvo}, it is reasonable to assume that the lowest densest packing in this plane group class is attained by the endecagon in $p2mg/cm$ and triangle in $p4$ for all $n$-gons.

	\subsection{Densest $p4gm$, $c2mm$, $pm$, and $p2mm$ packings}
	\label{sec:p4gm/c2mm/pm/p2mm}
	
	Densest disc packings in groups $p4gm$, $c2mm$, $pm$, and $p2mm$ have an equal packing density of $0.7853981$, constituting another class of plane groups. The equality between packing density in this plane group class was also true when a four-fold rotational symmetry was present in an $n$-gon. In fact, all four configurations are isometric, as seen in the example of an octagon in FIG.~\ref{fig:heptagonFig3}.
	
\begin{figure}[!t]
		\centering
		\includegraphics[trim={1 0 1 0},clip,height=0.11\textwidth]{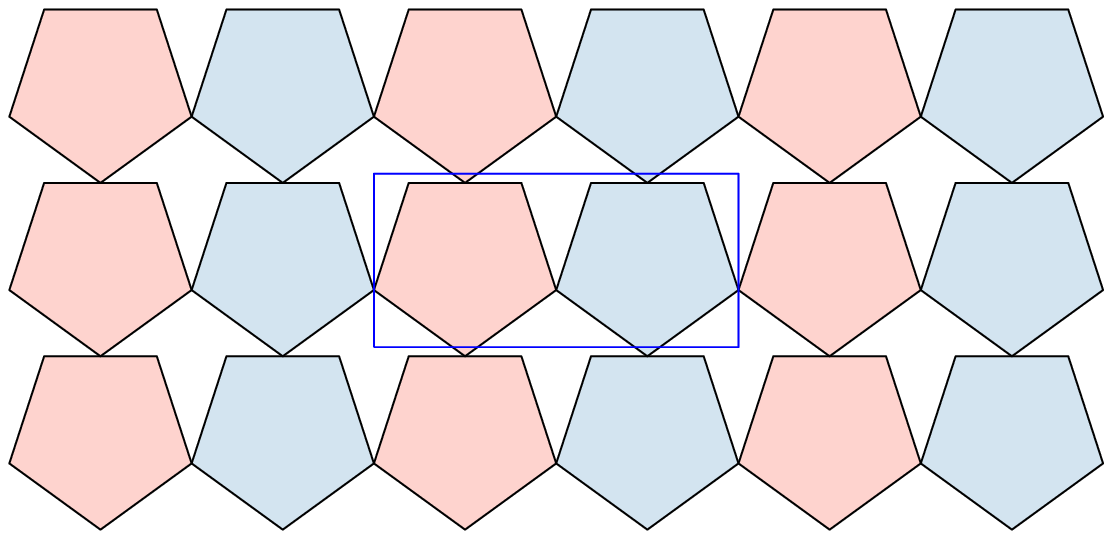}	
		\includegraphics[trim={1 0 1 0},clip,height=0.11\textwidth]{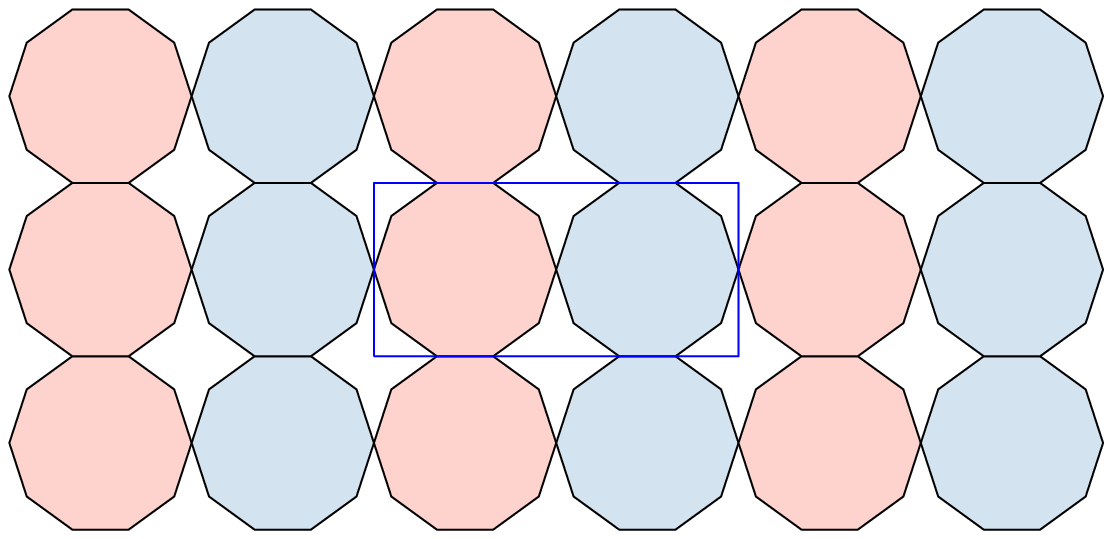}
		\caption{\label{fig:dekagonFig} 		
		Manually constructed $pm$ packings of a (left) pentagon and (right) decagon from their respective densest $p2mm$ packing configurations with approximately equal density as their $p2mm$ counterparts presented in FIG.~\ref{fig:heptagonFig3}.}
	\end{figure}

	For all examined $n$-gons, the densities of the densest $p2mm$ and $pm$ packings are equal. Visually comparing the $p2mm$ and $pm$, the densest packings in FIG.~\ref{fig:heptagonFig3}, the configurations are clearly nonisometric except for $n$-gons with four-fold rotational symmetry. However, the densest $p2mm$ packing of an arbitrary $n$-gon can be easily converted to a $pm$ packing of the same density as is demonstrated in FIG.~\ref{fig:dekagonFig}. This construction suggests that for $n$-gons without four-fold rotational symmetry, the $pm$ density landscape contains at least two nonisometric global maxima and provides multiple densest $pm$ packing configurations. These rewritten $pm$ configurations are isometric to their corresponding densest $p2mm$ packings for $n$-gons with central symmetry.
	
	The difference between these two alternative densest $pm$ configurations for centrally symmetric $n$-gons can be observed by noticing contact edge length. The total length of edges with nonzero contact of the decagon in FIG.~\ref{fig:dekagonFig} with its surrounding decagons is clearly higher when compared to the manually constructed configuration in FIG.~\ref{fig:heptagonFig3}. In the example of the pentagon, the polygons in the $pm$ configuration in FIG.~\ref{fig:heptagonFig3} are slightly rotated compared to FIG.~\ref{fig:dekagonFig} where one of the edges of pentagons is parallel to the basic vector of the primitive cell.
	
Further, our results suggest that for $n$-gons with a central symmetry, the densities of densest $c2mm$, $p2mm$, and $pm$ packings are equal. In fact, there are configurations where all three packing configurations are isometric. Moreover, given a $p2mm$/$pm$, it is possible to construct a $c2mm$ structure with the same density. However, these $c2mm$ packings are not optimal for centrally nonsymmetric $n$-gons.

The densest $p4gm$ packings in our experiments are higher for all centrally nonsymmetric $n$-gons than their corresponding $p2mm/pm$ packings. On the other hand, comparing $p4gm$ and $c2mm$ configurations, $n$-gons with a $4k-1$ fold rotational symmetry attained slightly higher packing densities in the $p4gm$ group than in $c2mm$, and for $n$-gons with a $4k+1$ fold rotational symmetry $c2mm$ packings were higher than those in $c2mm$. For centrally symmetric $n$-gons, the $p4gm$ densest packings were lower or equal to their $c2mm$/$p2mm$/$pm$ densest configurations.

\begin{figure}[t]
		\centering
		\includegraphics[trim={20 0 20 0},clip,width=\linewidth]{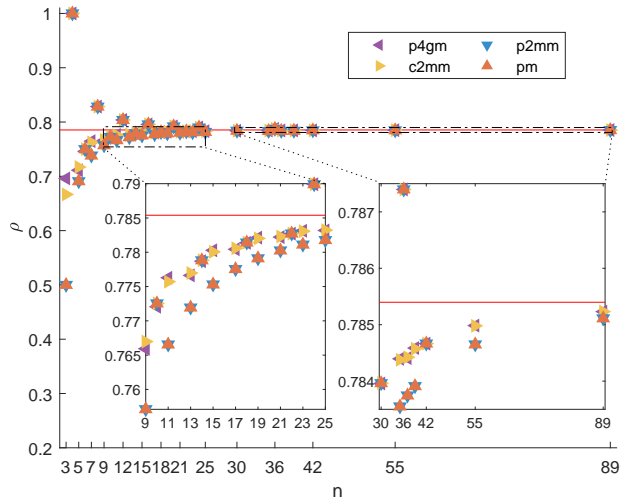}
		\caption{ Densities of the densest packings of examined $n$-gons in the $p4gm / c2mm / p2mm /pm$plane groups class. The red line denotes the density of the densest disc packing in this plane group class.}
		\label{fig:p4gmEvo}
	\end{figure}

The lowest densities of the densest $p4gm$, $c2mm$, $p2mm$, and $pm$ packings were attained in the case of the regular triangle, and uniform square tilings attained the highest densities in all four plane groups. From the evolution of densities as the number of vertices increases, presented in FIG.~\ref{fig:p4gmEvo}, it is reasonable to assume that the extremal values of densities of densest $p4gm$, $c2mm$, $p2mm$, and $pm$ packings among $n$-gons are attained for the regular triangle and square.
	
	\begin{figure*}
		\centering
		
		\begin{minipage}{3.4cm}
			%\vspace{-6cm}
			\centering
			\includegraphics[trim={0 0 0 0},clip,height=2.1cm]{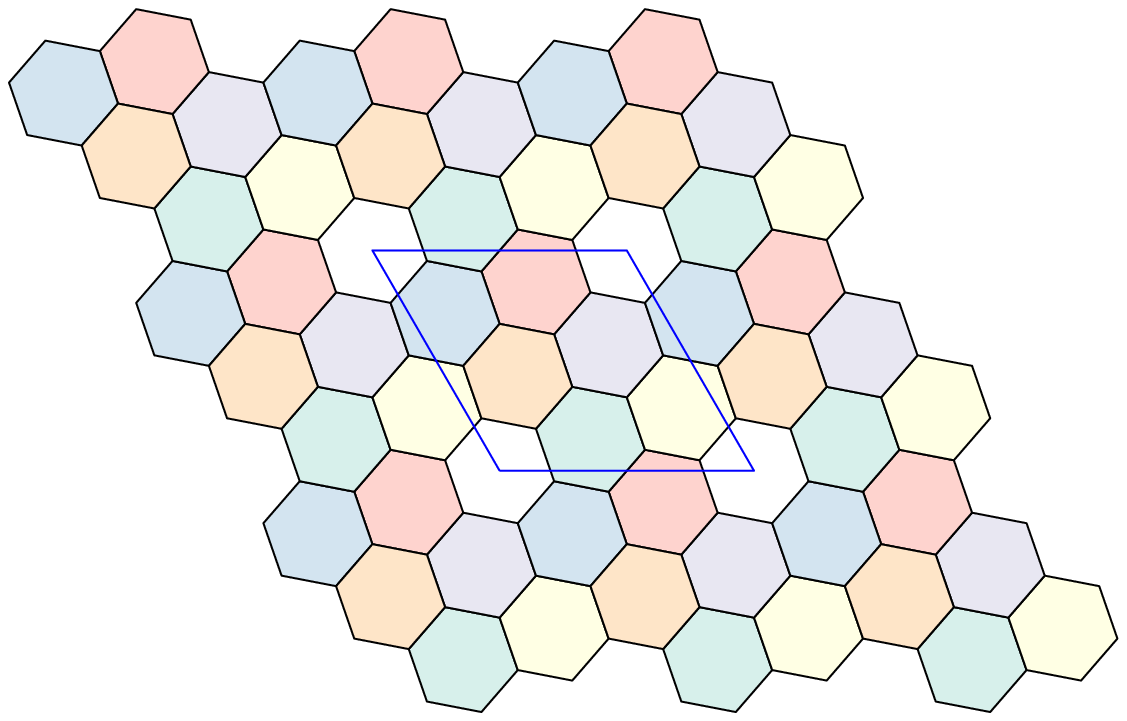}
			
			\includegraphics[trim={0 0 0 0},clip,height=2.1cm]{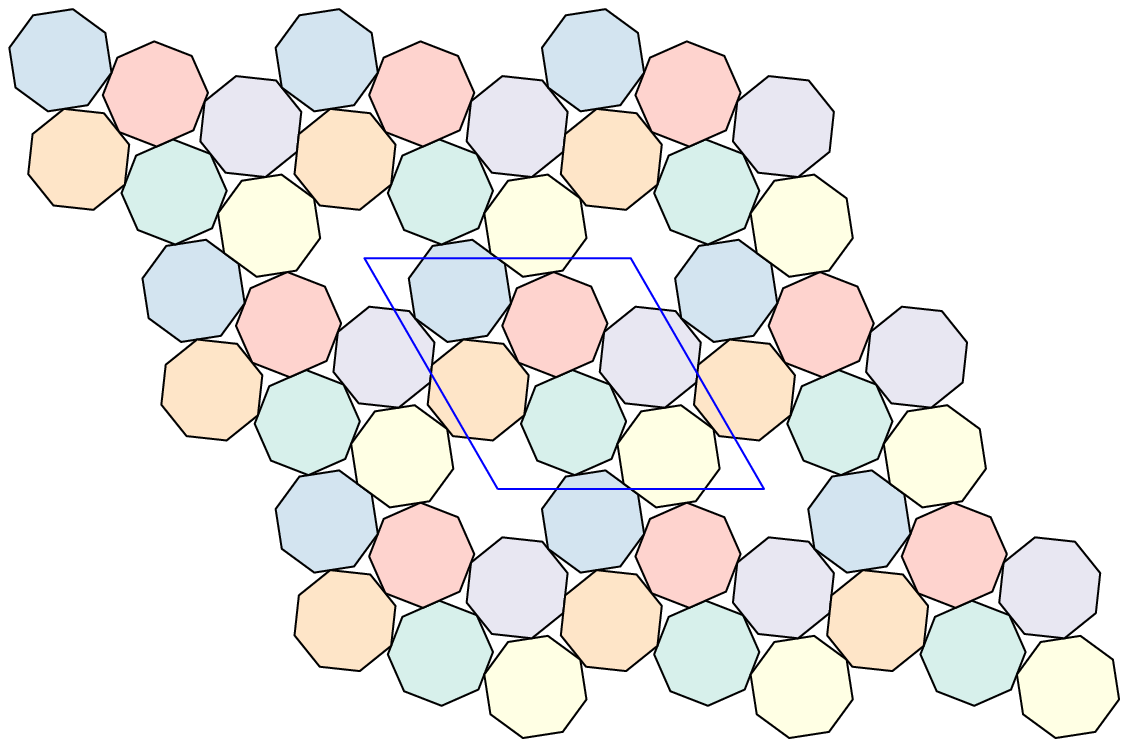}
			
			\includegraphics[trim={0 0 0 0},clip,height=2.1cm]{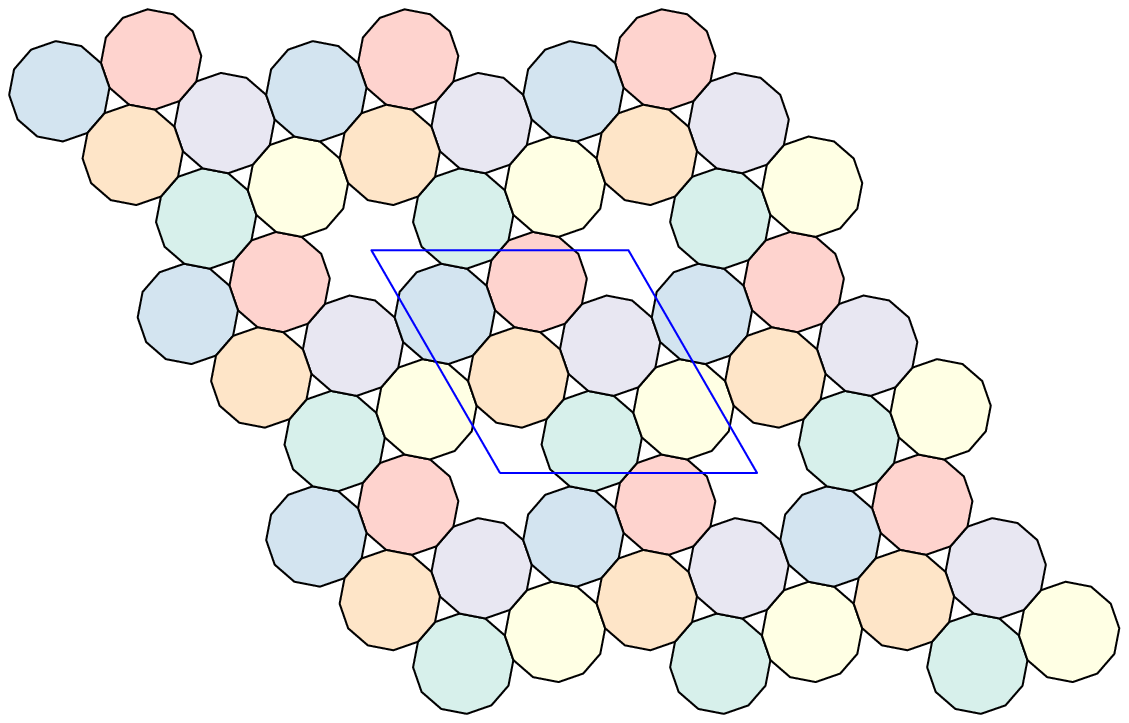}
			%\vspace{-0.5cm}
			
			\subcaption{p6}
			%\vspace{-1.4cm} 
		\end{minipage}	
		\begin{minipage}{3.4cm}
			%\vspace{-6cm}
			\centering
			\includegraphics[trim={0 0 0 0},clip,height=2.1cm]{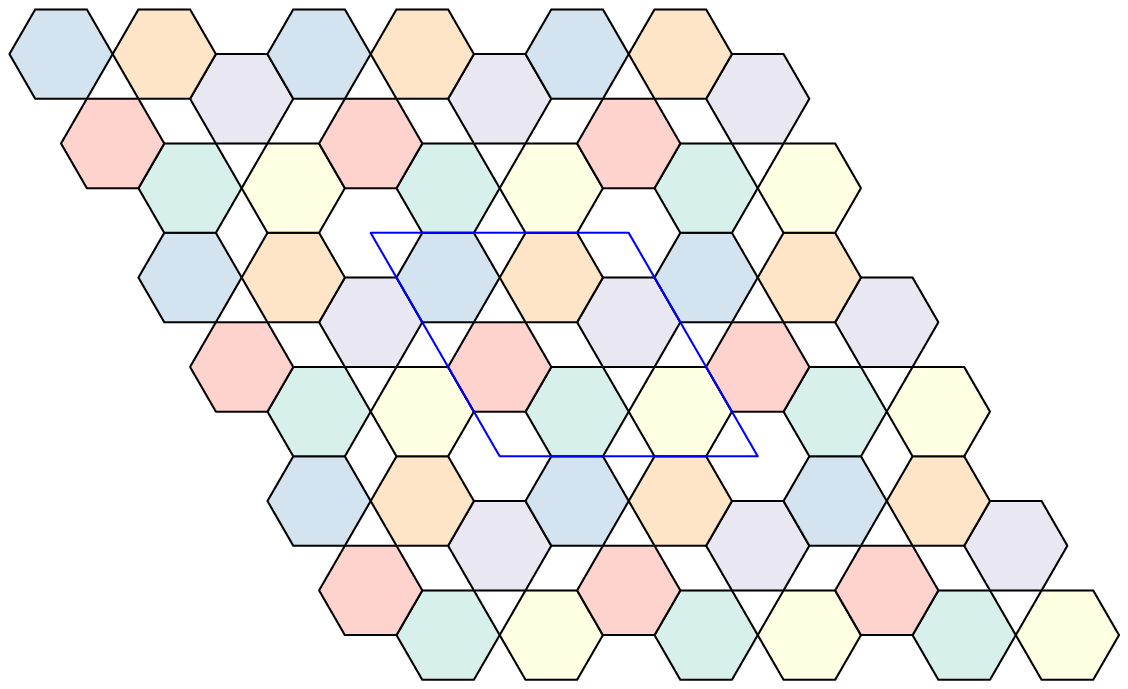}
			
			\includegraphics[trim={0 0 0 0},clip,height=2.1cm]{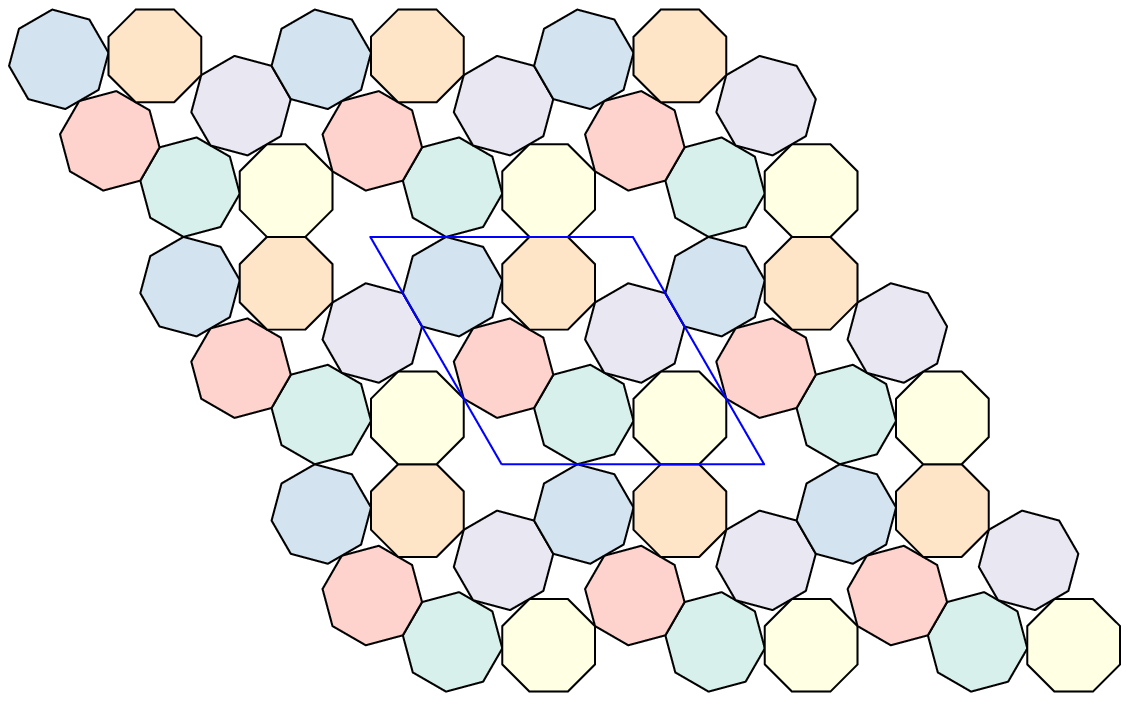}
			
			\includegraphics[trim={0 0 0 0},clip,height=2.1cm]{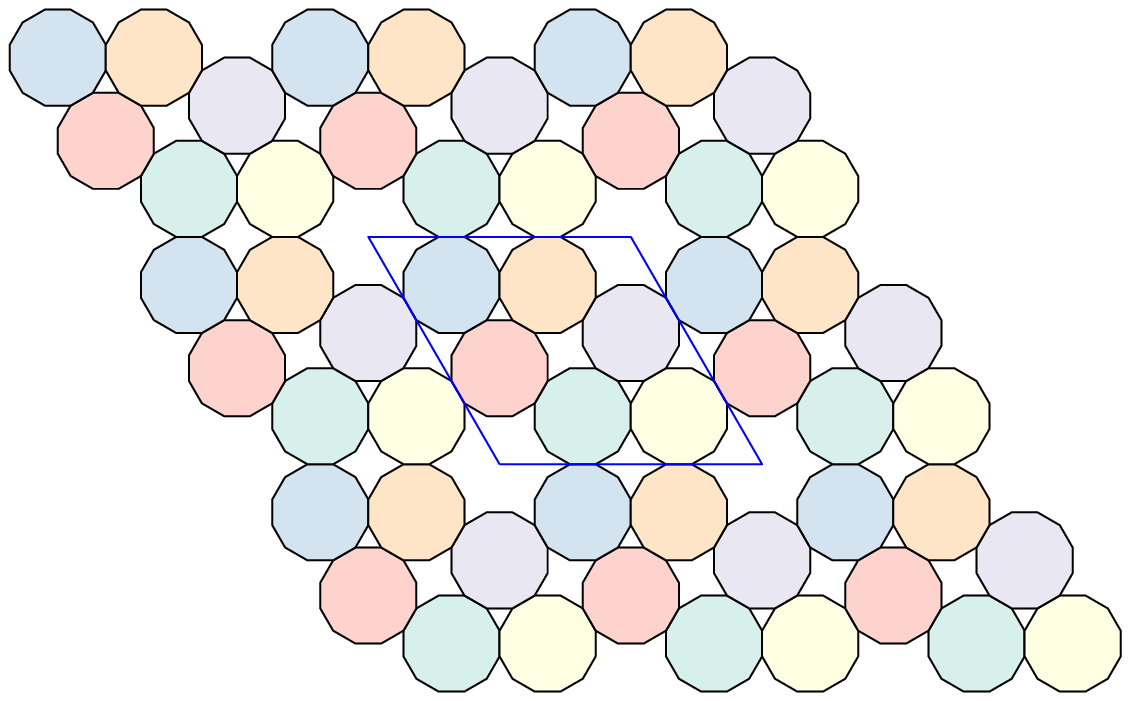}

			\subcaption{p31m}			
		\end{minipage}
		\begin{minipage}{3.4cm}
			%\vspace{-6cm}
			\centering
			\includegraphics[trim={0 0 0 0},clip,height=2.1cm]{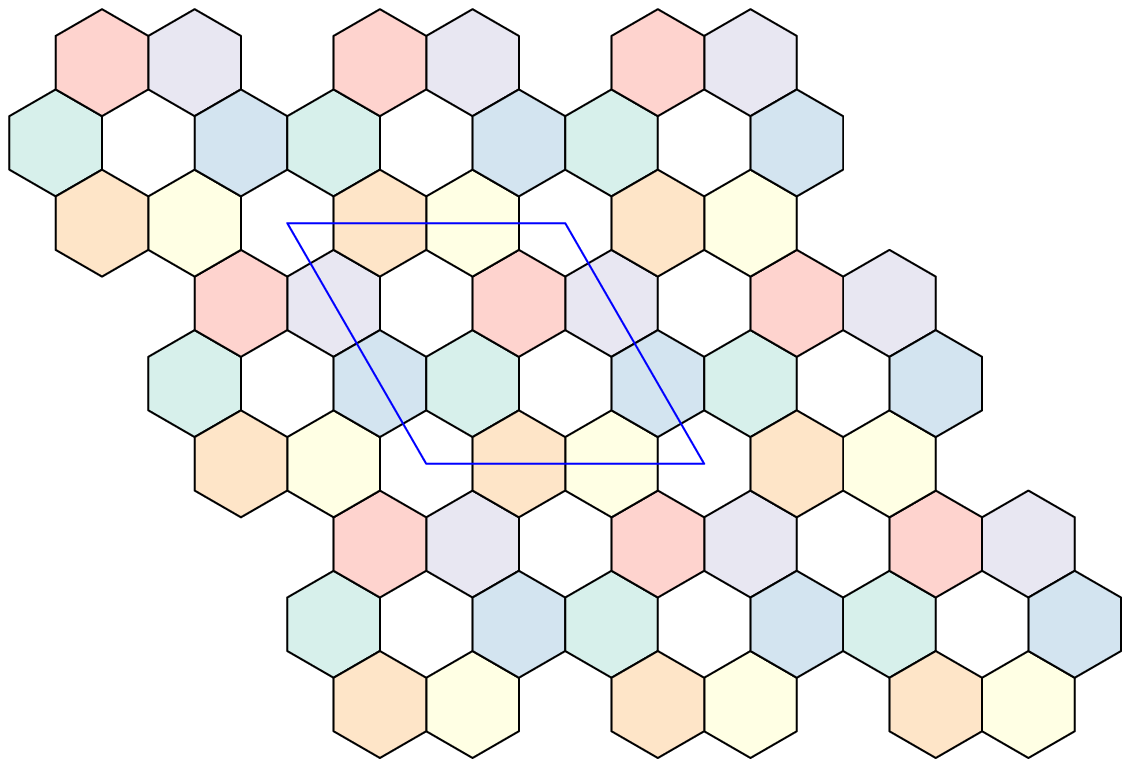}
			
			\includegraphics[trim={0 0 0 0},clip,height=2.1cm]{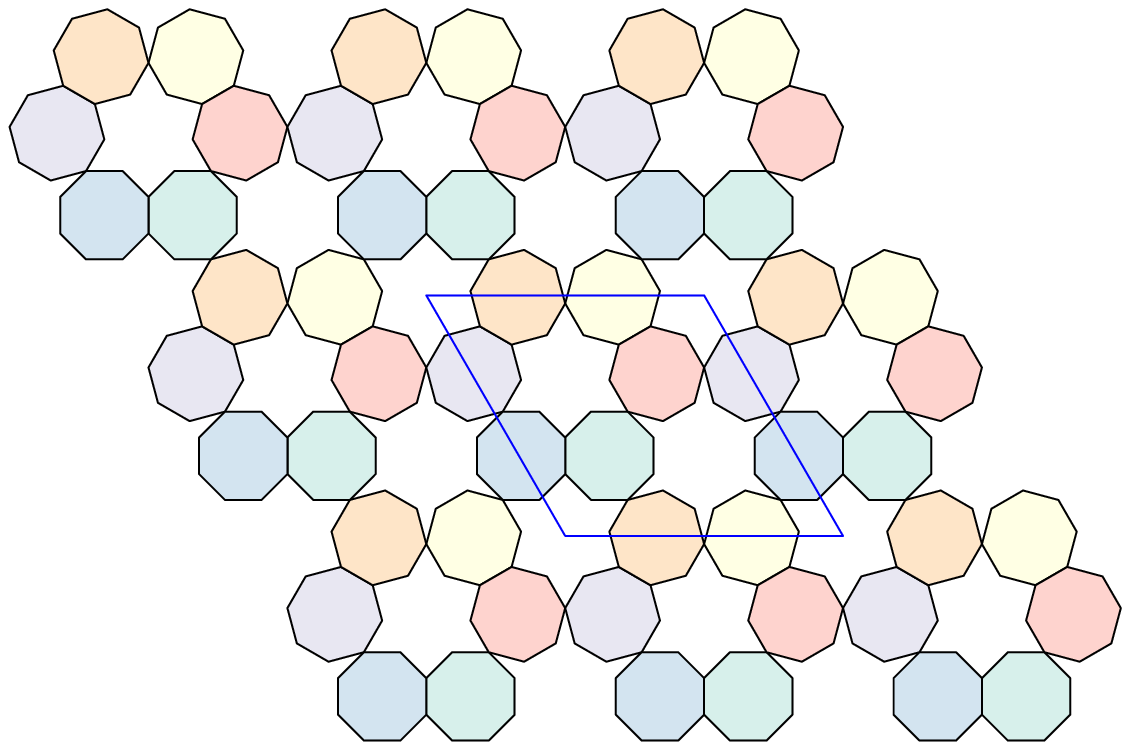}
			
			\includegraphics[trim={0 0 0 0},clip,height=2.1cm]{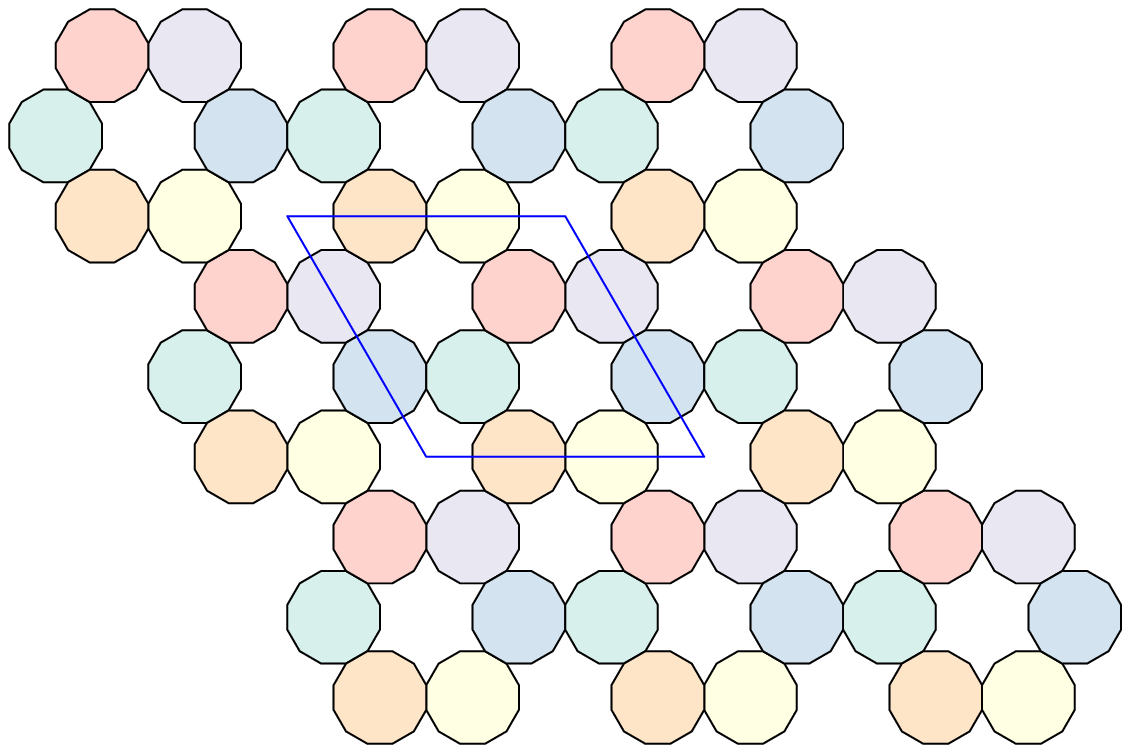}	
			
			\subcaption{p3m1} 
		\end{minipage}
		\begin{minipage}{3.4cm}
			%\vspace{-6cm}
			\centering
			\includegraphics[trim={0 0 0 0},clip,height=2.1cm]{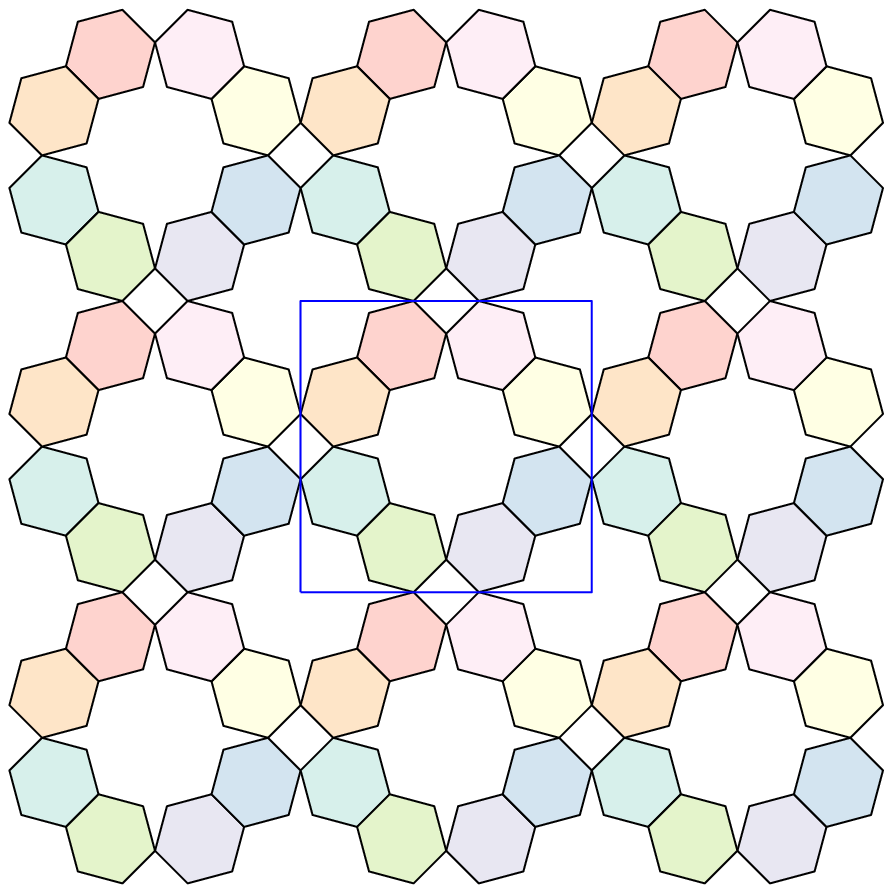}
			
			\includegraphics[trim={0 0 0 0},clip,height=2.1cm]{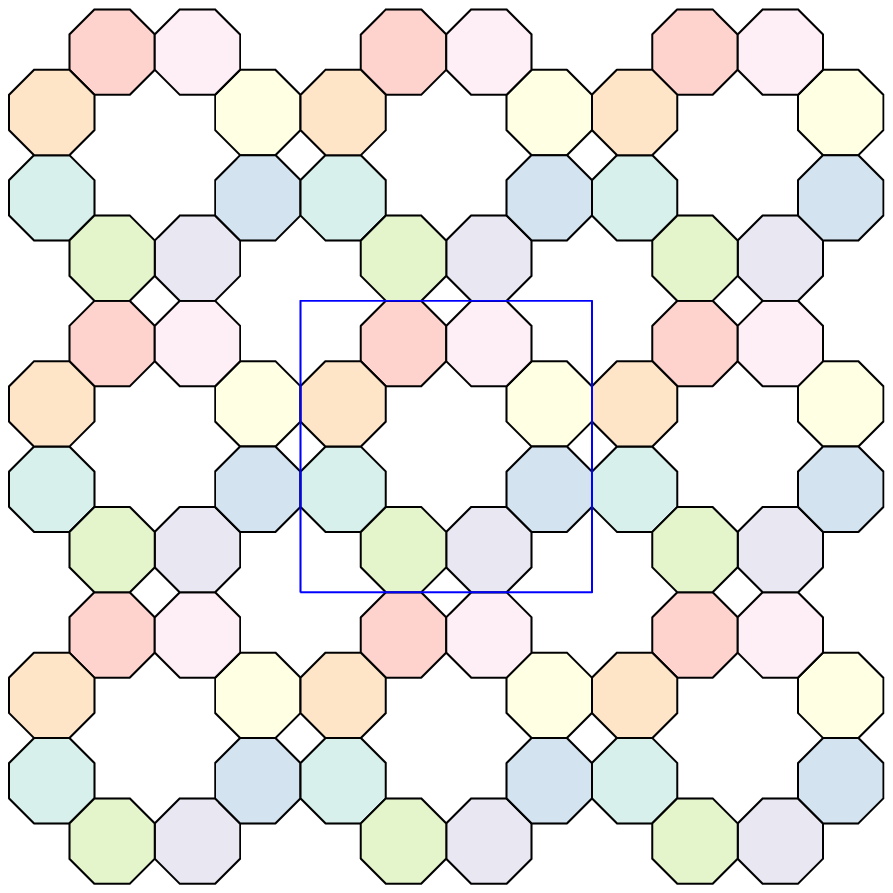}
			
			\includegraphics[trim={0 0 0 0},clip,height=2.1cm]{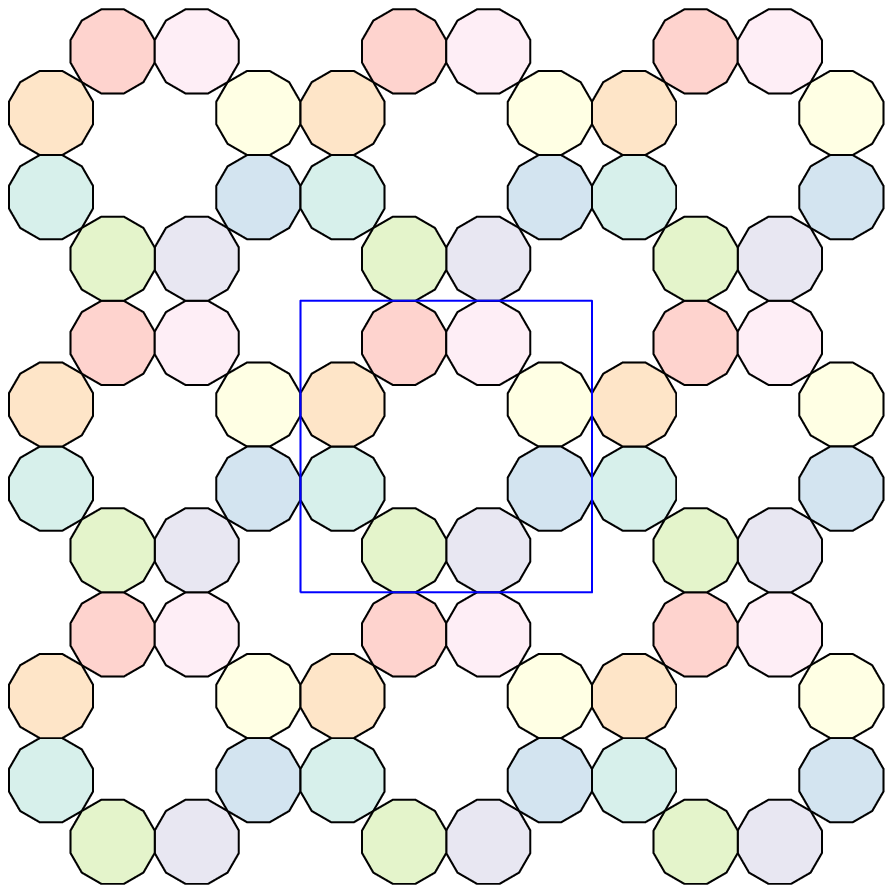}
			
			\subcaption{p4mm} 
		\end{minipage}
		\begin{minipage}{3.4cm}
			%\vspace{-6cm}
			\centering
			\includegraphics[trim={0 0 0 0},clip,height=2.1cm]{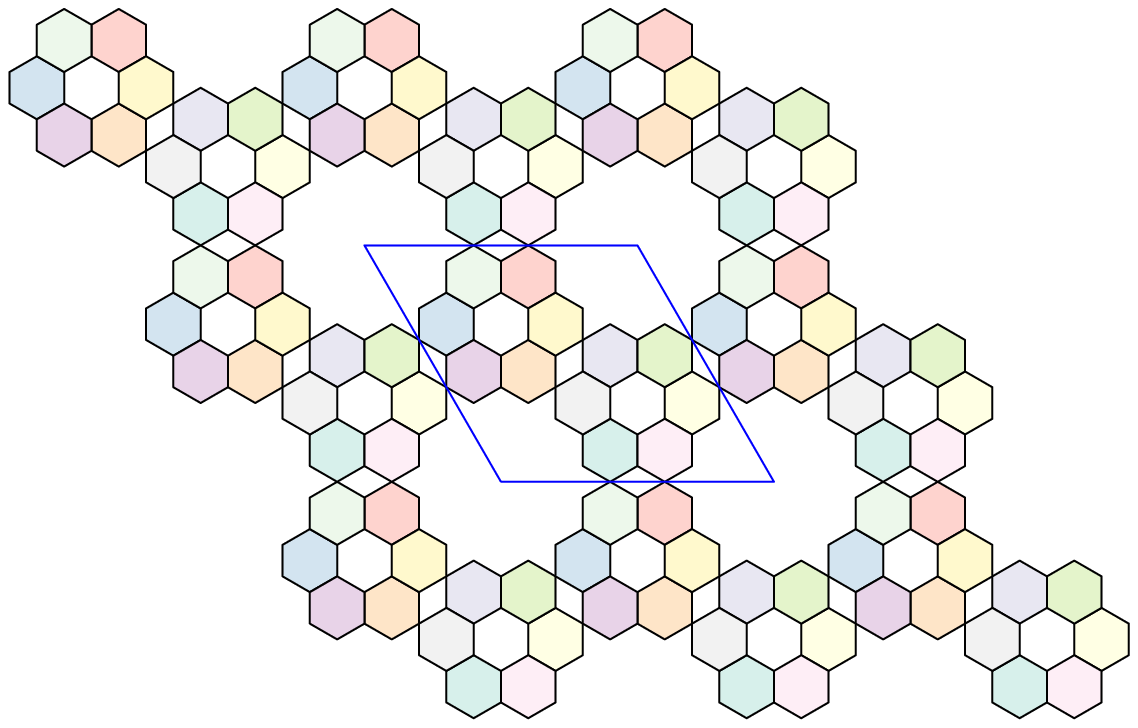}
			
			\includegraphics[trim={0 0 0 0},clip,height=2.1cm]{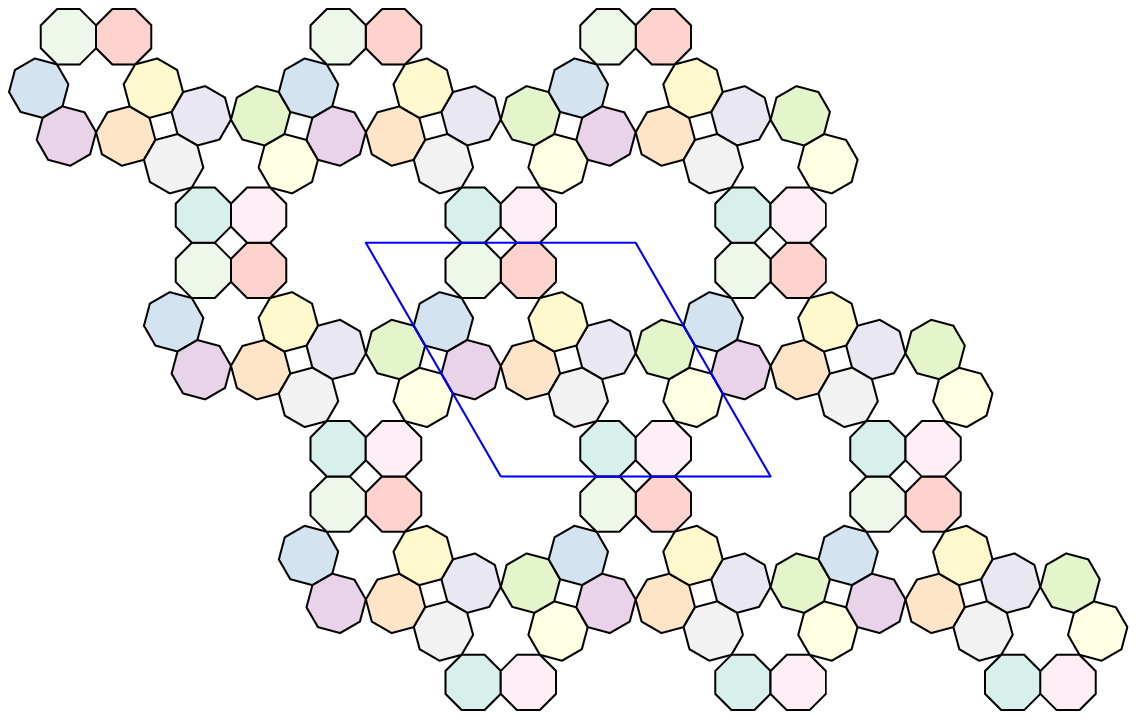}
			
			\includegraphics[trim={0 0 0 0},clip,height=2.1cm]{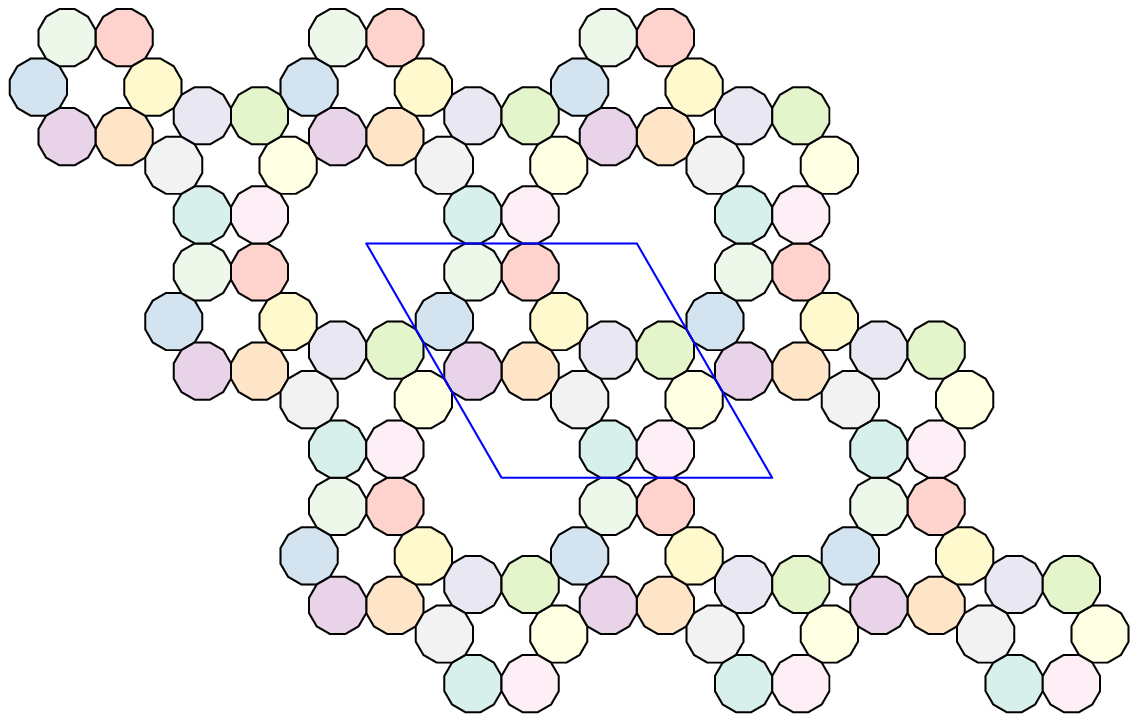}
			
			\subcaption{p6mm} 
		\end{minipage}
		
		\caption{\label{fig:dodekagonP6} Densest configurations of (from top to bottom) hexagon, octagon, and dodecagon in plane groups $p6$, $p31m$, $p3m1$, $p4mm$, and $p6mm$ with the following densities: hexagon in $p6 \approxeq 0.85714$, $p31m \approxeq 0.71999$, $p3m1 \approxeq 0.66666$, $p4mm \approxeq 0.52148$ and $p6mm \approxeq 0.47999$; octagon in $p6 \approxeq 0.76438$, $p31m \approxeq 0.71565$, $p3m1 \approxeq 0.57980$, $p4mm \approxeq 0.56854$ and $p6mm \approxeq 0.48235$; dodecagon in $p6 \approxeq 0.79560$, $p31m \approxeq 0.74613$, $p3m1 \approxeq 0.61880$, $p4mm \approxeq 0.53589$ and $p6mm \approxeq 0.49742$. The blue parallelogram denotes the primitive cell of the respective configuration. Colors represent symmetry operations modulo lattice translations.}
	\end{figure*}

	\subsection{Densest $p6$, $p31m$, $p3m1$, $p4mm$, and $p6mm$ packings}
	\label{sec:p6/p31m/p3m1/p4mm/p6mm}
	
	The last class corresponds to the plane groups for which the optimal densest packings of a disc share no common symmetries. All five plane groups can be therefore regarded as separate classes by themselves. Moreover, the disc attains the lowest density values when the packing configurations are restricted to $p6$, $p31m$, $p3m1$, $p4mm$, and $p6mm$ plane groups.
	
	Similarly to previous classes, we did not observe any clear trend with the packing density oscillating around the respective densest disc packing density. Although the densest $p6$ packing of a disc is lower than the packings in previous plane group classes we examined, this is true for all $n$-gons for $n > 18$, where the distance of densest $n$-gon $p6$ packing density to the $p6$ disc packing is sufficiently small to be separated from the previous class. For instance, $p6$ packing configurations of $5$, $9$, and $18$-gon have higher densities than their corresponding densest packings in the $p4gm/c2mm/p2mm/pm$ class, visually demonstrated in the colored rank table in FIG.~\ref{fig:tabelFig6}.
	
	However, we observed correspondences between the densest plane group packings in this density plane group class and previous classes. For a hexagon, we computed the following ratios of packing densities of respective densest $p2/p2gg/pg/p3/p1$, $p6$, and $p3m1$ packing configurations denoted $\mathcal{K}_{{p2/p2gg/pg/p3/p1}_{\textbf{max}}}$, $\mathcal{K}_{{p6}_{\textbf{max}}}$ and $\mathcal{K}_{{p3m1}_{\textbf{max}}}$,
	\begin{equation*}
		\frac{\rho\left(\mathcal{K}_{{p2/p2gg/pg/p3/p1}_{\textbf{max}}}\right)}{\rho\left(\mathcal{K}_{{p6}_{\textbf{max}}}\right)}=\frac{7}{6}
	\end{equation*}
	and
	\begin{equation*}
		\frac{\rho\left(\mathcal{K}_{{p2/p2gg/pg/p3/p1}_{\textbf{max}}}\right)}{\rho\left(\mathcal{K}_{{p3m1}_{\textbf{max}}}\right)}=\frac{3}{2}.
	\end{equation*}
	By numerical comparison, these ratios approximately hold for all$n$-gons with six-fold rotational symmetry in Table~\ref{tab:table1}. 
	
	Second, by comparing densest $p4mg/c2mm/pm$- $/p2mm$ and $p4mm$ packing configurations of an octagon denoted $\mathcal{K}_{{p4mg/c2mm/pm/p2mm}_{\textbf{max}}}$ and $\mathcal{K}_{{p4mm}_{\textbf{max}}}$, we obtained the following packing density ratio,
	\begin{equation*}
		\frac{\rho\left(\mathcal{K}_{{p4mg/c2mm/pm/p2mm}_{\textbf{max}}}\right)}{\rho\left(\mathcal{K}_{{p4mm}_{\textbf{max}}}\right)} = \frac{3+2\sqrt{2}}{4}.
	\end{equation*}	
	We compared this value against ratios of densest $p4mg/c2mm/pm/p2mm$ and $p4mm$ packings in Table~\ref{tab:table1} and observed an approximate equality for all $n$-gons with eight-fold rotational symmetry. Interestingly, a local eight-rotational symmetry is present in the octagonal octamers of the $p4mm$ packing configuration, shown in FIG.~\ref{fig:dodekagonP6}.
	
	Lastly, 
	we obtained the following packing density ratios of densest $p2mg/cm/p4$, $p31m$, and $p6mm$ packing configurations of a dodecagon denoted $\mathcal{K}_{{p2mg/cm/p4}_{\textbf{max}}}$, $\mathcal{K}_{{p31m}_{\textbf{max}}}$ and $\mathcal{K}_{{p6mm}_{\textbf{max}}}$,
	
	Lastly, we obtained the following packing density ratios of densest $p2mg/cm/p4$, $p31m$, and $p6mm$ packing configurations of a dodecagon denoted $\mathcal{K}_{{p2mg/cm/p4}_{\textbf{max}}}$, $\mathcal{K}_{{p31m}_{\textbf{max}}}$ and $\mathcal{K}_{{p6mm}_{\textbf{max}}}$,
	\begin{equation*}
		\frac{\rho\left(\mathcal{K}_{{p2mg/cm/p4}_{\textbf{max}}}\right)}{\rho\left(\mathcal{K}_{{p31m}_{\textbf{max}}}\right)}
		= \frac{2\sqrt{3}}{3}
	\end{equation*}
	and
	\begin{equation*}
		\frac{\rho\left(\mathcal{K}_{{p2mg/cm/p4}_{\textbf{max}}}\right)}{\rho\left(\mathcal{K}_{{p6mm}_{\textbf{max}}}\right)}
		= \sqrt{3}.
	\end{equation*}
	Compared to the $p2mg/cm/p4$, $p31m$ and $p6mm$ densities in Table~\ref{fig:dodekagonP6}, these packing density ratios are approximately equal for all $n$-gons with $12$-fold rotational symmetry. Additionally, by visual examination of the $p2mg/cm/p4$, $p31m$ and $p6mm$ configurations of a dodecagon in FIG.~\ref{fig:dodekagonP6} and FIG.~\ref{fig:heptagonFig2}, all three configurations contain local three-fold and four-fold rotational symmetries in the dodecagonal trimers and tetramers.
	
	Consequently, these relationships suggest that $24$-fold rotational symmetry of an $n$-gon constitutes the optimal configurations of a disc in all plane groups. Considering that the minimal symmetry containing six-fold, eight-fold, and $12$-fold rotational symmetries is $24$-fold rotational symmetry and that the symmetries of densest $p2/p2gg/pg/p3/p1$, $p2mg/cm/p4$ and $p4mg/c2mm/pm/p2mm$ packings of $n$-gons with $12$-fold rotational symmetry coincide with the symmetries of corresponding densest plane group packings of a disc \footnote[1] . 
	
\begin{figure}[t]
		\centering
		\includegraphics[trim={20 0 20 0},clip,width=\linewidth]{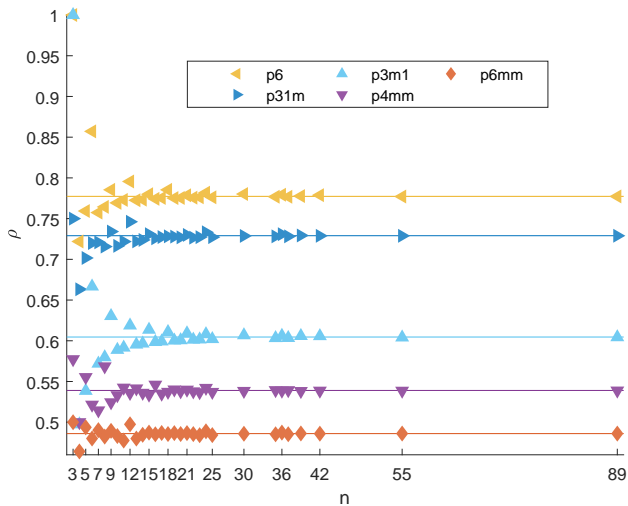}
		\caption{Densities of the densest packings of examined $n$-gons in the $p6 / p31m / p3m1 / p4mm / p6mm$ plane groups class. The colored lines denote the density of the corresponding densest plane group disc packings.}
		\label{fig:p6Evo}
	\end{figure}

	Concerning the extrema of the densest plane group packings, the regular triangle attained the highest packing density in all five plane groups, where in $p6$ and $p3m1$, we have regular triangular tilings. ETRPA obtained the lowest packing densities in our experiments in the case of a square in all five groups. The evolution of packing density as the number of vertices increases, shown in FIG.~\ref{fig:p6Evo}, suggests that the triangle and pentagon attain the extremal densities among all $n$-gons in this plane group class.
	
	\section{Discussion and conclusions }
	\label{sec:conclusions}
	
	Using the ETRPA, we obtained and analyzed the densest packings in all $17$ plane groups of various $n$-gons. Although ETRPA is a stochastic search algorithm and the plane group packing problem has multiple local and global optima, our results indicate that we indeed acquired approximations of densest $n$-gon packing configurations subject to CSG restrictions.
	
	Since two-dimensional CSGs are inherently periodic and combined with the fact that the $p2$ plane group is a local optimum in the space of all packings \cite{kallus2016local}, further confirmed by our results, it can be stated that the $p2$ plane group realizes the densest plane group packings for all $n$-gons among all $17$ plane group. Moreover, depending on the symmetries of an $n$-gon, the densest known packings can be realized in multiple plane groups. In the densest known configuration of a pentagon and heptagon, one of the mirror symmetry planes of the polygon is parallel to one of the primitive cell's basic vectors. This symmetry plane, in fact, coincides with one of the glide reflection planes in the $p2gg$ group and is orthogonal to the $pg$ glide reflection plane. Thus, the densest known configurations can be realized using glide reflections in $p2gg$ and $pg$ plane groups. 
	
	In the densest known packing configuration of the regular enneagon, the mirror symmetry plane of a polygon is slightly rotated in relation to the primitive cell basic vectors \cite{de2011} and, therefore, cannot be constructed using a glide reflection. We have observed this property for all $n$-gons containing a three-fold rotational symmetry but no central symmetry.
	
	\begin{figure}[t]
		\centering
		\includegraphics[trim={70 70 60 60},clip,height=0.3\textwidth]{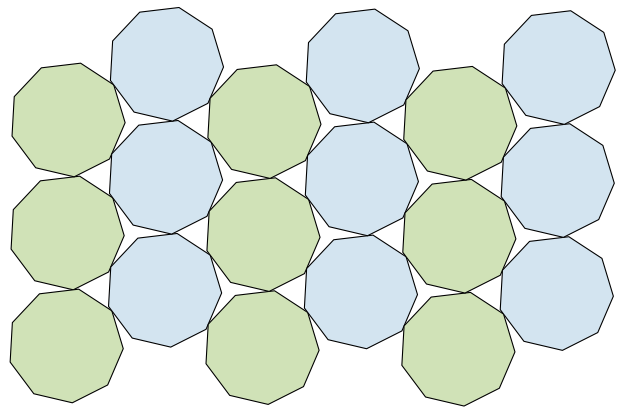}
		
		\caption{\label{fig:packing9gonP2toP1} $p3/p1$ packing of the enneagon with a packing density of approximately $0.87358$, manually constructed from the enneagon's densest $p2$ packing configuration in FIG.~\ref{fig:heptagonFig1}.}
	\end{figure}
	
	One of the methods to construct the densest double-lattice configurations of a convex polygon is based on finding the minimum area of a type of inscribed parallelogram called a half-length parallelogram \cite{mount1991}. In line with our results, the enneagon's diameter given by the minimum area half-length parallelogram has a nonzero slope and the local optimum corresponds to the densest $pg$ packing, contrary to the pentagon and heptagon cases where the diameter corresponding to the global optima coincides with the mirror symmetry of both polygons \cite{kallus2016local}.
	
	An intuition to why the $n$-gons with three-fold rotational symmetry are exceptional can be obtained from the optimal configuration of a disc as an approximation of the enneagon. The densest packing configuration of a disc is also referred to as a hexagonal close-packed configuration and can be constructed using the plane group $p3$ with the hexagonal crystal system. The hexagonal crystal system induces a six-fold rotational symmetry on the crystal lattice with the three-fold rotational symmetry as a subgroup. We observed this six-fold rotational packing symmetry in all $n$-gons with six-fold rotational symmetry.
	
	Additionally, there is a quasi-six-fold rotational packing symmetry \cite{duparcmeur1995dense} in the densest $p2$ configurations of all centrally nonsymmetric $n$-gons with a three-fold rotational symmetry which is not present in densest $p2$ configurations of $n$-gons without three-fold rotational symmetry. This quasi-six-fold rotational symmetry is realized by lattice translations in a $p2$ configuration. Two polygons related to each other by a lattice translation are also related by a three-fold rotational symmetry. Moreover, because of this relation, it is possible to construct a $p3/p1$ packing, as is demonstrated on a packing with $p3$ and $p1$ symmetries constructed from the densest $p2$ packing of an enneagon in FIG.~\ref{fig:packing9gonP2toP1}. Here the blue polygons are unchanged polygons from $p2$ packing of enneagon in FIG.~\ref{fig:heptagonFig1}. However, this packing has lower density than the densest $p3/p1$ packing of an enneagon.
	
	Since the crystallographic restriction theorem does not allow higher than $6$-fold rotational symmetries in a crystal \cite{elliott1998}, our results strongly suggest that centrally non-symmetric $n$-gons with a $3$-fold rotational symmetry are an exception to general densest plane group configuration symmetries of $n$-gons, and it is reasonable to state the following conjectures.
	
	Since the crystallographic restriction theorem does not allow higher than six-fold rotational symmetries in a crystal \cite{elliott1998}, our results strongly suggest that centrally non-symmetric $n$-gons with a three-fold rotational symmetry are an exception to general densest plane group configuration symmetries of $n$-gons, and it is reasonable to state the following conjectures.

	\begin{Conjecture}
		Densities of the densest $p2$, $pg$, and $p2gg$ packings are equal for all, but centrally nonsymmetric $n$-gons with three-fold rotational symmetry and $n\geq 9$, densities of the denses $p2$, $pg$, $p2gg$, and $p1$ packings are equal for all centrally symmetric $n$-gons, and densities of the densest $p2$, $pg$, $p2gg$, $p1$, and $p3$ packings are equal for all $n$-gons containing a six-fold rotational symmetry.
	\end{Conjecture}
	
	\begin{Conjecture}	
		Densities of the densest $p2mg$ and $cm$ packings are equal for all but $n$-gons with a $12k-1$ and $12k+1$ rotational symmetry where $k \in \mathbb{N}$ and densities of densest $p2mg$, $cm$, and $p4$ packings are equal for all $n$-gons containing a $12$-fold rotational symmetry. 
	\end{Conjecture}
	
	\begin{Conjecture}
		Densities  of the densest $pm$ and $p2mm$ packings are equal for all $n$-gons, densities of the densest $c2mm$, $pm$ and $p2mm$ packings are equal for all centrally symmetric $n$-gons, and densities of the densest $p4gm$, $c2mm$, $pm$, and $p2mm$ packings are equal for all $n$-gons containing a $4$-fold rotational symmetry.
		
		Densities of the densest $pm$ and $p2mm$ packings are equal for all $n$-gons, densities of the densest $c2mm$, $pm$ and $p2mm$ packings are equal for all centrally symmetric $n$-gons, and densities of the densest $p4gm$, $c2mm$, $pm$, and $p2mm$ packings are equal for all $n$-gons containing a four-fold rotational symmetry. 
	\end{Conjecture}

	FIG.~\ref{fig:tabelFig6} visualized non-trivial patterns of densest packings of regular
	polygons for the $17$ crystallographic plane groups. The next
	interesting problem is understanding these patterns via group-subgroup
	relations. For example, all five Bravais classes of two-dimensional lattices were
	often studied in a discrete way and only recently were uniquely
	parameterized as subspaces in a common continuous space of all
	lattices up to rigid motion \cite{kurlin2023mathematics,bright2023geographic}.
	
	\begin{acknowledgments}
		The authors thank for funding this research the Leverhulme Trust via the Leverhulme Research Centre for Functional Materials Design, the EPSRC via grants EP/R018472/1 and EP/X018474/1, and the RAEng via fellowship IF2122$\backslash$186.
	\end{acknowledgments}

	\bibliography{PRE_paper}

\begin{thebibliography}{35}%
\makeatletter
\providecommand \@ifxundefined [1]{%
 \@ifx{#1\undefined}
}%
\providecommand \@ifnum [1]{%
 \ifnum #1\expandafter \@firstoftwo
 \else \expandafter \@secondoftwo
 \fi
}%
\providecommand \@ifx [1]{%
 \ifx #1\expandafter \@firstoftwo
 \else \expandafter \@secondoftwo
 \fi
}%
\providecommand \natexlab [1]{#1}%
\providecommand \enquote  [1]{``#1''}%
\providecommand \bibnamefont  [1]{#1}%
\providecommand \bibfnamefont [1]{#1}%
\providecommand \citenamefont [1]{#1}%
\providecommand \href@noop [0]{\@secondoftwo}%
\providecommand \href [0]{\begingroup \@sanitize@url \@href}%
\providecommand \@href[1]{\@@startlink{#1}\@@href}%
\providecommand \@@href[1]{\endgroup#1\@@endlink}%
\providecommand \@sanitize@url [0]{\catcode `\\12\catcode `\$12\catcode
  `\&12\catcode `\#12\catcode `\^12\catcode `\_12\catcode `\%12\relax}%
\providecommand \@@startlink[1]{}%
\providecommand \@@endlink[0]{}%
\providecommand \url  [0]{\begingroup\@sanitize@url \@url }%
\providecommand \@url [1]{\endgroup\@href {#1}{\urlprefix }}%
\providecommand \urlprefix  [0]{URL }%
\providecommand \Eprint [0]{\href }%
\providecommand \doibase [0]{https://doi.org/}%
\providecommand \selectlanguage [0]{\@gobble}%
\providecommand \bibinfo  [0]{\@secondoftwo}%
\providecommand \bibfield  [0]{\@secondoftwo}%
\providecommand \translation [1]{[#1]}%
\providecommand \BibitemOpen [0]{}%
\providecommand \bibitemStop [0]{}%
\providecommand \bibitemNoStop [0]{.\EOS\space}%
\providecommand \EOS [0]{\spacefactor3000\relax}%
\providecommand \BibitemShut  [1]{\csname bibitem#1\endcsname}%
\let\auto@bib@innerbib\@empty
%</preamble>
\bibitem [{\citenamefont {Torquato}(2018)}]{torquato2018}%
  \BibitemOpen
  \bibfield  {author} {\bibinfo {author} {\bibfnamefont {S.}~\bibnamefont
  {Torquato}},\ }\bibfield  {title} {\bibinfo {title} {Perspective: Basic
  understanding of condensed phases of matter via packing models},\ }\href@noop
  {} {\bibfield  {journal} {\bibinfo  {journal} {The Journal of chemical
  physics}\ }\textbf {\bibinfo {volume} {149}},\ \bibinfo {pages} {020901}
  (\bibinfo {year} {2018})}\BibitemShut {NoStop}%
\bibitem [{\citenamefont {Boles}\ \emph {et~al.}(2016)\citenamefont {Boles},
  \citenamefont {Engel},\ and\ \citenamefont {Talapin}}]{boles2016}%
  \BibitemOpen
  \bibfield  {author} {\bibinfo {author} {\bibfnamefont {M.~A.}\ \bibnamefont
  {Boles}}, \bibinfo {author} {\bibfnamefont {M.}~\bibnamefont {Engel}},\ and\
  \bibinfo {author} {\bibfnamefont {D.~V.}\ \bibnamefont {Talapin}},\
  }\bibfield  {title} {\bibinfo {title} {Self-assembly of colloidal
  nanocrystals: From intricate structures to functional materials},\
  }\href@noop {} {\bibfield  {journal} {\bibinfo  {journal} {Chemical reviews}\
  }\textbf {\bibinfo {volume} {116}},\ \bibinfo {pages} {11220} (\bibinfo
  {year} {2016})}\BibitemShut {NoStop}%
\bibitem [{\citenamefont {S{\'a}nchez-Gonz{\'a}lez}\ \emph
  {et~al.}(2022)\citenamefont {S{\'a}nchez-Gonz{\'a}lez}, \citenamefont
  {Tsang}, \citenamefont {Troyano}, \citenamefont {Craig},\ and\ \citenamefont
  {Furukawa}}]{sanchez2022}%
  \BibitemOpen
  \bibfield  {author} {\bibinfo {author} {\bibfnamefont {E.}~\bibnamefont
  {S{\'a}nchez-Gonz{\'a}lez}}, \bibinfo {author} {\bibfnamefont {M.~Y.}\
  \bibnamefont {Tsang}}, \bibinfo {author} {\bibfnamefont {J.}~\bibnamefont
  {Troyano}}, \bibinfo {author} {\bibfnamefont {G.~A.}\ \bibnamefont {Craig}},\
  and\ \bibinfo {author} {\bibfnamefont {S.}~\bibnamefont {Furukawa}},\
  }\bibfield  {title} {\bibinfo {title} {Assembling metal--organic cages as
  porous materials},\ }\href@noop {} {\bibfield  {journal} {\bibinfo  {journal}
  {Chemical Society Reviews}\ } (\bibinfo {year} {2022})}\BibitemShut {NoStop}%
\bibitem [{\citenamefont {Liang}\ and\ \citenamefont {Dill}(2001)}]{liang2001}%
  \BibitemOpen
  \bibfield  {author} {\bibinfo {author} {\bibfnamefont {J.}~\bibnamefont
  {Liang}}\ and\ \bibinfo {author} {\bibfnamefont {K.~A.}\ \bibnamefont
  {Dill}},\ }\bibfield  {title} {\bibinfo {title} {Are proteins well-packed?},\
  }\href@noop {} {\bibfield  {journal} {\bibinfo  {journal} {Biophysical
  journal}\ }\textbf {\bibinfo {volume} {81}},\ \bibinfo {pages} {751}
  (\bibinfo {year} {2001})}\BibitemShut {NoStop}%
\bibitem [{\citenamefont {Tarnai}\ \emph {et~al.}(1995)\citenamefont {Tarnai},
  \citenamefont {Gaspar},\ and\ \citenamefont {Szalai}}]{tarnai1995}%
  \BibitemOpen
  \bibfield  {author} {\bibinfo {author} {\bibfnamefont {T.}~\bibnamefont
  {Tarnai}}, \bibinfo {author} {\bibfnamefont {Z.}~\bibnamefont {Gaspar}},\
  and\ \bibinfo {author} {\bibfnamefont {L.}~\bibnamefont {Szalai}},\
  }\bibfield  {title} {\bibinfo {title} {Pentagon packing models for"
  all-pentamer" virus structures},\ }\href@noop {} {\bibfield  {journal}
  {\bibinfo  {journal} {Biophysical journal}\ }\textbf {\bibinfo {volume}
  {69}},\ \bibinfo {pages} {612} (\bibinfo {year} {1995})}\BibitemShut
  {NoStop}%
\bibitem [{\citenamefont {Zoppi}\ \emph {et~al.}(2012)\citenamefont {Zoppi},
  \citenamefont {Bauert}, \citenamefont {Siegel}, \citenamefont {Baldridge},\
  and\ \citenamefont {Ernst}}]{zoppi2012}%
  \BibitemOpen
  \bibfield  {author} {\bibinfo {author} {\bibfnamefont {L.}~\bibnamefont
  {Zoppi}}, \bibinfo {author} {\bibfnamefont {T.}~\bibnamefont {Bauert}},
  \bibinfo {author} {\bibfnamefont {J.~S.}\ \bibnamefont {Siegel}}, \bibinfo
  {author} {\bibfnamefont {K.~K.}\ \bibnamefont {Baldridge}},\ and\ \bibinfo
  {author} {\bibfnamefont {K.-H.}\ \bibnamefont {Ernst}},\ }\bibfield  {title}
  {\bibinfo {title} {Pentagonal tiling with buckybowls: pentamethylcorannulene
  on cu (111)},\ }\href@noop {} {\bibfield  {journal} {\bibinfo  {journal}
  {Physical Chemistry Chemical Physics}\ }\textbf {\bibinfo {volume} {14}},\
  \bibinfo {pages} {13365} (\bibinfo {year} {2012})}\BibitemShut {NoStop}%
\bibitem [{\citenamefont {Barth}(2007)}]{barth2007}%
  \BibitemOpen
  \bibfield  {author} {\bibinfo {author} {\bibfnamefont {J.~V.}\ \bibnamefont
  {Barth}},\ }\bibfield  {title} {\bibinfo {title} {Molecular architectonic on
  metal surfaces},\ }\href@noop {} {\bibfield  {journal} {\bibinfo  {journal}
  {Annu. Rev. Phys. Chem.}\ }\textbf {\bibinfo {volume} {58}},\ \bibinfo
  {pages} {375} (\bibinfo {year} {2007})}\BibitemShut {NoStop}%
\bibitem [{\citenamefont {Cui}\ \emph {et~al.}(2018)\citenamefont {Cui},
  \citenamefont {Ebrahimi}, \citenamefont {Macleod},\ and\ \citenamefont
  {Rosei}}]{cui2018}%
  \BibitemOpen
  \bibfield  {author} {\bibinfo {author} {\bibfnamefont {D.}~\bibnamefont
  {Cui}}, \bibinfo {author} {\bibfnamefont {M.}~\bibnamefont {Ebrahimi}},
  \bibinfo {author} {\bibfnamefont {J.~M.}\ \bibnamefont {Macleod}},\ and\
  \bibinfo {author} {\bibfnamefont {F.}~\bibnamefont {Rosei}},\ }\bibfield
  {title} {\bibinfo {title} {Template-driven dense packing of pentagonal
  molecules in monolayer films},\ }\href@noop {} {\bibfield  {journal}
  {\bibinfo  {journal} {Nano Letters}\ }\textbf {\bibinfo {volume} {18}},\
  \bibinfo {pages} {7570} (\bibinfo {year} {2018})}\BibitemShut {NoStop}%
\bibitem [{\citenamefont {Wang}\ and\ \citenamefont {Sui}(1999)}]{wang1999}%
  \BibitemOpen
  \bibfield  {author} {\bibinfo {author} {\bibfnamefont {H.-W.}\ \bibnamefont
  {Wang}}\ and\ \bibinfo {author} {\bibfnamefont {S.-f.}\ \bibnamefont {Sui}},\
  }\bibfield  {title} {\bibinfo {title} {Pentameric two-dimensional
  crystallization of rabbit c-reactive protein on lipid monolayers},\
  }\href@noop {} {\bibfield  {journal} {\bibinfo  {journal} {Journal of
  structural biology}\ }\textbf {\bibinfo {volume} {127}},\ \bibinfo {pages}
  {283} (\bibinfo {year} {1999})}\BibitemShut {NoStop}%
\bibitem [{\citenamefont {Theil}(2006)}]{theil2006}%
  \BibitemOpen
  \bibfield  {author} {\bibinfo {author} {\bibfnamefont {F.}~\bibnamefont
  {Theil}},\ }\bibfield  {title} {\bibinfo {title} {A proof of crystallization
  in two dimensions},\ }\href@noop {} {\bibfield  {journal} {\bibinfo
  {journal} {Communications in Mathematical Physics}\ }\textbf {\bibinfo
  {volume} {262}},\ \bibinfo {pages} {209} (\bibinfo {year}
  {2006})}\BibitemShut {NoStop}%
\bibitem [{\citenamefont {Weinan}\ and\ \citenamefont {Li}(2009)}]{weinan2009}%
  \BibitemOpen
  \bibfield  {author} {\bibinfo {author} {\bibfnamefont {E.}~\bibnamefont
  {Weinan}}\ and\ \bibinfo {author} {\bibfnamefont {D.}~\bibnamefont {Li}},\
  }\bibfield  {title} {\bibinfo {title} {On the crystallization of 2d hexagonal
  lattices},\ }\href@noop {} {\bibfield  {journal} {\bibinfo  {journal}
  {Communications in mathematical physics}\ }\textbf {\bibinfo {volume}
  {286}},\ \bibinfo {pages} {1099} (\bibinfo {year} {2009})}\BibitemShut
  {NoStop}%
\bibitem [{\citenamefont {Woodley}\ \emph {et~al.}(2020)\citenamefont
  {Woodley}, \citenamefont {Day},\ and\ \citenamefont {Catlow}}]{woodley2020}%
  \BibitemOpen
  \bibfield  {author} {\bibinfo {author} {\bibfnamefont {S.~M.}\ \bibnamefont
  {Woodley}}, \bibinfo {author} {\bibfnamefont {G.~M.}\ \bibnamefont {Day}},\
  and\ \bibinfo {author} {\bibfnamefont {R.}~\bibnamefont {Catlow}},\
  }\bibfield  {title} {\bibinfo {title} {Structure prediction of crystals,
  surfaces and nanoparticles},\ }\href@noop {} {\bibfield  {journal} {\bibinfo
  {journal} {Philosophical Transactions of the Royal Society A}\ }\textbf
  {\bibinfo {volume} {378}},\ \bibinfo {pages} {20190600} (\bibinfo {year}
  {2020})}\BibitemShut {NoStop}%
\bibitem [{\citenamefont {Courant}(1965)}]{courant1965least}%
  \BibitemOpen
  \bibfield  {author} {\bibinfo {author} {\bibfnamefont {R.}~\bibnamefont
  {Courant}},\ }\bibfield  {title} {\bibinfo {title} {The least dense lattice
  packing of two-dimensional convex bodies},\ }\href@noop {} {\bibfield
  {journal} {\bibinfo  {journal} {Communications on Pure and Applied
  Mathematics}\ }\textbf {\bibinfo {volume} {18}},\ \bibinfo {pages} {339}
  (\bibinfo {year} {1965})}\BibitemShut {NoStop}%
\bibitem [{\citenamefont {T{\'o}th}(2013)}]{toth2013}%
  \BibitemOpen
  \bibfield  {author} {\bibinfo {author} {\bibfnamefont {L.~F.}\ \bibnamefont
  {T{\'o}th}},\ }\href@noop {} {\emph {\bibinfo {title} {Lagerungen in der
  Ebene auf der Kugel und im Raum}}},\ Vol.~\bibinfo {volume} {65}\ (\bibinfo
  {publisher} {Springer-Verlag},\ \bibinfo {year} {2013})\BibitemShut {NoStop}%
\bibitem [{\citenamefont {Kuperberg}\ and\ \citenamefont
  {Kuperberg}(1990)}]{kuperbergSup}%
  \BibitemOpen
  \bibfield  {author} {\bibinfo {author} {\bibfnamefont {G.}~\bibnamefont
  {Kuperberg}}\ and\ \bibinfo {author} {\bibfnamefont {W.}~\bibnamefont
  {Kuperberg}},\ }\bibfield  {title} {\bibinfo {title} {Double-lattice packings
  of convex bodies in the plane},\ }\href@noop {} {\bibfield  {journal}
  {\bibinfo  {journal} {Discrete \& Computational Geometry}\ }\textbf {\bibinfo
  {volume} {5}},\ \bibinfo {pages} {389} (\bibinfo {year} {1990})}\BibitemShut
  {NoStop}%
\bibitem [{\citenamefont {Hales}\ and\ \citenamefont {Kusner}(2016)}]{hales}%
  \BibitemOpen
  \bibfield  {author} {\bibinfo {author} {\bibfnamefont {T.}~\bibnamefont
  {Hales}}\ and\ \bibinfo {author} {\bibfnamefont {W.}~\bibnamefont {Kusner}},\
  }\bibfield  {title} {\bibinfo {title} {Packings of regular pentagons in the
  plane},\ }\href@noop {} {\bibfield  {journal} {\bibinfo  {journal} {arXiv
  preprint arXiv:1602.07220}\ } (\bibinfo {year} {2016})}\BibitemShut {NoStop}%
\bibitem [{\citenamefont {Torda}\ \emph {et~al.}(2022)\citenamefont {Torda},
  \citenamefont {Goulermas}, \citenamefont {P{\'u}{\v{c}}ek},\ and\
  \citenamefont {Kurlin}}]{torda2022entropic}%
  \BibitemOpen
  \bibfield  {author} {\bibinfo {author} {\bibfnamefont {M.}~\bibnamefont
  {Torda}}, \bibinfo {author} {\bibfnamefont {J.~Y.}\ \bibnamefont
  {Goulermas}}, \bibinfo {author} {\bibfnamefont {R.}~\bibnamefont
  {P{\'u}{\v{c}}ek}},\ and\ \bibinfo {author} {\bibfnamefont {V.}~\bibnamefont
  {Kurlin}},\ }\bibfield  {title} {\bibinfo {title} {Entropic trust region for
  densest crystallographic symmetry group packings},\ }\href@noop {} {\bibfield
   {journal} {\bibinfo  {journal} {arXiv preprint arXiv:2202.11959}\ }
  (\bibinfo {year} {2022})}\BibitemShut {NoStop}%
\bibitem [{\citenamefont {Hahn}(2005)}]{brock2016}%
  \BibitemOpen
  \bibinfo {editor} {\bibfnamefont {T.}~\bibnamefont {Hahn}},\ ed.,\ \href@noop
  {} {\emph {\bibinfo {title} {International tables for crystallography volume
  A: Space-group symmetry}}},\ \bibinfo {edition} {5th}\ ed.\ (\bibinfo
  {publisher} {Springer},\ \bibinfo {year} {2005})\BibitemShut {NoStop}%
\bibitem [{\citenamefont {Rogers}(1964)}]{rogers}%
  \BibitemOpen
  \bibfield  {author} {\bibinfo {author} {\bibfnamefont {C.~A.}\ \bibnamefont
  {Rogers}},\ }\href@noop {} {\emph {\bibinfo {title} {Packing and
  covering}}},\ \bibinfo {number} {54}\ (\bibinfo  {publisher} {Cambridge
  University Press},\ \bibinfo {year} {1964})\BibitemShut {NoStop}%
\bibitem [{\citenamefont {Geman}\ and\ \citenamefont
  {Geman}(1984)}]{geman1984}%
  \BibitemOpen
  \bibfield  {author} {\bibinfo {author} {\bibfnamefont {S.}~\bibnamefont
  {Geman}}\ and\ \bibinfo {author} {\bibfnamefont {D.}~\bibnamefont {Geman}},\
  }\bibfield  {title} {\bibinfo {title} {Stochastic relaxation, gibbs
  distributions, and the bayesian restoration of images},\ }\href@noop {}
  {\bibfield  {journal} {\bibinfo  {journal} {IEEE Transactions on pattern
  analysis and machine intelligence}\ ,\ \bibinfo {pages} {721}} (\bibinfo
  {year} {1984})}\BibitemShut {NoStop}%
\bibitem [{\citenamefont {Mardia}\ \emph {et~al.}(2008)\citenamefont {Mardia},
  \citenamefont {Hughes}, \citenamefont {Taylor},\ and\ \citenamefont
  {Singh}}]{mardia2008}%
  \BibitemOpen
  \bibfield  {author} {\bibinfo {author} {\bibfnamefont {K.~V.}\ \bibnamefont
  {Mardia}}, \bibinfo {author} {\bibfnamefont {G.}~\bibnamefont {Hughes}},
  \bibinfo {author} {\bibfnamefont {C.~C.}\ \bibnamefont {Taylor}},\ and\
  \bibinfo {author} {\bibfnamefont {H.}~\bibnamefont {Singh}},\ }\bibfield
  {title} {\bibinfo {title} {A multivariate von mises distribution with
  applications to bioinformatics},\ }\href@noop {} {\bibfield  {journal}
  {\bibinfo  {journal} {Canadian Journal of Statistics}\ }\textbf {\bibinfo
  {volume} {36}},\ \bibinfo {pages} {99} (\bibinfo {year} {2008})}\BibitemShut
  {NoStop}%
\bibitem [{\citenamefont {Kullback}\ and\ \citenamefont
  {Leibler}(1951)}]{kullback1951}%
  \BibitemOpen
  \bibfield  {author} {\bibinfo {author} {\bibfnamefont {S.}~\bibnamefont
  {Kullback}}\ and\ \bibinfo {author} {\bibfnamefont {R.~A.}\ \bibnamefont
  {Leibler}},\ }\bibfield  {title} {\bibinfo {title} {On information and
  sufficiency},\ }\href@noop {} {\bibfield  {journal} {\bibinfo  {journal} {The
  annals of mathematical statistics}\ }\textbf {\bibinfo {volume} {22}},\
  \bibinfo {pages} {79} (\bibinfo {year} {1951})}\BibitemShut {NoStop}%
\bibitem [{\citenamefont {Amari}(1998)}]{amari}%
  \BibitemOpen
  \bibfield  {author} {\bibinfo {author} {\bibfnamefont {S.-I.}\ \bibnamefont
  {Amari}},\ }\bibfield  {title} {\bibinfo {title} {Natural gradient works
  efficiently in learning},\ }\href@noop {} {\bibfield  {journal} {\bibinfo
  {journal} {Neural computation}\ }\textbf {\bibinfo {volume} {10}},\ \bibinfo
  {pages} {251} (\bibinfo {year} {1998})}\BibitemShut {NoStop}%
\bibitem [{Note1()}]{Note1}%
  \BibitemOpen
  \bibinfo {note} {See Supplemental Material at \protect \href
  {https://milotorda.net/wp-content/uploads/supplemental_material.pdf}{milotorda.net}
  for a complete table of densest plane group packings of $n$-gons, including
  configuration parameters and visualizations of respective structures and for
  exact densities of a discs densest plane group packings and additional
  information on relationships between densest plane group packings of a disc
  and $n$-gons with $24$-fold rotational symmetry}\BibitemShut {NoStop}%
\bibitem [{\citenamefont {Elliott}(1998)}]{elliott1998}%
  \BibitemOpen
  \bibfield  {author} {\bibinfo {author} {\bibfnamefont {S.}~\bibnamefont
  {Elliott}},\ }\href@noop {} {\emph {\bibinfo {title} {The physics and
  chemistry of solids}}}\ (\bibinfo  {publisher} {Wiley},\ \bibinfo {year}
  {1998})\BibitemShut {NoStop}%
\bibitem [{\citenamefont {Atkinson}\ \emph {et~al.}(2012)\citenamefont
  {Atkinson}, \citenamefont {Jiao},\ and\ \citenamefont {Torquato}}]{atkinson}%
  \BibitemOpen
  \bibfield  {author} {\bibinfo {author} {\bibfnamefont {S.}~\bibnamefont
  {Atkinson}}, \bibinfo {author} {\bibfnamefont {Y.}~\bibnamefont {Jiao}},\
  and\ \bibinfo {author} {\bibfnamefont {S.}~\bibnamefont {Torquato}},\
  }\bibfield  {title} {\bibinfo {title} {Maximally dense packings of
  two-dimensional convex and concave noncircular particles},\ }\href@noop {}
  {\bibfield  {journal} {\bibinfo  {journal} {Physical Review E}\ }\textbf
  {\bibinfo {volume} {86}},\ \bibinfo {pages} {031302} (\bibinfo {year}
  {2012})}\BibitemShut {NoStop}%
\bibitem [{\citenamefont {de~Graaf}\ \emph {et~al.}(2011)\citenamefont
  {de~Graaf}, \citenamefont {van Roij},\ and\ \citenamefont
  {Dijkstra}}]{de2011}%
  \BibitemOpen
  \bibfield  {author} {\bibinfo {author} {\bibfnamefont {J.}~\bibnamefont
  {de~Graaf}}, \bibinfo {author} {\bibfnamefont {R.}~\bibnamefont {van Roij}},\
  and\ \bibinfo {author} {\bibfnamefont {M.}~\bibnamefont {Dijkstra}},\
  }\bibfield  {title} {\bibinfo {title} {Dense regular packings of irregular
  nonconvex particles},\ }\href@noop {} {\bibfield  {journal} {\bibinfo
  {journal} {Physical Review Letters}\ }\textbf {\bibinfo {volume} {107}},\
  \bibinfo {pages} {155501} (\bibinfo {year} {2011})}\BibitemShut {NoStop}%
\bibitem [{\citenamefont {Fejes}(1942)}]{fejes1942dichteste}%
  \BibitemOpen
  \bibfield  {author} {\bibinfo {author} {\bibfnamefont {L.}~\bibnamefont
  {Fejes}},\ }\bibfield  {title} {\bibinfo {title} {{\"U}ber die dichteste
  kugellagerung},\ }\href@noop {} {\bibfield  {journal} {\bibinfo  {journal}
  {Mathematische Zeitschrift}\ }\textbf {\bibinfo {volume} {48}},\ \bibinfo
  {pages} {676} (\bibinfo {year} {1942})}\BibitemShut {NoStop}%
\bibitem [{\citenamefont {Kallus}\ and\ \citenamefont
  {Kusner}(2016)}]{kallus2016local}%
  \BibitemOpen
  \bibfield  {author} {\bibinfo {author} {\bibfnamefont {Y.}~\bibnamefont
  {Kallus}}\ and\ \bibinfo {author} {\bibfnamefont {W.}~\bibnamefont
  {Kusner}},\ }\bibfield  {title} {\bibinfo {title} {The local optimality of
  the double lattice packing},\ }\href@noop {} {\bibfield  {journal} {\bibinfo
  {journal} {Discrete \& computational geometry}\ }\textbf {\bibinfo {volume}
  {56}},\ \bibinfo {pages} {449} (\bibinfo {year} {2016})}\BibitemShut
  {NoStop}%
\bibitem [{\citenamefont {Kallus}(2015)}]{kallus2015pessimal}%
  \BibitemOpen
  \bibfield  {author} {\bibinfo {author} {\bibfnamefont {Y.}~\bibnamefont
  {Kallus}},\ }\bibfield  {title} {\bibinfo {title} {Pessimal packing shapes},\
  }\href@noop {} {\bibfield  {journal} {\bibinfo  {journal} {Geometry \&
  Topology}\ }\textbf {\bibinfo {volume} {19}},\ \bibinfo {pages} {343}
  (\bibinfo {year} {2015})}\BibitemShut {NoStop}%
\bibitem [{\citenamefont {Rogers}(1951)}]{rogers1951}%
  \BibitemOpen
  \bibfield  {author} {\bibinfo {author} {\bibfnamefont {C.~A.}\ \bibnamefont
  {Rogers}},\ }\bibfield  {title} {\bibinfo {title} {The closest packing of
  convex two-dimensional domains},\ }\href@noop {} {\bibfield  {journal}
  {\bibinfo  {journal} {Acta Mathematica}\ }\textbf {\bibinfo {volume} {86}},\
  \bibinfo {pages} {309} (\bibinfo {year} {1951})}\BibitemShut {NoStop}%
\bibitem [{\citenamefont {Mount}(1991)}]{mount1991}%
  \BibitemOpen
  \bibfield  {author} {\bibinfo {author} {\bibfnamefont {D.~M.}\ \bibnamefont
  {Mount}},\ }\bibfield  {title} {\bibinfo {title} {The densest double-lattice
  packing of a convex polygon},\ }\href@noop {} {\bibfield  {journal} {\bibinfo
   {journal} {Discrete and Computational Geometry: Papers from the DIMACS
  Special Year}\ }\textbf {\bibinfo {volume} {6}},\ \bibinfo {pages} {245}
  (\bibinfo {year} {1991})}\BibitemShut {NoStop}%
\bibitem [{\citenamefont {Duparcmeur}\ \emph {et~al.}(1995)\citenamefont
  {Duparcmeur}, \citenamefont {Gervois},\ and\ \citenamefont
  {Troadec}}]{duparcmeur1995dense}%
  \BibitemOpen
  \bibfield  {author} {\bibinfo {author} {\bibfnamefont {Y.~L.}\ \bibnamefont
  {Duparcmeur}}, \bibinfo {author} {\bibfnamefont {A.}~\bibnamefont
  {Gervois}},\ and\ \bibinfo {author} {\bibfnamefont {J.}~\bibnamefont
  {Troadec}},\ }\bibfield  {title} {\bibinfo {title} {Dense periodic packings
  of regular polygons},\ }\href@noop {} {\bibfield  {journal} {\bibinfo
  {journal} {Journal de Physique I}\ }\textbf {\bibinfo {volume} {5}},\
  \bibinfo {pages} {1539} (\bibinfo {year} {1995})}\BibitemShut {NoStop}%
\bibitem [{\citenamefont {Kurlin}(2023)}]{kurlin2023mathematics}%
  \BibitemOpen
  \bibfield  {author} {\bibinfo {author} {\bibfnamefont {V.~A.}\ \bibnamefont
  {Kurlin}},\ }\bibfield  {title} {\bibinfo {title} {Mathematics of
  2-dimensional lattices},\ }\href@noop {} {\bibfield  {journal} {\bibinfo
  {journal} {Foundations of Computational Mathematics}\ } (\bibinfo {year}
  {2023})}\BibitemShut {NoStop}%
\bibitem [{\citenamefont {Bright}\ \emph {et~al.}(2023)\citenamefont {Bright},
  \citenamefont {Cooper},\ and\ \citenamefont {Kurlin}}]{bright2023geographic}%
  \BibitemOpen
  \bibfield  {author} {\bibinfo {author} {\bibfnamefont {M.~J.}\ \bibnamefont
  {Bright}}, \bibinfo {author} {\bibfnamefont {A.~I.}\ \bibnamefont {Cooper}},\
  and\ \bibinfo {author} {\bibfnamefont {V.~A.}\ \bibnamefont {Kurlin}},\
  }\bibfield  {title} {\bibinfo {title} {Geographic-style maps for
  2-dimensional lattices},\ }\href@noop {} {\bibfield  {journal} {\bibinfo
  {journal} {Acta Crystallographica Section A}\ } (\bibinfo {year}
  {2023})}\BibitemShut {NoStop}%
\end{thebibliography}%
	
\end{document}